\documentclass[12pt]{article}
\author{Constantin N. Beli}
\title{Relative integral spinor norm groups over dyadic local fields}
\usepackage{amsfonts,amsmath,amssymb}

\def\D{\Delta}

\def\({\overline} \def\){\underline}
\def\<{\cdot} \def\go{\mathfrak}
\def\>{~~~~~~~}
\def\be{\begin{equation}}
\def\ee{\end{equation}}

 \def\ti{\times} 
\def\oo{{\cal O}} \def\vol{\go v} \def\sss{{\go s}} \def\nnn{{\go n}}
\def\fff{\dot{F}} \def\ooo{{\oo^\ti}} \def\oos{\oo^{\ti 2}}
\def\mo{{\rm mod}~}  
  \def\ffs{\fff^2}
 \def\p{\go p} \def\*{\sharp}
   \def\0{}
 \def\1{^{-1}}  
  \def\[{\prec} \def\]{\succ}
\def\bm{\left(\begin{array}} \def\em{\end{array}\right)}
   
 \def\N{{\rm N}} 
  \def\rep{{\rightarrow\!\!\!\! -}}

 \def\m2{~(\mo 2)} \def\no{\noindent}
 \def\btm{\begin{thm}}
\def\etm{\end{tm}}
 \def\blem{\begin{lem}}
\def\elem{\end{lem}}

\newtheorem{theorem}{Theorem}[section]
\newtheorem{proposition}[theorem]{Proposition}
\newtheorem{lemma}[theorem]{Lemma}
\newtheorem{definition}{Definition}
\newtheorem{corollary}[theorem]{Corollary}

\newtheorem{bof}[theorem]{}
\newtheorem{teorema}{Theorem}

\def\qed{\mbox{$\Box$}\vspace{\baselineskip}}
\def\pf{$Proof.~$} 
\def\bco{\begin{corollary}} \def\eco{\end{corollary}}
\def\bpr{\begin{proposition}} \def\epr{\end{proposition}}
\def\bdf{\begin{definition}} \def\edf{\end{definition}}
\def\btm{\begin{theorem}} \def\etm{\end{theorem}}
\def\blm{\begin{lemma}} \def\elm{\end{lemma}}
\def\bff{\begin{bof}\rm} \def\eff{\end{bof}}
\def\btr{\begin{teorema}} \def\etr{\end{teorema}}

\def\de{\newcommand} \de\tm[1]{{\no\bf Theorem~#1}} 
 \def\mb{\mathbb} 
\def\RR{{\mb R}} \def\QQ{{\mb Q}}  \def\ZZ{{\mb Z}}

\de\lm[1]{{\no\bf Lemma~#1}}
\de\df[1]{{\no\bf Definition~#1}} \de\co[1]{{\no\bf Corollary~#1}}
\de\tp[1]{\te (#1 )} \de\ts[1]{\te (O^-(#1 ))} \de\ty[1]{\te
(O(#1 ))} \de\tx[1]{\te (#1 )}
 \de\upn[2]{(1+\p^{#1})\ffs\cap\N (#2 )} \de\xt[2]{\te (#1 /#2 )}
 \de\ups[1]{(1+\p^{#1})\ffs} \de\upo[1]{(1+\p^{#1} )\ooo^2}
\de\upon[2]{(1+\p^{#1})\ooo^2\cap\N (#2 )}
\de\lr[1]{\longrightarrow^{\!\!\!\!\!\!\!\! #1}}
\de\lf[1]{\longleftarrow^{\!\!\!\!\!\!\!\! #1}}
\DeclareMathOperator\ord{ord} \DeclareMathOperator\rank{rank}
\DeclareMathOperator\Ima{Im}

\begin{document}
\maketitle

\begin{abstract}
We give explicit formulas for the relative spinor norm groups in the
case when the base field is dyadic.
\end{abstract}

\section{Introduction and background to BONG theory}

If $M,N$, with $N\subseteq M$, are two quadratic lattices over a
non-archimedian field of characteristic zero, then the set $X(M/N)$ of
the relative integral rotations, defined as
$$X(M/N)=\{ \sigma\in O^+(FM)\mid N\subseteq\sigma (M)\},$$
was introduced by Schulze-Pillot [SP] and Hsia [H2] as a
generalization of the group of usual (absolute) rotations of a
quadratic lattice $L$, $O^+(L)$. Calculating the spinor norms of
$O^+(L)$ and $X(M/N)$ is an important part of the theory of spinor
genera. More precisely, $\theta (O^+(L))$ is used for classification
purposes, while $\theta (X(M/N))$ is used for representation purposes.

The fact that $\theta (X(M/N))$ is a group was proved in [HSX]. In the
same paper explicit formulas for $\theta (X(M/N))$ were provided in
the nondyadic case. The 2-adic case (unramified finite extensions of
$\QQ_2$) was relegated to a separate source, with the mention that
``several pages are required merely to state the resulting
calculations''. Formulas for $\theta (X(M/N)$ in the 2-adic case are
provided in Shao's PhD thesis [S], but again, they are very
complicated, in both statements and proofs.

The common procedure used in dealing with $\theta (X(M/N))$, which
appears in both [HSX], for the purpose of proving that
$\theta (X(M/N))$ is a group, and in [X4, \S2], consists of two kinds
of reduction steps. The first kind applies when the norms and the
scales of $M$ and $N$ are not the same, i.e. when we have the strict
inclusions $\nnn M\supset\nnn N$ and/or $\sss M\supset\sss N$. Then
one proves that either $\theta (X(M/N))=\fff$ or there is a lattice
$M'$, $N\subseteq M'\subset M$, with $\nnn M'\subset\nnn M$ and/or
$\sss M'\subset\sss M$ such that $X(M/N)=X(M'/N)O^+(M)$. This has the
effect of decreasing the norm and/or the scale of $M$. After applying
repeatedly this procedure, we reduce to the case when $\nnn M=\nnn N$
and/or $\sss M=\sss N$. Then the second kind of reduction applies,
when the ranks of $M$ and $N$ are decreased. Namely, one proves that
either $\theta (X(M/N))=\fff$ or $X(M/N)=X(M^*/N^*)O^+(M)$ for some
lattices $M^*$ and $N^*$ of smaller ranks. In most cases we have
$M=K\perp M^*$ and $N=K\perp N^*$, where $K$ is modular of rank
$\leq 2$, with $\nnn K=\nnn M=\nnn N$ and $\sss K=\sss M=\sss N$. The
one exception is in case (ii) of the proof of [HSX, Theorem 2,1], when
$\nnn N=\nnn M=\sss N\subset\sss M$. In this case $N=\oo x\perp N^*$,
where $x\in N$ such that $Q(x)\oo =\nnn N$, but we don't have a
similar splitting $M=\oo x\perp M^*$. Instead, $M^*$ is defined by
$M^*=x^\perp\cap M$.

Our approach is similar, but instead of Jordan splittings we use good
BONGs. The first kind of reduction applies when $\nnn M\supset\nnn N$
and involves decreasing the norm of $M$. Namely, if
$M'=\{ x\in M\,\mid\, Q(x)\in\p\nnn M\}$ and $M'$ is a lattice, then
$N\subseteq M'\subset M$ and $\nnn M'\subset M$ and we have
$X(M/N)=X(M'/N)$. If $M'$ is not a lattice, then
$\theta (X(M/N))=\fff$. The second kind of reduction applies when
$\nnn M=\nnn N$. Then we prove that either $\theta (X(M/N))=\fff$ or
every first element of any good BONG of $N$ is also a first element in
a good BONG of $M$ and for every first element $x$ of a good BONG of
$N$ we have $X(M/N)=X(M^*/N^*)O^+(M)$, where $M^*=pr_{X^\perp}M$ and
$N^*=pr_{X^\perp}N$.
\bigskip

We now give a brief summary of the notations and results we use,
including a review of the BONG theory.

Throughtout the paper we use mostly the notation from [OM]. We have a
dyadic field $F$ and we denote by $\oo$ the ring of integers, by $\p$
the prime ideal, by $\pi$ a generator of $\p$ by
$\ord :F\to\ZZ\cup\{\infty\}$ the valuation map and by $e$ the
ramification index $e(F/\QQ_2)$, i.e. $e=\ord 2$.

We denote by $(\cdot,\cdot )_\p :\fff/\ffs\times\fff/\ffs\to\{\pm 1\}$
the Hilbert symbol. If $c\in\fff$ we denote by $\N (c)$ the norm group
$\N (F(\sqrt c)/F)=\{ x\in\fff\,\mid\, (x,c)_\p =1\}$.

The quadratic defect of an element $a\in F$ is the ideal
${\mathfrak d}(a)=\bigcap_{x\in F}(a-x^2)\oo$. We denote by
$\Delta =1-4\rho$, with $\rho\in\ooo$, a unit of quadratic defect
$4\oo$.

\bff We defined the order of the relative quadratic defect as the map
$$d:\fff/\ffs :\ZZ\cup\{\infty\},\qquad
d(a)=\ord a^{-1}{\mathfrak d}(a).$$ 
It has the following properties:

(1) $\Ima d=\{ 0,1,3,5,\ldots,2e-1,2e,\infty\}$.

(2) $d(a)=\infty$ iff $a=1$. $d(a)=2e$ iff $a=\Delta$. $d(a)=0$ iff
$\ord a$ is odd.

(3) $d$ satisfies the domination principle
$d(ab)\geq\min\{ d(a),d(b)\}$ $\forall a,b\in\fff$.

(4) If $d(a)+d(b)>2e$, then $(a,b)_\p =1$.

(5) if $d(a)<\infty$, then there is $b\in\fff$ with $d(b)=2e-d(a)$
such that $(a,b)_\p =-1$.

(6) $d(-1)\geq e$. (We have
${\mathfrak d}(-1)\subseteq (-1-1^2)\oo =2\oo$, so
$d(-1)=\ord{\mathfrak d}(-1)\geq\ord 2=e$.)
\eff

If $k$ is a positive integer, then
$\ups k=\{ a\in\fff\,\mid\, d(a)\geq k\}$. By convention, we set
$\ups\alpha =\{ a\in\fff\,\mid\, d(a)\geq\alpha\}$ for every
$\alpha\in\RR\cup\{\infty\}$. Then $\ups\alpha =\ups k$, where $k$ is
the smallest element in $\Ima d$ (possibly, $k=\infty$) such that
$k\geq\alpha$. Also $\ups\alpha =\ups{\lceil\alpha\rceil}$
$\forall\alpha\in\RR$.

In particular, $\ups\alpha =\fff$ if $\alpha\leq 0$, $\ups 1=\ooo\ffs$,
$\ups{2e}=\langle\Delta\rangle\ffs$ and $\ups\alpha =\ffs$ if
$\alpha >2e$.

\blm Let $a,c\in\fff$ and let $\alpha\in\RR\cup\{\infty\}$.

(i) We have $\ups\alpha\subseteq\N (c)$ iff $\alpha +d(c)>2e$.

(ii) We have $\upn\alpha a\subseteq\N (c)$ iff
$\alpha +\max\{ d(ac),d(c)\} >2e$.
\elm
\pf (i) is just [B2, Lemma 1.2(iii)].

(ii) By [B2, Lemma 1.3(i)], $\upn\alpha a\subseteq\N (c)$ iff
$\ups\alpha\subseteq\N (ac)$ or $\ups\alpha\subseteq\N (c)$, which, by
(i), is equivalent to $\alpha +d(ac)>2e$ or $\alpha +d(c)>2e$, i.e. to
$\alpha +\max\{ d(ac),d(c)\} >2e$. \qed

If $V$ is a quadratic space, its corresponding quadratic form and the
resulting bilinear form will be denoted by $Q$ and $B$. If
$x_1,\ldots,x_n$ is an orthogonal basis for $V$, with $Q(x_i)=a_i$,
then we say that $V\cong [a_1,\ldots,a_n]$ relative to the basis
$x_1,\ldots,x_n$. 

If $L$ is a quadratic lattice, we denote by $\sss L$, $\nnn L$ and
$\vol L$ its scale, norm and volume.

For convenience, if $c_1,c_2,\ldots$ is a sequence in $\fff$, we
denote $c_{i,j}:=c_ic_{i+1}\cdots c_j$ when $i\leq j$ and we also put
$c_{i,i-1}:=1$. 

An element $x$ of a quadratic lattice $L$ is called a norm generator
if $Q(x)\oo =\nnn L$. Then a basis of norm generators (BONG) of $L$ is
defined recursively as a sequence $x_1,\ldots,x_n$ such that $x_1$ is
a norm generator of $L$ and $x_2,\ldots,x_n$ is a BONG of
$pr_{x_1^\perp}L$. A lattice is uniquely determined by a BONG so we
will write $L={\prec x_1,\ldots,x_n\succ}$ if $x_1,\ldots,x_n$ is a BONG
of $L$. Moreover, if $Q(x_i)=a_i$, then we say that
$L\cong{\prec a_1,\ldots,a_n\succ}$ relative to the BONG $x_1,\ldots,x_n$.

\blm If $L\cong{\prec a_1,\ldots,a_n\succ}$ relative to some BONG and
$\ord a_i=R_i$, then

(i) In $\fff/\oos$ we have $\det L=a_{1,n}$

(ii) $\ord\vol L=R_1+\cdots +R_n$.
\elm
\pf (i) is just [B2, Lemma 2.1]. And for (ii) we take orders in the
relation $\vol L=(\det L)\oo =a_{1,n}\oo$. \qed

If $L$ is binary, we define $a(L)\in\fff/\oos$ by $a(L)=a^{-2}\det L$,
where $a\in\fff$ with $a\oo =\nnn L$ and $R(L)\in\ZZ$ by
$R(L)=\ord a(L)=\ord (\nnn L)^{-2}\vol L$. If
$L\cong{\prec a_1,a_2\succ}$, with $\ord a_i=R_i$, then $a(L)=a_2/a_1$
and $R(L)=R_2-R_1$. The invariant $a(L)$ uniquely determine the class
of $L$ up to scaling.

The set of all possible values of $a(L)$ for a binary lattice $L$ is
denoted by ${\mathcal A}={\mathcal A}_F$. We have
${\mathcal A}=\{ a\in\fff/\ooo\,\mid\, 4a\in\oo,\,
{\mathfrak d}(-a)\in\oo\}$. If $R=\ord a$, then $a\in\mathcal A$ iff
$R+2e\geq 0$ and $R+d(-a)\geq 0$. The equality $R+2e=0$ holds iff
$a=-1/4$ or $-\Delta/4$ (in $\fff/\oos$). The equality $R+d(-a)\geq 0$
holds iff $a=-\Delta/4$. Also $R$ cannot be odd and negative. (See
[B2, Lemma 3.5].)

As seen in [B2, \S3], we have two maps
$g:{\mathcal A}\to Sgp(\ooo/\oos )$  and 
$G:{\mathcal A}\to Sgp(\fff/\ffs )$ such that for every binary lattice
$L$ we have $\theta (O^+(L))=G(a(L))$ (see [B2, Lemma 3,7]) and if
$\eta\in\ooo$, then $L\cong L^\eta$ iff $\eta\in g(a(L))$ (see [B2,
Lemma 3.11 and 3.12]). If $L\cong\prec a_1,a_2\succ$, these write as
$\theta (O^+(L))=G(a_2/a_1)$ and $L\cong\prec\eta a_1,\eta a_2\succ$
iff $\eta\in g(a_2/a_1)$. 

Here $Sgp(\ooo/\oos )$ and $Sgp(\fff/\ffs )$ denote the sets of all
subgroups of $\ooo/\oos$ and $\fff/\ffs$. They may be identified with
subgroups of $\ooo$ containing $\oos$ and subgroups of $\fff$
containing $\ffs$, respectively.

Note that the formulas for $\theta (O^+(L))$ in the binary case, which
are expressed by $\theta (O^+(L))=G(a_2/a_1)$ were obtained in the
case when $L$ is modular by Hsia and Xu (see
[H1, Propositions B, C and D] and [X1]) and in the case when $L$ is not
modular by Xu (see [X2, Theorems 2.1, 2.2 and 2.3]).

A BONG $x_1,\ldots,x_n$ of a lattice $L$ is called good if
$\ord Q(x_i)\leq\ord Q(x_{i+2})$ for $1\leq i\leq n-2$. By [B2, Lemmas
4.3(ii), 3.6 and 3.5], we have the following lemma, which gives
necessary and sufficient conditions for an orthogonal basis of a
quadratic space to be a good BONG of a lattice.

\blm If $x_1,\ldots,x_n$ is an orthogonal basis of a quadratic space
$V$, $Q(x_i)=a_i$ and $\ord a_i=R_i$, then there is a lattice
$L\cong{\prec a_1,\ldots,a_n\succ}$ relative to the good BONG
$x_1,\ldots,x_n$ iff $R_i\leq R_{i+2}$ for $1\leq i\leq n-2$ and 
$a_{i+1}/a_i\in{\mathcal A}$ for $1\leq i\leq n-1$.

The second condition writes equivalently as
$$R_{i+1}-R_i+2e\geq 0\text{ and }R_{i+1}-R_i+d(-a_{i,i+1})\geq 0,\,
\forall 1\leq i\leq n-1.$$

Moreover, $R_{i+1}-R_i+2e=0$ iff $a_{i+1}/a_i=-1/4$ or $-\Delta/4$ in
$\fff/\oos$ and we have $R_{i+1}-R_i+d(-a_{i,i+1})=0$ iff
$a_{i+1}/a_i=-\Delta/4$ in $\fff/\oos$. And we have $a_{i+1}/a_i=-1/4$
or $-\Delta/4$ iff
${\prec a_i,a_{i+1}\succ}\cong\frac 12\pi^{R_1}A(0,0)$ or
$\frac 12\pi^{R_1}A(2,2\rho )$, respectively.

Also $R_{i+1}-R_i$ cannot be odd and negative.
\elm

If $L\cong{\prec a_1,\ldots,a_n\succ}$ relative to some good BONG, then
the orders of $a_1,\ldots,a_n$ are independent of the choice of the
good BONG. So we have the invariants $R_i(L)$, $R_i(L)=\ord a_i$, for
$1\leq i\leq n$.

For the rest of the section $L\cong{\prec a_1,\ldots,a_n\succ}$ relative
to a good BONG $x_1,\ldots,x_n$, $R_i(L)=R_i$ and
$a_i=\pi^{R_i}\varepsilon_i$, with $\varepsilon_i\in\ooo$.

\blm (i) $\nnn L=\p^{R_1}$ and $\sss L=\p^{(R_1+R_2)/2}$.

(ii) For any $1\leq i\leq n-1$ we have
$L={\prec x_1,\ldots,x_i\succ}\perp{\prec x_{i+1},\ldots,x_\succ}$ iff
$R_i\leq R_{i+1}$.

Conversely, if $x_1,\ldots,x_i$ and $x_{i+1},\ldots,x_n$ are good
BONGs for some lattices $L'$ and $L''$, then $x_1,\ldots,x_n$ is a
good BONG for $L:=L'\perp L''$ iff $R_i\leq R_{i+1}$,
$R_{i-1}\leq R_{i+1}$ (if $i\geq 2$) and $R_i\leq R_{i+2}$
(if $i\leq n-2$).

(iii) $L^\#\cong{\prec a_n^{-1},\ldots,a_1^{-1}\succ}$ relative to the
good BONG $x_n^\#,\ldots,x_1^\#$, where
$x^\# :=Q(x)^{-1}x$. Consequently,
$R_i(L^\# )=\ord a_{n+1-i}^{-1}=-R_{n+1-i}$. 

(iv) If $1\leq i\leq j\leq n$ and
${\prec x_i,\ldots,x_j\succ}\cong{\prec b_i,\ldots,b_j\succ}$ 
relative to a good BONG $y_i,\ldots,y_j$, then
$L\cong{\prec a_1,\ldots,a_{i-1},b_i,\ldots,b_j,a_{j+1},\ldots,a_n\succ}$
relative to the good BONG
$x_1,\ldots,x_{i-1},y_i,\ldots,y_j,x_{j+1},\ldots,x_n$. In particular,
if $1\leq i\leq n-1$ and $\eta\in g(a_{i+1}/a_i)$, then
$L\cong{\prec a_1,\ldots,a_{i-1},\eta a_i,\eta _{i+1},
a_{i+2},\ldots,a_n\succ}$ relative to the good BONG.

(v) $O({\prec x_i,\ldots,x_j\succ})\subseteq O(L)$ for every
$1\leq i\leq j\leq n$. Consequently, $\tau_{x_i}\in O^-(L)$, so
$a_i\in\theta (O^-(L))$, for all $1\leq i\leq n$ and
$G(a_{i+1}/a_i)\subseteq\theta (O^+(L))$ for all $1\leq i\leq n-1$.
\elm
\pf The main statements of (i)-(v) are [B2, Corollary 4.4(iv), (i) and
(v), Lemma 4.8, Lemma 4.9(ii) and (i)]. The second statement of (iv)
follows by taking $j=i+1$ in the first one, noting that
${\prec a_i,a_{i+1}\succ}\cong{\prec\eta a_i,\eta a_{i+1}\succ}$ for
all $\eta\in g(a_{i+1}/a_i)$. And for the second and third statement of
(v) we note that
$\tau_{x_i}\in{O^-(\prec x_i\succ )}\subseteq O^-(L)$, so
$a_i=\theta (\tau_{x_i})\in\theta (O^-(L))$ and
$G(a_{i+1}/a_i)={\theta (O^+(\prec x_i,x_{i+1}\succ ))}
\subseteq\theta (O^+(L))$. \qed

We say that $L$ has property A if $R_i<R_{i+2}$ for
$1\leq i\leq n-2$. If $L$ has property A, we say that it has property
B if for $1\leq i\leq n-1$, in each of the following two cases.
$$R_{i+1}-R_i\text{ is odd and }\leq 2e+1,$$
$$R_{i+1}-R_i\text{ is even and }d(-a_{i,i+1})\leq e-(R_{i+1}-R_i)/2$$
we have $R_i-R_{i-1}>2e$ (if $i\geq 2$) and $R_{i+2}-R_{i+1}>2e$ (if
$i\leq n-2$).

\bpr (i) If $L$ has property A then every BONG of $L$ is good.

(ii) If $L$ has property A, but not property B, then
$\theta (O^+(L))=\fff$

(iii) If $L$ doesn't have property A, then
$\theta (O^+(L))\supseteq\ooo\ffs$, i.e. $\theta (O^+(L))=\ooo\ffs$ or
$\fff$.
\epr
\pf For (i) see [B2, Corollary 4.2(i)].  For (ii) see [B2, Lemma 4.11].

Unfortunately, the proof of (iii), which is stated at the beginning of
[B2, \S7], is missing from [B2]. Although (iii) can be proved in a
 way similar to the proof of the inclusion
$\ups{[(R_{i+2}-R_i)/2]}\subseteq\theta (O^+(L))$ from [B2, \S5], for
convenience, we will improvise a short proof which uses some results
from \S4 and \S6. We take first the case $n=3$. The fact that
condition A fails means that $R_1=R_3$. By Lemma 4.4 there is a
lattice $M\cong{\prec\pi^{-2}a_1,a_2,a_3\succ}$ such that $M'=L$. By
Lemma 6.2, this implies that $O^+(M)\subseteq O^+(L)$. If
$R_i(M)=S_i$, then $S_1=\ord\pi^{-2}a_1=R_1-2$, $S_2=R_2$ and
$S_3=R_3=R_1$. Then $S_1<S_3$, so $M$ has property A. By
[B2, Theorem 1] (i.e. Theorem 1.7 below), we get
$\theta (O^+(L))\supseteq\theta (O^+(M))\supseteq
\ups{[(S_3-S_1)/2]}=\ups 1=\ooo\ffs$.

For the general case, if $L$ doesn't have property A, then
$R_i=R_{i+2}$ for some $i$. Hence, if
$K\cong{\prec a_i,a_{i+1},a_{i+2}\succ}$, then, by the case $n=3$ and
Lemma 1.5(v), we have
$\ooo\ffs\subseteq\theta (O^+(K))\subseteq\theta (O^+(L))$. \qed 

{\bf Remark.} In [B2, Definition 7], property A is given in terms of
Jordan decompositions. The condition that $R_i<R_{i+2}$ for all
$1\leq i\leq n-2$ is just an equivalent statement, by
[B2, Lemma 4.3(i)]. Property A, with the original definition, appears
in the proof of [HSX, Theorem 2.1]. There a weaker statement is
proved. Namely, that if property A fails, then
$[\fff :\theta (O^+(L))]\leq 2$. In [B1] the stronger conclusion, that
$\theta (O^+(L)\supseteq\ooo\ffs$, was obtained by similar arguments.
Namely, we prove that if $L$ doesn't have property A, then either $L$
splits a lattice of the form $H\cong\pi^rA(0,0)$, in which case
$\ooo\ffs =\theta (O^+(H))\subseteq\theta (O^+(L))$, or
$\theta (O^+(L))=\fff$. The proof uses some results by Hsia and Xu
from [H1] and [X3], including Xu's minimal norm splittings. (See the
beginning of [B1, \S I.1.4].)

We now state the main results of [B2], with formulas for
$\theta (O^+(L))$. For the case when $L$ has property A we have
[B2, Theorem 1]: 

\btm If $L$ has property $A$ then
$$\theta (O^+(L)=G(a_2/a_1)\cdots G(a_n/a_{n-1})\ups\alpha,$$
where $\alpha =\min\{[(R_{i+2}-R_i)/2]\,\mid\,1\leq i\leq n-2\}$.
\etm

\bff Note that $\ups\alpha$ is the product of the groups
$\ups{[(R_{i+2}-R_i)/2]}$ for $1\leq i\leq n-2$. By the Remark at the
end of [B2, \S5], if $R_{i+2}-R_i$ is odd, then
$\ups{[(R_{i+2}-R_i)/2]}\subseteq G(a_{i+1}/a_i)G(a_{i+2}/a_{i+1})$,
so the factor $\ups{[(R_{i+2}-R_i)/2]}$ may be removed from the formula
for $\theta (O^+(L))$. So in Theorem 1.7 $\alpha$ may be replaced by
$\alpha':=\min\{ (R_{i+2}-R_i)/2\,\mid\, R_{i+2}-R_i\text{ is even}\}$,
as stated in [B2, \S5].

In fact we have a more precise result. If the largest factor of
$\ups\alpha=\prod\ups{[(R_{i+2}-R_i)/2]}$ occurs at an index $i$ such
that $R_{i+2}-R_i$ is odd, then $\ups\alpha
=\ups{[(R_{i+2}-R_i)/2]}\subseteq G(a_{i+1}/a_i)G(a_{i+2}/a_{i+1})$, so
the factor $\ups\alpha$ can be removed from the formula for
$\theta (O^+(L))$.
\eff

If $L$ doesn't have property A, then, by Proposition 1.6(iii),
$\theta (O^+(L))$ can only be $\fff$ or $\ooo\ffs$. The next result,
which is [B2, Theorem 3], decides between the two possibilities.

\btm We have $\theta (O^+(L))\subseteq\ooo\ffs$ iff
$G(a_{i+1}/a_i)\subseteq\ooo\ffs$ for all $1\leq i\leq n-1$ and
$(R_{i+1}-R_i)/2+e$ is even for all $1\leq i\leq n-2$ such that
$R_i=R_{i+2}$.
\etm

{\bf Remark.} By [B2, Lemma 7.2], condition
$G(a_{i+1}/a_i)\subseteq\ooo\ffs$ implies that
$R_{i+1}-R_i=\ord a_{i+1}/a_i$ is even. Hence if
$\theta (O^(L))\subseteq\ooo\ffs$, then
$R_1\equiv\cdots\equiv R_n\pmod 2$
\medskip

For the proof of Theorems 1.7 and 1.9 in [B2] we used the following
key lemmas, which we will need here as well. See
[B2, Lemmas 6.6 and 7.3].

\blm If $L$ has property B, then for every $x\in L$ with
$Q(x)=Q(x_1)=a_1$ there is some $\sigma\in O^+(L)$ with
$\sigma (x_1)=x$. 
\elm

\blm If $R_i=R_{i+2}$ and
$G(a_{i+1}/a_i),G(a_{i+2}/a_{i+1})\subseteq\ooo\ffs$, then
${\prec x_1,x_2,x_3\succ} =H\perp\oo v$, where
$H\cong\pi^{(R_i+R_{i+1})/2}A(0,0)$ and
$Q(v)=-\pi^{R_i}\varepsilon_{1,3}$. Moreover, $M=H\perp K$, where
$K\cong{\prec a_1\ldots,a_{i-1},-\pi^{R_i}
\varepsilon_{i,i+2},a_{i+3},\ldots a_n\succ}$ relative to the good BONG
$x_1\ldots,x_{i-1},v,x_{i+3}\ldots x_n$. 
\elm

We now state the main result of [B4], the classification theorem for
dyadic quadratic lattices in terms of BONGs.

The invariants $\alpha_i(L)$, for $1\leq i\leq n-1$ are defined as
\begin{multline*}\alpha_i(L)=\min (\{(R_{i+1}-R_i)/2+e\}
\cup\{ R_{i+1}-R_j+d(-a_{j,j+1})\,\mid\, 1\leq j\leq i\}\\
\cup\{ R_{j+1}-R_i+d(-a_{j,j+1})\,\mid\, i\leq j\leq n-1\} ).
\end{multline*}

\btm Let $L\cong{\prec a_1,\ldots,a_n\succ}$ and
$K\cong{\prec b_1,\ldots,b_n\succ}$ relative to good BONGs such that
$FL\cong FK$. Let $R_i=R(L)$, $S_i=R(K)$, $\alpha_i=\alpha_i(L)$
and $\beta_i=\alpha_i(K)$. Then $L\cong K$ iff the following
conditions hold.

(i) $R_i=S_i$ for $1\leq i\leq n$.

(ii) $\alpha_i=\beta_i$ for $1\leq i\leq n-1$.

(iii) $d(a_{1,i}b_{1,i})\geq\alpha_i$ for $1\leq i\leq n-1$.
  
(iv) $[b_1,\ldots,b_{i-1}]\rep [a_1,\ldots,a_i]$ for every
$1<i<n$ such that $\alpha_{i-1}+\alpha_i>2e$.
\etm

\bco If $L\cong K$ and $R_{i+1}-R_i>2e$ for some $1\leq i\leq n-1$,
then $[a_1,\ldots,a_i]\cong [b_1,\ldots,b_i]$. In particular,
$a_{1,i}=b_{1,i}$.

Same thing happens if $i=0$ or $n$.
\eco
\pf By condition (iii) of Theorem 1.12 and [B4, Corollary 2.8(ii)], we
have $d(a_{1,i}b_{1,i})\geq\alpha_i>2e$. This implies
$d(a_{1,i}b_{1,i})=\infty$, i.e.  $a_{1,i}b_{1,i}=1$, so
$a_{1,i}=b_{1,i}$ in $\fff/\ffs$. If $i=1$, this means $a_1=b_1$, so
$[a_1]\cong [b_1]$ trivially. If $i>1$, since also
$\alpha_{i-1}\geq 0$, we have $\alpha_{i-1}+\alpha_i>2e$. By condition
(iv) of Theorem 1.12, we get
$[b_1,\ldots,b_{i-1}]\rep [a_1,\ldots,a_i]$, which, by [OM, 63:21], is
equivalent to
$[a_1,\ldots,a_i]\cong [b_1,\ldots,b_{i-1},a_{1,i}b_{1,i-1}]$. But in
$\fff/\ffs$ we have $a_{1,i}b_{1,i}=1$, so
$a_{1,i}b_{1,i-1}=b_i$. Hence $[a_1,\ldots,a_i]\cong [b_1,\ldots,b_i]$.

If $i=0$ then both $[a_1,\ldots,a_0]$ and $[b_1,\ldots,b_0]$ are the
trivial zero space and we have $a_{1,0}=b_{1,0}=1$. And if $i=n$, then 
$[a_1,\ldots,a_n]\cong [b_1,\ldots,b_n]\cong FL$  and
$a_{1n}=b_{1,n}=\det FL$. \qed

We also have the representation theorem announced in [B3, Theorem 4.5]
and proved in [B5]. However, since these results are not published, we
will not use them here. In the Appendix we will show how Theorem 7.3
can be improved if we assume the results from [B5]. 

\section{The maps $g$ and $G$ revisited}

In [B2, Definitions 4 and 6] we introduced the maps
$g:\mathcal A\to Sgp(\ooo/\oos )$ and
$G:\fff/\oos\to Sgp(\fff/\fff^2)$. Here 
$Sgp(\ooo/\oos )$ and $Sgp(\fff/\fff^2)$ denote the sets of all
subgroups of $\ooo/\oos$ and $\fff/\fff^2$, which we identify with
subgroups of $\ooo$ containing $\oos$ and subgroups of $\fff$
containing $\fff^2$, respectively.

Since $\ooo\cap\ffs=\oos$, we have an isomorphism
$\ooo\ffs/\ffs\cong\ooo/\oos$, which induces a bijection
$Sgp(\ooo\ffs/\ffs )\to Sgp(\ooo/\oos )$, given by
$H\mapsto H\cap\ooo$. Its inverse is given by $H\mapsto H\ffs$.

In this section we define some related functions, $\bar g$ and
$\bar G$, which help us write the main results in a more convenient
way and, in general, are more user friendly. We start by noting that
we have an injective map $\fff/\oos\to\fff/\fff^2\times\RR$ given by
$a\mapsto (a,\ord a)$. An element $(a,R)\in\fff/\fff^2\times\RR$
belongs to the image of this map iff $R\in\ZZ$ and
$\ord a\equiv R\pmod 2$. We extend the definition of $g$ and $G$ on
the whole $\fff/\fff^2\times\RR$.

\bdf We define
$$\bar{\mathcal A}=\{ (a,R)\in\fff/\ffs\times\RR\,\mid\,
R+2e\geq 0,\, R+d(-a)>0\}\cup\{ (-\Delta,-2e)\}.$$
\edf

\bff{\bf Remarks.} If $a\neq -1$, then $d(-a)<\infty$, i.e.
$d(-a)\leq 2e$, so $R+d(-a)>0$ implies that $R+2e\geq 0$. Thus the
condition $R+2e\geq 0$ from the definition of $\bar{\mathcal A}$ is
superfluous. So if $d(-a)<2e$, then $(a,R)\in\bar{\mathcal A}$ iff
$R+d(-a)>0$. If $d(-a)=2e$, i.e. if $a=-\Delta$, then
$(a,R)\in\bar{\mathcal A}$ iff either $R+d(-a)=R+2e>0$ or $R=-2e$ (the
case $(a,R)=(-\Delta,-2e)$.) Thus $(-\Delta,R)\in\bar{\mathcal A}$ iff
$R+2e\geq 0$. And if $d(-a)=\infty$, i.e. if $a=-1$, the condition
$R+d(-a)>0$ is superfluous, so $(a,R)\in\bar{\mathcal A}$ iff
$R+2e>0$. In conclusion:
\begin{multline*}
\bar{\mathcal A}=\{ (a,R)\in\fff/\ffs\times\RR\,\mid\, d(-a)<2e,\,
R+d(-a)>0\}\\
\cup\{ (a,R)\in\fff/\ffs\times\RR\,\mid\, d(-a)\geq 2e,\, R+2e\geq 0\}.
\end{multline*}
Also note that if $R+d(-a)>0$ and $a\neq -1$, then $d(-a)\leq 2e$, so
$R+2e\geq R+d(-a)>0$. Hence $\{ (a,R)\in\fff/\ffs\times\RR\,\mid\,
R+2e\geq 0,\, R+d(-a)>0\} =\{ (a,R)\in\fff/\ffs\times\RR\,\mid\,
R+2e>0,\, R+d(-a)>0\}\cup\{ (-1,-2e)\}$, which gives another
description of $\bar{\mathcal A}$,
$$\bar{\mathcal A}=\{ (a,R)\in\fff/\ffs\times\RR\,\mid\,
R+2e>0,\, R+d(-a)>0\}\cup\{ (-1,-2e),(-\Delta,-2e)\}.$$
\eff

\blm The preimage of $\bar{\mathcal A}$ under the embedding
$\fff/\oos\to\fff/\ffs\times\RR$ is $\mathcal A$.
\elm
\pf Let $a\in\fff/\oos$, with $\ord a=R$. Its image under the
embedding $\fff/\oos\to\fff/\ffs\times\RR$ is $(a,R)$. If
$a\in\mathcal A$ then $R+2e\geq 0$ and $R+d(-a)\geq 0$. Moreover, the
inequality $R+d(-a)\geq 0$ is strict unless $a=-\Delta/4$. Hence
either $R+2e\geq 0$ and $R+d(-a)>0$ or $a=-\Delta/4$. In the first
case $(a,R)\in\bar{\mathcal A}$, by definition. In the second case
$a=-\Delta$ in $\fff/\ffs$ and $R=\ord a=-2e$, so
$(a,R)=(-\Delta,-2e)\in\bar{\mathcal A}$.

Conversely, if $(a,R)\in\bar{\mathcal A}$ then either $R+2e\geq 0$ and
$R+d(-a)>0$ or $(a,R)=(-\Delta,-2e)$. In both cases we have
$R+2e\geq 0$ and $R+d(-a)\geq 0$ so, by definition,
$a\in\mathcal A$. \qed

\bdf We define
$\hat g,\, \bar g :\fff/\fff^2\times\RR\to Sgp(\fff/\fff^2)$
as follows:
$$\hat g(a,R)=\upn{\alpha (a,R)}{-a},\text{ where }
\alpha (a,R)=\min\{ R/2+e,R+d(-a)\},$$
$\bar g(a,R)=\ooo\fff^2$ if $(a,R)=(-1,-2e)$ and
$\bar g(a,R)=\hat g(a,R)$ otherwise.

Note that
$\hat g(-1,-2e)=\upn 01=\fff\supset\ooo\ffs=\bar g(-1,-1e)$, so
$\hat g(a,R)\supseteq\bar g(a,R)$ for all $(a,R)$, with equality if
$(a,R)\neq (-1,-2e)$.
\edf

\blm (i) We have $\bar g(a,R)\subseteq\ooo\ffs$ iff
$(a,R)\in\bar{\mathcal A}$ or $a=-\Delta$. Equivalently,
$\bar g(a,R)\subseteq\ooo\ffs$ iff either $a=-\Delta$ or $R+2e\geq 0$
and $R+d(-a)>0$.

(ii) If $(a,R)\notin\bar{\mathcal A}$, then $\bar g(a,R)={\rm N}(-a)$.
\elm
\pf Let $(a,R)\in\fff/\ffs\times\RR$. We use the second expression for
$\mathcal A$ from 2.1.

If $a=-\Delta$, then $\bar g(a,R)\subseteq\N (-a)=\ooo\ffs$, so we are
done. Suppose now that $(a,R)\in\mathcal A$. If $(a,R)=(-1,-2e)$, then
$\bar g(a,R)=\ooo\ffs$ by definition and if $(a,R)=(-\Delta,-2e)$, then
we are in the case $a=-\Delta$, already treated. So we are left with
the case when $R+2e>2e$, so $R/2+e>0$, and $R+d(-a)>0$. Then
$\alpha (a,R)>0$, which implies that
$\bar g(a,R)\subseteq\ups{\alpha (a,R)}\subseteq\ooo\ffs$. So we have
the ``if'' part of (i).

If $(a,R)\notin\mathcal A$, then $(a,R)\neq (-1,-2e)$, so
$\bar g(a,R)=\hat g(a,R)$. By  2.1, we have
$(a,R)\notin\{ (a,R)\in\fff/\ffs\times\RR\,\mid\, R+2e>0,\, R+d(-a)>0\}$.
Then either $R+2e\leq 0$, so $R/2+e\leq 0$, or $R+d(-a)\leq 0$. In
both cases, $\alpha (a,R)\leq 0$, so
$\bar g(a,R)=\upn{\alpha (a,R)}{-a}=\fff\cap\N (-a)=\N (-a)$, and so
we have (ii).

Finally, if $a\neq -\Delta$ and $(a,R)\notin\mathcal A$, then, by
(ii), we have $\bar g(a,R)=\N (-a)$. But $-a\neq\Delta$, so
$\N (-a)\not\subseteq\N (\Delta )=\ooo\ffs$. Hence we have the ``only
if'' part of (i).\qed 

\blm We have
$$\hat g(a,R)=\begin{cases}\ups{R/2+e}&\text{if }d(-a)>e-R/2\\
\upn{R+d(-a)}{-a}&\text{if }d(-a)\leq e-R/2\end{cases}.$$

If $d(-a)\leq e-R/2$, then $\bar g(a,R)=\hat g(a,R)$ and
$[\ups{R+d(-a)}:\bar g(a,R)]=2$.
\elm
\pf If $d(-a)>$ or $\leq e-R/2$, then $R+d(-a)>$ or $\leq R/2+e$, so
$\alpha (a,R)=R/2+e$ or $R+d(-a)$, accordingly.

If $d(-a)>e-R/2$ we get $\hat g(a,R)=\upn{R/2+e}{-a}$. But we also
have $d(-a)+(R/2+e)>2e$, so, by Lemma 1.2(i),
$\hat{R/2+e}\subseteq{\rm N}(-a)$. Hence $\hat g(a,R)=\ups{R/2+e}$.

If $d(-a)\leq e-R/2$, then $\hat g(a,R)=\upn{R+d(-a)}{-a}$. But
$d(-a)+(R+d(-a))\leq 2e$, so ${\rm N}(-a)\not\subseteq{R+d(-a)}$.
Hence we have a strict inclusion $\hat g(a,R)\subset\ups{R+d(-a)}$.
Then $\dfrac{\ups{R+d(-a)}}{\hat g(a,R)}$ is a nontrivial group that
embeds into $\dfrac\fff{{\rm N}(-a)}\cong\ZZ/2\ZZ$. It follows that
$[\ups{R+d(-a)}:\hat g(a,R)]=2$. Also from $d(-a)\leq e-R/2<\infty$
we get that $a\neq -1$, so $\bar g(a,R)=\hat g(a,R)$. \qed

\bpr If $a\in\mathcal A$ and $\ord a=R$, then
$g(a)=\hat g(a,R)\cap\ooo=\bar g(a,R)\cap\ooo$ 
\epr
\pf By [B2, Definition 6], we have $g(a)=\oos$ if $R>2e$ and
$$g(a)=\begin{cases}\upo{R/2+e}&\text{if }d(-a)>e-R/2\\
\upon{R+d(-a)}{-a}&\text{if }d(-a)\leq e-R/2\end{cases}$$
if $R\leq 2e$. But the formula from the case $R\leq 2e$ applies also
in the case $R>2e$. Indeed, if $R>2e$, then $d(-a)\geq 0>e-R/2$ so the
formula $g(a)=\upo{R/2+e}$ applies. But $R/2+e>2e$, so
$\upo{R/2+e}=\oos$. Then $g(a)=\hat g(a,R)\cap\ooo$ follows from Lemma
. and the fact that $\upo\alpha =\ups\alpha\cap\ooo$.

If $(a,R)\neq (-1,-2e)$, then
$\hat g(a,R)\cap\ooo=\bar g(a,R)\cap\ooo$ follows from
$\hat g(a,R)=\bar g(a,R)$. And if
$(a,R)=(-1,-2e)$, then $\hat g(a,R)=\fff$ and $\bar g(a,R)=\ooo\ffs$,
so $\hat g(a,R)\cap\ooo=\bar g(a,R)\cap\ooo =\ooo$. \qed

\blm (i) For every $a\in\fff/\ffs$ the map $\bar g(a,\cdot )$ is
decreasing, i.e. if $R\leq S$, then $\bar g(a,R)\supseteq\bar g(a,S)$.

(ii) If $d(-a)\leq -R$, then $\bar g(a,R)=\N (-a)$.

(iii) If $R\leq 0$, then $a\in\bar g(a,R)$.

(iv) $\bar g(a,R)\subseteq\ups R$. In particular, if $R>2e$, then
$\bar g(a,R)=\ffs$.
\elm
\pf (i) If $R\leq S$, then $\alpha (a,R)\leq\alpha (a,S)$, so
$\upn{\alpha (a,R)}{-a}\supseteq\upn{\alpha (a,S)}{-a}$. i.e.
$\hat g(a,R)\supseteq\hat g(a,S)$. Hence $\hat g(a,\cdot )$ is
decreasing.

Since $\bar g$ and $\hat g$ coincide everywhere but at $(-1,-2e)$,
this implies that $\bar g(a,\cdot )=\hat g(a,\cdot )$ is decreasing
for $a\neq -1$ and for $a=-1$ we only need to prove that
$\bar g(-1,R)\supseteq\bar g(-1,-2e)$ when $R<-2e$ and
$\bar g(-1,-2e)\supseteq\bar g(-1,S)$ when $S>-2e$. We have
$\bar g(-1,-2e)=\ooo\ffs$. If $R<-2e$ then $(-1,R)\notin\mathcal A$,
so, by Lemma 2.3(ii), $\bar g(-1,R)=\N (1)=\fff\supset\ooo\ffs$. And
if $S>-2e$ then $(-1,S)\in\mathcal A$, so, by Lemma 2.3(i),
$\bar g(-1,S)\subseteq\ooo\ffs$. This concludes the proof.

(ii) If $d(-a)\leq -R$, then $\alpha{\alpha (a,R)}\leq R+d(-a)\leq 0$,
so $\hat g(a,R)=\upn{\alpha (a,R)}{-a}=\fff\cap\N (-a)=\N (-a)$. Since
$d(-a)\leq -R<\infty$, we have $a\neq -1$, so
$\bar g(a,R)=\hat g(a,R)$.

If $R>2e$, then $R/2+e>2e$ and $R+d(-a)>2e$, so $\alpha (a,R)>2e$. It
follows that $\bar g(a,R)=\hat g(a,R)\subseteq\ups\alpha =\ffs$. Hence
$\bar g(a,R)=\ffs$.

(iii) We have $-1\in\ooo\ffs=\bar g(-1,-2e)$. Suppose that
$(a,R)\neq (-1,-2e)$. Since $R\leq 0$, we have
$d(-1)\geq e\geq R/2+e\geq\alpha (a,R)$ and
$d(-a)\geq R+d(-a)\geq\alpha (a,R)$, so, by the domination principle,
$d(a)\geq\alpha (a,R)$, i.e. $a\in\ups{\alpha (a,R)}$. We also have
$(a,-a)_\p =1$, so $a\in\N (-a)$. Hence
$a\in\upn{\alpha (a,R)}{-a}=\hat g(a,R)=\bar g(a,R)$.

(iv) If $R>2e$, then $\ups R=\ffs$, but also $R/2+e>2e$ and
$R+d(-a)\geq R>2e$, so $\alpha (a,R)>2e$, which implies
$\bar g(a,R)\subseteq\ups{\alpha (a,R)}=\ffs$, i.e.
$\bar g(a,R)=\ffs$.

If $R\leq 2e$, then $R/2+e\geq R$ and $R+d(-a)\geq R$, so
$\alpha (a,R)\geq R$. It follows that
$\bar g(a,R)\subseteq\ups{\alpha (a,R)}\subseteq\ups R$. \qed

\blm If $a,b\in\fff/\ffs$ and $x\in\RR$ then the following are
equivalent.

(i) $\min\{ x,d(ab)\} +\max\{ d(a),d(b)\} >2x$.

(ii) $d(a)+d(b)>2x$.
\elm
\pf We may assume that, say, $d(a)\leq d(b)$. Then (i) writes as
$\min\{ x,d(ab)\} +d(b)>2x$, which is equivalent to $x+d(b)>2x$,
i.e. $d(b)>x$, and $d(ab)+d(b)>2x$.

Suppose first that $d(b)\leq x$. Then (i) fails and
$d(a)+d(b)\leq 2d(b)\leq 2x$, so (ii) fails, as well. Hence
(i)$\iff$(ii) holds.

Suppose now that $d(b)>x$. Then (i) is equivalent to $d(ab)+d(b)>2x$,
i.e. to $d(ab)>2x-d(b)$ and (ii) is equivalent to $d(a)>2x-d(b)$. But
$2d(b)>2x$, so $d(b)>2x-d(b)$. It follows that
$d(ab)>2x-d(b)\iff d(a)>2x-d(b)$, which concludes the proof. (By the
domination principle, we have $d(ab)\geq\min\{ d(a),d(b)\}$ and
$d(a)\geq\min\{ d(ab),d(b)\}$.) \qed 

\blm If $a,c\in\fff/\ffs$ and $R\in\RR$, we have
$\hat g(a,R)\subseteq\N (c)$ iff $R+d(-ac)+d(c)>2e$ and 
$\bar g(a,R)\subseteq\N (c)$ iff $R+d(-ac)+d(c)>2e$ or
$(a,R)=(-1,-2e)$ and $c=\Delta$.
\elm
\pf By Lemma 1.2(ii), we have
$\hat g(a,R)=\upn{\alpha (a,R)}{-a}\subseteq\N (c)$ iff
$\alpha (a,R)+\max\{ d(-ac),d(c)\}>2e$. Since
$\alpha (a,R)=\min\{ R/2+e,R+d(-a)\} =R+\min\{ e-R/2,d(-a)\}$, this
writes equivalently as
$\min\{ e-R/2,d(-a)\} +\max\{ d(-ac),d(c)\}>2e-R=2(e-R/2)$. Since
in $\fff/\ffs$ we have $(-ac)\cdot c=-a$, by Lemma 2.7, applied when
$a,b,x$ are $-ac,c,e-R/2$, this is equivalent to
$d(-ac)+d(c)>2(e-R/2)=2e-R$, i.e. to $R+d(-ac)+d(c)>2e$. So we have
the first claim.

If $R+d(-ac)+d(c)>2e$, then
$\bar g(a,R)\subseteq\hat g(a,R)\subseteq\N (c)$. If
$(a,R)\neq (-1,-2e)$, then $\bar g(a,R)=\hat g(a,R)$ so the reverse is
also true. If $(a,R)=(-1,-2e)$, then
$\bar g(a,R)=\ooo\ffs\subseteq\N (c)$ iff $c=1$ or $\Delta$. If $c=1$,
then $R+d(-ac)+d(c)=\infty >2e$, while if $c=\Delta$, then
$R+d(-ac)+d(c)=-2e+d(\Delta )+d(\Delta )=2e$. Hence the only exception
is $(a,R)=(-1,-2e)$ and $c=\Delta$. \qed

\blm If $a,b\in\fff/\ffs$ and $R\in\RR$, then
$$\hat g(a,R)\hat g(b,R)=\ups{R+d(ab)}\hat g(a,R)
=\ups{R+d(ab)}\hat g(b,R).$$
Unless $R=-2e$ and $\{ a,b\} =\{ -1,-\Delta\}$, the same holds for
$\bar g$. 
\elm

\pf (i) We have two cases:

1. $d(-a),d(-b)>e-R/2$. By Proposition 2.5,
$\hat g(a,R)=\hat g(b,R)=\ups{R/2+e}$. We also have
$d(ab)\geq\min\{ d(-a),d(-b)\} >e-R/2$, so $R+d(ab)>R/2+e$ and so
$(1+\p^{R+d(ab)})\fff^2\subseteq (1+\p^{R/2+e})\fff^2$. Hence all
sides of the relation we want to prove are equal to $\ups{R/2+e}$.

2. $d(-a)$ or $d(-b)\leq e-R/2$. We may assume that, say,
$d(-b)\leq d(-a)$. Then $d(-b)\leq e-R/2$. Then, by Lemma 2.4,
$\hat g(b,R)\subset\ups{R+d(-b)}$ and $[\ups{R+d(-b)}:\hat g(b,R)]=2$,
so there is no intermediate subgroup between $\hat g(b,R)$ and
$\ups{R+d(-b)}$. Then for every subgroup $H$ of $\ups{R+d(-b)}$ we have
$\hat g(b,R)\subseteq H\hat g(b,R)\subseteq\ups{R+d(-b)}$, which
implies that $H\hat g(b,R)=\hat g(b,R)$ or $\ups{R+d(-b)}$. We have
$H\hat g(b,R)=\hat g(b,R)$ iff
$H\subseteq\hat g(b,R)=\upn{R+d(-b)}{-b}$. But
$H\subseteq\ups{R+d(-b)}$, so this is equivalent to
$H\subseteq\N (-b)$. In conclusion, $H\hat g(b,R)=\hat g(b,R)$ if
$H\subseteq\N (-b)$ and $=\ups{R+d(-b)}$ otherwise.

We apply the observation above for $H=\hat g(a,R)$ or
$\ups{R+d(ab)}$. First we verify that these are subgroups of
$\ups{R+d(-b)}$. By adding $R$ to
the inequalities $d(-b)\leq e-R/2,d(-a)$, we get
$R+d(-b)\leq\min\{ R/2+e,R+d(-a)\} =\alpha (a,R)$, so
$\hat g(a,R)\subseteq\ups{\alpha (a,R)}\subseteq\ups{R+d(-b)})$. We
also have $d(ab)\geq\min\{ d(-a),d(-b)\} =d(-b)$, so
$R+d(ab)\geq R+d(-b)$ and so $\ups{R+d(ab)}\subseteq\ups{R+d(-b)}$.
Then, by Lemmas 2.8 and 1.2(i), both $\hat g(a,R)\subseteq\N (-b)$ and
$\ups{R+d(ab)}\subseteq\N (-b)$ are equivalent to $R+d(ab)+d(-b)>2e$.
It follows that 
\begin{multline*}
\hat g(a,R)\hat g(b,R)=\ups{R+d(ab)}\hat g(b,R)\\
=\begin{cases}\upn{R+d(-b)}{-b}&\text{if }R+d(ab)+d(-b)>2e\\
\ups{R+d(-b)}&\text{if }R+d(ab)+d(-b)\leq 2e\end{cases}.
\end{multline*}

If $d(-a)=d(-b)$ then, in particular, $d(-b)\geq d(-a)$ so, by
reversing the roles of $a$ and $b$, we get
$\hat g(a,R)\hat g(b,R)=(1+\p^{R+d(ab)})\fff^2\bar g(a,R)$ and we are
done.

If $d(-a)>d(-b)$ then, by the domination principle,
$d(ab)=d(-b)\leq e-R/2$. Hence $R+d(ab)+d(-b)=R+2d(-b)\leq 2e$ so
$\hat g(a,R)\hat g(b,R)=\ups{R+d(-b)}=\ups{R+d(ab)}$. Since also
$\hat g(a,R)\subseteq\ups{R+d(-b)}=\ups{R+d(ab)}$, we also have
$\ups{R+d(ab)}\hat g(a,R)=\ups{R+d(ab)}$ and we are done. 

If $(R,a),(R,b)\neq (-1,-2e)$, then  $\bar g(a,R)=\hat g(a,R)$ and
$\bar g(b,R)=\hat g(b,R)$, so replacing $\hat g$ by $\bar g$ doesn't
change our statement. 

Suppose now that, say, $(b,R)=(-1,-2e)$. Then
$\bar g(b,R)=\ooo\fff^2$. If $a=b=-1$, then $\bar g(a,R)=\ooo\fff^2$
and $\ups{R+d(ab)})=\ups\infty=\fff^2$, so all sides of the relations
we want to prove are equal to $\ooo\fff^2$. Suppose now that
$a\neq -1$. Then $ab\neq 1$, so $d(ab)\leq 2e=-R$ and so
$R+d(ab)\leq 0$. It follows that $\ups{R+d(ab)}=\fff$ and so
$\ups{R+d(ab)}\bar g(a,R)=\ups{R+d(ab)}\bar g(b,R)=\fff$. On the other
hand, $a\neq -1$, so $d(-a)\leq 2e=-R$ and so $R+d(-a)\leq 2e$. By
Lemma 2.6(ii), $\bar g(a,R)=\N (-a)$. If $a\neq -\Delta$, then
$\N (-a)\not\subseteq\ooo\ffs$, so
$\bar g(a,R)\bar g(b,R)=\N (-a)\ooo\ffs =\fff$ and we are done.

Similarly, if $(a,R)=(-1,-2e)$, then
$\bar g(a,R)\bar g(b,R)=\ups{R+d(ab)}\bar g(a,R)
=\ups{R+d(ab)}\bar g(b,R)$, except when
$b=-\Delta$. This concludes the proof. \qed

\bco If $a,b\in\fff/\ffs$ and $R\geq S$, then
$$\hat g(a,R)\hat g(b,S)=\ups{R+d(ab)}\hat g(b,S).$$

The same holds for $\bar g$, unless $R=-2e$ and
$(a,b)=(-1,-\Delta )$ or $R=S=-2e$ and $(a,b)=(-\Delta,-1)$.
\eco
\pf By Lemma 2.6(i), $\hat g(b,R)\subseteq\hat g(b,S)$, so
$\hat g(b,R)=\hat g(b,R)\hat g(b,S)$. Then, when we multiply by
$\hat g(b,S)$ the relation
$\hat g(a,R)\hat g(b,R)=\ups{R+d(ab)}\hat g(b,R)$, from Lemma 2.9, we
get $\hat g(a,R)\hat g(b,S)=\ups{R+d(ab)}\hat g(b,S)$.

The proof also holds for $\bar g$, except when $R=-2e$ and
$\{ a,b\} =\{ -1,-\Delta\}$, in which case
$\hat g(a,R)\hat g(b,R)=\ups{R+d(ab)}\hat g(b,R)$ no longer
holds. However, if $(a,b)=(-\Delta,-1)$ and $S<R=-2e$, then
$R+d(ab)=-2e+d(\Delta )=0$, so $\ups{R+d(ab)}=\fff$. And, since
$S<-2e$, we have $(b,S)\notin\mathcal A$, so, by Lemma 2.6(ii),
$\bar g(b,S)=\N (-b)=\N (1)=\fff$. It follows that
$\hat g(a,R)\hat g(b,R)=\ups{R+d(ab)}\hat g(b,S)=\fff$. Hence for
$(a,b)=(-\Delta,-1)$ we have only one exception, when $R=S=-2e$. \qed

\bdf We define $G:\fff/\fff^2\times\RR\to Sgp(\fff/\fff^2)$ by
$$\bar G(a,R)=\langle a\rangle\bar g(a,f(R)),\text{ where }
f:\RR\to\RR,\, f(x)=\min\{ R/2-e,R-2e\}.$$
\edf

\bff Note that $f$ is strictly increasing and continuous, with
$\lim_{x\to\pm\infty}f(x)=\pm\infty$, so $f$ is bijective. We have
$x/2-e\geq$ or $\leq x-2e$ iff $x\leq 2e$. Thus $f(x)=x/2-e$ for
$x\leq 2e$ and $f(x)=x-2e$ for $x\geq 2e$.

For the inverse function, note that $f(2e)=0$. Suppose that
$f(x)=y$. Since $f$ is strictly increasing, if $y\leq 0$, then
$x\leq 2e$, so $y=f(x)=x/2-e$ and so $x=2y+2e$. And if $y\geq 0$, then
$x\geq 2e$, so $y=f(x)=x-2e$ and so $x=y+2e$. Hence $f^{-1}(y)=2y+2e$
for $y\leq 0$ and $f^{-1}(y)=y+2e$ for $y\geq 0$. For short,
$f^{-1}(y)=\min\{ 2y+2e,y+2e\}$

Also note that $f(-2e)=-2e$, so $f(x)=-2e$ iff $x=-2e$.
\eff

\blm (i) For every $a\in\fff/\ffs$ the map $\bar G(a,\cdot )$ is
decreasing.

(ii) If $R\leq 2e$, then the $\langle a\rangle$ factor can be dropped
in the formula for $\bar G(a,R)$. Namely, we have
$$\bar G(a,R)=\begin{cases}\bar g(a,R/2-e)&\text{if }R\leq 2e\\
\langle a\rangle\bar g(a,R-2e)&\text{if }R>2e\end{cases}.$$

(iii) We have $\bar G(-1,-2e)=\ooo\ffs$ and for $(a,R)\neq (-1,-2e)$
in the formula for $\bar G(a,R)$ we may replace $\bar g$ bar $\hat g$.

(iv) $\bar g(a,R)\subseteq\bar G(a,R)\subseteq\N (-a)$.
\elm
\pf (i) By Lemma 2.6(i), the map $\bar g(a,\cdot )$ is decresing and so
is $\langle a\rangle\bar g(a,\cdot )$. Since $\bar G(a,\cdot )$ is the
composion of the decreasing map $\langle a\rangle\bar g(a,\cdot )$ and
the increasing map $f$, it is decreasing.

(ii) By 2.11, we have
$\bar G(a,R)=\langle a\rangle\bar g(a,f(R))=
\langle a\rangle\bar g(a,R/2-e)$, if $R\leq 2e$, and
$=\langle a\rangle\bar g(a,R-2e)$, if $R>2e$. But if
$R\leq 2e$, then $R/2-e\leq 0$, so, by Lemma 2.6(iii),
$a\in\bar g(a,R/2-e)$. Hence in this case
$\langle a\rangle\bar g(a,R/2-e)=\bar g(a,R/2-e)$.

(iii)  By definition, $\bar g(a,f(R))=\hat g(a,f(R))$ for
$(a,f(R))\neq (-1,-2e)$. But, by 2.11, $f(R)=-2e$ iff $R=-2e$, so
$(a,f(R))=(-1,-2e)$ iff $(a,R)=(-1,-2e)$. Hence
$\bar g(a,f(R))=\hat g(a,f(R))$ for $(a,R)\neq (-1,-2e)$.

In the exceptional case, by (ii), we have
$\bar G(-1,-2e)=\bar g(-1,-2e)=\ooo\ffs$.

(iv) For $\bar g(a,R)\subseteq\bar G(a,R)$ it suffices to prove that
$\bar g(a,R)\subseteq\bar g(a,f(R))$. If $R\geq -2e$, then
$R\geq R/2-e$. Since also $R>R-2e$, we have $R\geq f(R)$, so, by Lemma
2.6(i),  $\bar g(a,R)\subseteq\bar g(a,f(R))$. Suppose now that
$R<-2e$. Then $f(R)<f(-2e)=-2e$. It follows that
$(a,R),(a,f(R))\notin\bar{\mathcal A}$ so, by Lemma 2.3(ii),
$\bar g(a,R)=\bar g(a,f(R))=\N (-a)$.

We have $\bar g(a,f(R))\subseteq\N (-a)$, by the definition of
$\bar g$. And $(a,-a)_\p =1$, so $a\in\N (-a)$. Thus
$\bar G(a,R)=\langle a\rangle\bar g(a,f(R))\subseteq\N (-a)$. \qed

\bco We have $\bar G(-1,-2e)=\ooo\ffs$ and if $(a,R)\neq (-1,-2e)$,
then: 

(i) If $R\leq 2e$, then
$\bar G(a,R)=\upn{\min\{ R/4+e/2,R/2+d(-a)-e\}}{-a}$. Altenatively,
$$\bar G(a,R)=\begin{cases}\ups{R/4+e/2}&\text{if }d(-a)>3e/2-R/4\\
\upn{R/2+d(-a)-e}{-a}&\text{if }d(-a)\leq 3e/2-R/4\end{cases}$$

(ii) If $R>2e$, then
$\bar G(a,R)=\langle a\rangle\upn{\min\{ R/2,R+d(-a)-2e\}}{-a}$.
Altenatively, 
$$\bar G(a,R)=\begin{cases}\langle a\rangle\ups{R/2}
&\text{if }d(-a)>2e-R/2\\
\langle a\rangle\upn{R+d(-a)-2e}{-a}
&\text{if }d(-a)\leq 2e-R/2\end{cases}$$

(iii) If $d(-a)\leq e-R/2$ or $(a,R)\notin\bar{\mathcal A}$ (e.g. when
$R<-2e$), then $\bar G(a,R)=\N (-a)$. Same happens if $\ord a$ is odd
and $R\leq 2e+2$.

If $d(-a)>e-R/2$ and $(a,R)\neq (-1,-2e)$, then
$\ups{R/2+e}\subseteq\bar G(a,R)$.

(iv) $\bar G(a,R)\subseteq\langle a\rangle\ups{R-2e}$. In particular,
if $R>4e$, then $\bar G(a,R)=\langle a\rangle\ffs$.
\eco
\pf By Lemma 2.12(iii), $\bar G(-1,-2e)=\ooo\ffs$ and if
$(a,R)\neq (-1,-2e)$, then
$\bar G(a,R)=\langle a\rangle\bar g(a,f(R))=
\langle a\rangle\hat g(a,f(R))$ and the factor $\langle a\rangle$ can
be dropped if $R\leq 2e$. By definition,
$\hat g(a,f(R))=\upn{\alpha (a,f(R))}{-a}$, with
$\alpha (a,f(R))=\min\{ f(R)/2+e,f(R)+d(-a)\}$ and, by Lemma 2.4,
$$\hat g(a,f(R))=\begin{cases}\ups{f(R)/2+e}&\text{if }d(-a)>e-f(R)/2\\
\upn{f(R)+d(-a)}{-a}&\text{if }d(-a)\leq e-f(R)/2\end{cases}.$$

If $R\leq 2e$, then $f(R)=R/2-e$, so $f(R)/2+e=R/4+e/2$,
$f(R)+d(-a)=R/2+d(-a)-e$ and $e-f(R)/2=3e/2-R/4$, so we get (i). If
$R\geq 2e$, then $f(R)=R-2e$, so $f(R)/2+e=R/2$,
$f(R)+d(-a)=R+d(-a)-2e$ and $e-f(R)/2=2e-R/2$, so we get (ii).

(iii) If $d(-a)\leq e-R/2$, then $0\leq e-R/2$, so $R\leq 2e$, and 
$d(-a)\leq -(R/2-e)$, then, by Lemmas 2.12(ii) and 2.6(ii),
$\bar G(a,R)=\bar g(a,R/2-e)=\N (-a)$. If
$(a,R)\notin\bar{\mathcal A}$, then by Lemma 2.3(ii),
$\bar g(a,R)=\N (-a)$, which, by Lemma 2.12(iv), implies
$\bar g(a,R)=\bar G(a,R)=\N (-a)$. If $\ord a$ is odd and $R=2e+2$,
then $d(-a)=0$, so $\min\{ R-2e,R/2+d(-a)-e\} =1$. Then, by (i), we
have $\bar G(a,2e+2)=\langle a\rangle\upn 1{-a}=
\langle a\rangle\ooo\ffs\cap\N (-a)=\N (-a)$. (Since $\ord a$ is odd,
we have $\langle a\rangle\ooo\ffs =\fff$.) If $R\leq 2e+2$, then, by
Lemma 2.12 (i) and (iv), we have $\N (a,R)=\bar G(a,2e+2)\subseteq\bar
G(a,R)\subseteq\N (a,R)$

If $d(-a)>e-R/2$ and $(a,R)\neq (-1,-2e)$, by Lemmas 2.4 and 2.12(iv),
we have $\ups{R/2+e}=\hat g(a,R)=\bar g(a,R)\subseteq\bar G(a,R)$.

(iv) If $R\leq 2e$, then $\ups{R-2e}=\fff$ so the statement is
trivial. And if $R>2e$, by Lemma 2.12(ii),
$\bar G(a,R)=\langle a\rangle\bar g(a,R-2e)$. But, by Lemma 2.6(iv),
$\bar g(a,R-2e)\subseteq\ups{R-2e}$ and if $R>4e$, i.e. if $R-2e>2e$,
then $\bar G(a,R-2e)=\ffs$. Hence the conclusion. \qed 

\bpr If $a\in\fff/\oos$ and $\ord a=R$ then $G(a)=\bar G(a,R)$.
\epr

\pf Note that $(a,R)=(-1,-2e)$ in $\fff/\fff^2\times\RR$ iff $a\in
-\fff^2$ and $\ord a=R=-2e$, which is equivalent to $a=-\frac 14$ in
$\fff/\oos$. So the exception $a=-\frac 14$ of [B2, Definition 4]
corresponds to the exception $(a,R)=(-1,-2e)$ of Corollary 2.13. In
this case we have $G(a)=\bar G(a,R)=\ooo\fff^2$. 

In all other cases we see that the formulas from [B2, Definition 4]
coincide with those from Corollary 2.13. The only not obvious case is
when $d(-a)>3e/2-R/4$ in case (III). Here we must prove that
$\ups{R/4+e/2}=\ups{e-[e/2-R/4]}$. This follows from
$\ups x=\ups{\lceil x\rceil}$ for $x=R/4+e/2$. We note that
$\lceil x\rceil =-[-x]$. Thus
$\lceil R/4+e/2\rceil =-[-e/2-R/4]=-(-e+[e/2-R/4])=e-[e/2-R/4]$. This
concludes the proof. \qed

\bco If $R\in\RR$ and $a,c\in\fff/\fff^2$, then
$\bar G(a,R)\subseteq\N (c)$ iff $f(R)+d(-ac)+d(c)>2e$ and
$(a,c)_\p =1$ or $(a,R)=(-1,-2e)$ and $c=\Delta$. If $R\leq 2e$, the
condition $(a,c)_\p =1$ is superfluous. Explicitely:

If $R\leq 2e$, then $\bar G(a,R)\subseteq\N (c)$ iff
$R/2+d(-ac)+d(c)>3e$ or $(a,R)=(-1,-2e)$ and $c=\Delta$.

If $R>2e$, then  $\bar G(a,R)\subseteq\N (c)$ iff
$R+d(-ac)+d(c)>4e$ and $(a,c)_\p =1$.
\eco
\pf We have 
$\bar G(a,R)=\langle a\rangle\bar g(a,f(R))\subseteq\N (c)$ iff
$\bar g(a,f(R))\subseteq\N (c)$ and $a\in\N (c)$, i.e. $(a,c)_\p =1$.
By Lemma 2.8, the first condition is equivalent to
$f(R)+d(-ac)+d(c)>2e$ or $(a,f(R))=(-1,-2e)$ and $c=\Delta$. But
$(a,f(R))=(-1,-2e)$ is equivalent to $(a,R)=(-1,-2e)$. (By 2.11,
$f(R)=-2e$ iff $R=-2e$.) And, by Lemma 2.12(ii), if $R\leq 2e$, then
the factor $\langle a\rangle$ may be removed from the formula for
$\bar G(a,R)$, so the condition $(a,c)_\p =1$ is superfluous. This
concludes the proof.

If $R\leq 2e$, which includes the exceptional case, and when the
condition $(a,c)_\p =1$ is superfluous, we have $f(R)=R/2-e$, so
$f(R)+d(-ac)+d(c)>2e$ writes as $R/2+d(-ac)+d(c)>3e$.

If $R>2e$, then $f(R)=R-2e$, so $f(R)+d(-ac)+d(c)>2e$ writes as
$R+d(-ac)+d(c)>4e$. \qed

\blm Let $a,b\in\fff/\ffs$ and let $R,S\in\RR$ with $R\geq S$. Assume
that we don't have $R=-2e$ and $(a,b)=(-1,-\Delta )$ or $R=S=-2e$
and $(a,b)=(-\Delta,-1)$. Then
$$\bar G(a,R)\bar G(b,S)=\langle ab\rangle
\ups{f(R) +d(ab)}\bar G (b,S).$$
Moreover, the factor $\langle ab\rangle$ is superfluous if $R\leq 2e$.
Explicitely:
$$\bar G(a,R)\bar G(S,R)=\begin{cases}
\ups{R/2+d(ab)-e}\bar G (b,S)&\text{if }R\leq 2e\\
\langle ab\rangle\ups{R+d(ab)-2e}\bar G (b,S)&\text{if }R>2e
\end{cases}.$$
Also the factor $\langle ab\rangle$ may be replaced by
$\langle a\rangle$.
\elm
\pf We have $\bar G(a,R)=\langle a\rangle\bar g(a,f(R))$ and
$\bar G(b,S)=\langle b\rangle\bar g(b,f(S))$, so
$$\bar G(a,R)\bar G(b,S)=\langle a,b\rangle\bar g(a,f(R))\bar g(b,f(S)).$$
Moreover, by Lemma 2.12 (ii), if $R\leq 2e$, so $S\leq R\leq 2e$, then
the $\langle a\rangle$ and $\langle b\rangle$ factors can be dropped,
and so does $\langle a,b\rangle$.

For the product $\bar g(a,f(R))\bar g(b,f(S))$ we use Corollary
2.10. By 2.11, $f$ is strictly incresing and $f(-2e)=-2e$ so the
hypothesis implies that $f(R)\geq f(S)$ and we dont have $f(R)=-2e$
and $(a,b)=(-1,-\Delta )$ or $f(R)=f(S)=-2e$ and $(a,b)=(-\Delta,-1)$. 
Hence the hypothesis of Corollary 2.10 is satisfied for $(a,f(R))$ and
$(b,f(S))$. Also the exceptional case doesn't occur so the general
formula 
$\bar g(a,f(R))\bar g(b,f(S))=\ups{f(R)+d(ab)}\bar g(b,f(S))$ holds.
Since $\langle a,b\rangle =\langle a\rangle\langle b\rangle
=\langle ab\rangle\langle b\rangle$ and
$\bar G(b,S)=\langle b\rangle\bar g(b,f(S))$, we have
$\langle a,b\rangle\bar g(b,f(S))=\langle a\rangle\bar G(b,S)
=\langle ab\rangle\bar G(b,S)$. Hence $\bar G(a,R)\bar G(b,S)
=\langle a,b\rangle\ups{f(R)+d(ab)}\bar g(b,f(S))$ writes both as
$\langle ab\rangle\ups{f(R)+d(ab)}\bar G(b,S)$ and
$\langle a\rangle\ups{f(R)+d(ab)}\bar G(b,S)$.

If $R\leq 2e$, then the factors $\langle a\rangle$ and
$\langle b\rangle$ can be ignored and we have
$f(R)+d(ab)=R/2+d(ab)-e$, so we get the claimed formula
$\ups{R/2+d(ab)-e}\bar G(a,S)$. And if $R>2e$, then
$f(R)+d(ab)=R+d(ab)-2e$, so we get the formulas
$\langle ab\rangle\ups{R+d(ab)-2e}\bar G(a,S)$ and
$\langle a\rangle\ups{R+d(ab)-2e}\bar G(a,S)$. \qed

\blm (i) We have $\bar G(a,R)\subseteq\ooo\ffs$ iff either $d(-a)=2e$
or $\ord a$ is even, $R\geq -2e$ and $d(-a)>e-R/2$.

In particular, if $R\geq 2e$, then $\bar G(a,R)\subseteq\ooo\ffs$ iff
$\ord a$ is even.

(ii) Depending on $R$, $\bar G(a,R)\subseteq\ooo\ffs$ is equivalent
to: $d(-a)=2e$ if $R<-2e$; $d(-a)\geq 2e$ if $R=2e$; $d(-a)>e-R/2$ if
$-2e<R\leq 2e$; $d(-a)>0$ if $R\geq 2e$. 
\elm
\pf Assume first that $R\leq 2e$. By Lemma 2.12(ii),
$\bar G(a,R)=\bar g(a,R/2-e)$, so, by Lemma 2.3,
$\bar G(a,R)\subseteq\ooo\ffs$ iff either $a=-\Delta$ or
$R/2-e+2e\geq 0$ and $R/2-e+d(-a)>0$. That is, iff either $d(-a)=2e$
or $R\geq -2e$ and $d(-a)>e-R/2$. The condition that $\ord a$ is even
is superfluous, as it follows from $d(-a)>e-R/2\geq 0$.

If $R\geq 2e$, then $a\in\bar G(a,R)$, so
$\bar G(a,R)\subseteq\ooo\ffs$ implies that $\ord a$ is
even. Conversely, if $\ord a$ is even, then $d(-a)>0=e-2e/2$. By the
case $R\leq 2e$, this implies $\bar G(a,2e)\subseteq\ooo\ffs$. Since
$R\geq 2e$, by Lemma 2.12(i), we get
$\bar G(a,R)\subseteq\bar G(a,2e)\subseteq\ooo\ffs$. The condition
$d(-a)>e-R/2$ is superfluous in this case, since it follows from
$d(-a)>0$ and $R\geq 2e$.

(ii) In the cases $R<-2e$ and $R\geq 2e$ the statement follows
directly from (i). If $-2e\leq R\leq 2e$, then
$\bar G(a,R)\subseteq\ooo\ffs$ iff either $d(-a)=2e$ or
$d(-a)>e-R/2$. If $R=-2e$, this means either $d(-a)=2e$ or
$d(-a)>e-(-2e)/2=2e$, i.e. $d(-a)\geq 2e$. If $-2e<R\leq 2e$, then
$e-R/2<e-(-2e)/2=2e$ so if $d(-a)=2e$, then also $d(-a)>e-R/2$. Hence
the case $d(-a)=2e$ is superfluous and we are left with the condition
$d(-a)>e-R/2$. \qed

\bco If $\bar G(a_i,R_i)\subseteq\ooo\ffs$ for $1\leq i\leq s$ and
$R\in\RR$, with $R\geq -2e$ and $R\geq R_i$ for $1\leq i\leq s$, then
$\bar G((-1)^{s-1}a_{1,s},R)\subseteq\ooo\ffs$.
\eco
\pf For every $i$ we have $R\geq R_i$ so, by Lemma 2.12(i),
$\bar G(a_i,R)\subseteq\bar G(a_i,R_i)\subseteq\ooo\ffs$. By Lemma
2.17(ii), this implies that $d(-a_i)$ is $\geq 2e$, $>e-R/2$ or $>0$,
corresponding to $R=-2e$, $-2e<R\leq 2e$ and $R\geq 2e$,
respectively. By the domination principle, this implies that
$d(-(-1)^{s-1}a_{1,s})=d((-1)^sa_{1,s}$ is $\geq 2e$, $>e-R/2$ or
$>0$, accordingly. Then, again by Lemma 2.17(ii),
$\bar G((-1)^{s-1}a_{1,s},R)\subseteq\ooo\ffs$. \qed

\blm (i) We have $\Delta\in\bar g(a,R)$ iff $\ord a$ is even and
$R\leq 2e$. We have $\Delta\in\bar g(a,R)$ iff either $\ord a$ is even
and $R\leq 2e$ or $a=\Delta$.

(ii) If $a\in\mathcal A$, then $\Delta\in g(a)$ iff $\ord a$ is even
and $\leq 2e$.
\elm
\pf (i) By Lemma 2.12(iv) a necessary condition for both
$\Delta\in\bar g(a,R)$ and $\Delta\in\bar G(a,R)$ is
$\Delta\in\N (-a)$, which is equivalent to $2\mid\ord a$. (We have
$\Delta\in\N (-a)$ iff $-a\in\N (\Delta )=\ooo\ffs$.) Assume now that
$\ord a$ is even. Then $d(-a)>0=e-(2e)/2$, so, by Lemma 2.4,
$\bar g(a,2e)=\ups{(2e)/2+e}=\ups{2e}\ni\Delta$. Then, Lemma 2.12(ii),
we have $\bar 2.12(ii)$, we have
$\bar G(a,4e)=\langle a\rangle\bar g(a,2e)\ni\Delta$. By Lemmas 2.6(i)
and 2.12(i), if $R\leq 2e$, then
$\Delta\in\bar g(a,2e)\subseteq\bar g(a,R)$ and if $R\leq 4e$, then
$\Delta\in\bar G(a,4e)\subseteq\bar g(a,R)$. By Lemma 2.6(iv), if
$R>2e$, then $\bar g(a,R)=\ffs$, so $\Delta\notin\bar g(a,R)$. And,
by Corollary  2.13(iv), if $R>4e$, then
$\bar G(a,R)=\langle a\rangle\ffs$, so $\Delta\in\bar G(a,R)$ iff
$a=\Delta$.

(ii) By Proposition 2.5, $g(a)=\bar g(a,\ord a)\cap\ooo$. Since
$\Delta\in\ooo$, we have $\Delta\in g(a)$ iff
$\Delta\in\bar g(a,\ord a)$, which, by (i), happens iff $\ord a$ is
even and $\ord a\leq 2e$. \qed

\section{Preliminary results} 

We start by proving the necessity of condition (i) of the
representation theorem announced in [B3, Theorem 4.5] when $i=1,2$.

\blm Let $L,K$ be quadratic lattices with $K\subseteq L$ and let
$R_i=R_i(L)$ and $S_i=R_i(K)$.

(i) $R_1\leq S_1$.

(ii) If $\rank K\geq 2$, then either $R_2\leq S_2$ or $\rank L\geq 3$
and $R_2+R_3\leq S_1+S_2$.

(iii) If $\rank K\geq 2$, then $R_1+R_2\leq S_1+S_2$.
\elm

\pf Let $L\cong{\prec a_1,\ldots,a_m\succ}$ and
$K\cong{\prec b_1,\ldots,b_n\succ}$ relative to the good BONGs
$x_1,\dots,x_m$ and $y_1,\ldots,y_n$.

(i) It follows by taking orders in $\nnn K\subseteq\nnn L$.

Next we note that (iii) follows from (i) and (ii). Indeed, by (i) we
have $R_1\leq S_1$ so if $R_2\leq S_2$, then $R_1+R_2\leq S_1+S_2$. If
$R_2>S_2$, then, by (ii), we have $R_1+R_2\leq R_2+R_3\leq S_1+S_2$.
So we only have to prove (ii).

By [B2, Lemma 2.7(iii)], there is a lattice
$J\cong{\prec b_1,b_2\succ}$ relative to the BONG $y_1,y_2$ and we
have $J\subseteq K\subseteq L$. Moreover, $R_i(J)=R_i(K)=S_i$ for
$i=1,2$.

We use induction on $m:=\rank L$. If $m=2$, then $FL=FJ$ so
$J\subseteq L$ implies $L^\#\subseteq K^\#$. By (i) and
Lemma 1.5(iii) we get $-S_2=R_1(K^\# )\leq R_1(L^\# )=-R_2$ so
$R_2\leq S_2$. 

Suppose now that $n\geq 3$.  By Lemma 1.3(ii), we have
$S_1+S_2=\ord\vol J$. In particular, (iii) writes as
$\ord\vol J\geq R_1+R_2$ .

Let  $L'=pr_{x_1^\perp}L$. Then $L'\cong{\prec a_2,\ldots,a_m\succ}$
relative to the good BONG $x_2,\ldots,x_n$ and the sequence
$R_1(L'),\ldots,R_{m-1}(L')$ is $R_2,\ldots,R_m$. In particular,
$\nnn L'=\p^{R_2}$. By the (iii) part of the induction hypothesis, for
every binary sublattice $J'\subseteq L'$ we have
$\ord\vol J'\geq R_1(L')+R_2(L')=R_2+R_3$.

Recall that if $J$ has arbitrary rank $k$ and $v_1,\ldots,v_k$ is
a basis, then $\vol J$ is the ideal generated by
$d_B(v_1,\ldots,v_k):=\det (B(v_i,v_j))_{1\leq i,j\leq k}$. In the
binary case, if $x,y$ is a basis of $J$ then $\vol J$ is the ideal
generated by $d_B(x,y)=Q(x)Q(y)-B(x,y)^2$.

Let $z,t\in pr_{x_1^\perp}L=L'$ be the projections of $x,y$ on
$x_1^\perp$. Then $x=\alpha x_1+z$ and $y=\beta x_1+t$ for some
$\alpha,\beta\in F$ and $z,t\in L'$. If $\alpha =\beta =0$ then
$x,y\in L'$  so $J\subseteq L'$, which implies
$S_1+S_2=\ord\vol J\geq R_2+R_3$ so we are done. So we may assume
that, say, $\alpha\neq 0$ and, moreover, $\ord\alpha\leq\ord\beta$. It
follows that $\alpha^{-1}\beta\in\oo$ so we have another basis for
$J$, $x,y'$, where
$y'=y-\alpha^{-1}\beta x=t-\alpha^{-1}\beta z\in L'$.  

We have
$$\begin{aligned}
d_B(x,y')&=Q(\alpha x_1+z)Q(y')-B(\alpha x_1+z,y')^2\\
&=(Q(\alpha x_1)+Q(z))Q(y')-B(z,y')^2=Q(\alpha x_1)Q(y')+d_B(z,y')
\end{aligned}$$
It follows that $S_1+S_2=\ord d_B(x,y')\geq\min\{\ord Q(\alpha
x_1)Q(y'),\ord d_B(x,y')\}$.

Now $Q(x)=Q(\alpha x_1)+Q(z)$ and we have $x\in J$ and $z\in L'$, so
$\ord Q(x)\geq\ord\nnn J=S_1$ and $\ord Q(z)\geq\ord\nnn
L'=R_2$. Hence $\ord Q(\alpha x_1)\geq\min\{\ord Q(x),\ord Q(z)\}
\geq\min\{ S_1,R_2\}$. On the other hand, $y'$ belongs to both
$J$ and $L'$ so $\ord Q(y')$ is both $\geq\ord\nnn J=S_1$ and
$\geq\ord\nnn L'=R_2$ so $\ord Q(y')\geq\max\{ S_1,R_2\}$. It follows
that
$$\ord Q(\alpha x_1)Q(y')=\ord Q(\alpha x_1)+\ord Q(y')\geq
\min\{ S_1,R_2\} +\max\{ S_1,R_2\} =S_1+R_2.$$

We also have $J':=\oo z+\oo y'\subseteq L'$ so
$\ord d_B(z,y')=\ord\vol (L')\geq R_2+R_3$ if $J'$
is binary and nonegenerate and $d_B(z,y')=0$ otherwise. In both cases
$\ord d_B(z,y')\geq R_2+R_3$. In conclusion,
$$S_1+S_2\geq\min\{ S_1+R_2,R_2+R_3\},$$
so either $S_1+S_2\geq S_1+R_2$, i.e. $S_2\geq R_2$, or
$S_1+S_2\geq R_2+R_3$. This concludes the proof. \qed

\bco If $\rank L\geq 2$ and $R_i(L)=R_i$ then $R_1+R_2=\min\{\ord\vol
J\mid\, J\subseteq L,\, J\text{ binary}\}$.
\eco
\pf If $J\subseteq L$ is binary and $R_i(J)=S_i$, then, by Lemma
3.1(iii), $\ord\vol J=S_1+S_2\geq R_1+R_2$. If $x_1,\ldots,x_n$ is a
good BONG for $L$ then, by [B2, Lemma 2.7(iii)], we have
$J:={\prec x_1,x_2\succ}\subseteq L$. We 
have $R_i(J)=\ord Q(x_i)=R_i$ for $i=1,2$ so $\ord\vol J=R_1+R_2$. Hence
the conclusion. \qed

\blm Let $L$ be a lattice with $\rank L\geq 2$ and let
$R_i(L)=R_i$. Suppose that $K$, with $R_i(K)=S_i$, is maximal amongst
all sublattices of $L$ with $FK=FL$ and $R_1+R_2<S_1+S_2$. Then
$S_1+S_2\leq R_1+R_2+2$. 
\elm
\pf Since $S_1+S_2\neq R_1+R_2$, we have $K\neq L$ so the inclusion
$K\subseteq L$ is strict. We consider an intermediate lattice
$K\subset K'\subseteq L$ with $[K':K]=\p$. If $R_i(K')=S'_i$, then, by
Lemma 3.1(iii), we have $R_1+R_2\leq S'_1+S'_2$ and, by the maximality
of $K$, we cannot have  $R_1+R_2<S'_1+S'_2$. Thus
$S'_1+S'_2=R_1+R_2$. By Corollary 3.2, there is a binary lattice
$J'\subseteq K'$ such that $\ord\vol J'=S'_1+S'_2=R_1+R_2$.

If $J=J'\cap K$, then we have an injective map $J'/J\to K'/K$, so we
have $[J':J]\mid [K':K]=\p$. Hence $[J':J]=\oo$
or $\p$. It follows that $\vol J=[J':J]^2\vol J'=\vol J'$ or
$\p^2\vol J'$. We have $\ord\vol J'=R_1+R_2$ so, by taking orders, we
get $\ord\vol J=R_1+R_2$ or $R_1+R_2+2$. But $J\subseteq K$ and
$R_i(K)=S_i$ so, by Corollary 3.2, $S_1+S_2\leq\ord\vol J$. If
$\ord\vol J=R_1+R_2$, we get $S_1+S_2\leq R_1+R_2$, which contradicts
the hypothesis. Thus $S_1+S_2\leq\ord\vol J=R_1+R_2+2$. \qed

\bco Let $L$ and $J$ be quadratic lattices of ranks $m\geq 2$ and $n$,
with $J\subseteq L$. Let $R_i(L)=R_i$ and $R_i(J)=T_i$. If either
$\rank J\leq 1$ or $\rank J\geq 2$ and $R_1+R_2<T_1+T_2$, then there
is some lattice $K$ with $J\subseteq K\subseteq L$, $FK=FL$ and
$R_1+R_2<S_1+S_2\leq R_1+R_2+2$.
\eco
\pf We first reduce to the case when $m=n$, i.e. when $FJ=FL$. If
$m>n$, then let $J'$ be a lattice over the orthogonal complement of
$FJ$ in $FL$, so that $FJ\perp FJ'=FL$. Let
$J\cong{\prec b_1,\ldots,b_n\succ}$ and
$J'\cong{\prec b_{n+1},\ldots,b_m\succ}$ relative to some good
BONGs. Let $s\gg 0$. Then
$\p^s J'\cong{\prec\pi^{2s}b_{n+1},\ldots,\pi^{2s}b_m\succ}$ relative
to some good BONG and we have $T_i:=\ord\pi^{2s}b_i=2s+\ord b_i\gg 0$
for $n+1\leq i\leq m$. Then $T_{n-1}\leq T_{n+1}$ (if $n\geq 2$),
$T_n\leq T_{n+1}$ (if $n\geq 1$) and $T_n\leq T_{n+2}$ (if $n\geq 1$
and $m\geq n+2$). By Lema 1.5(ii), we have
$\bar J:=J\perp\p^sJ'\cong{\prec b_1,\ldots,b_n,
\pi^{2s}b_{n+1},\ldots,\pi^{2s}b_m\succ}$ relative to some good BONG
and so $R_i(\bar J)=T_i$ $\forall i$. If $n\geq 2$ then
$R_1+R_2<T_1+T_2$ follows from the hypothesis. If $n=1$ it follows
from $T_2\gg 0$ and for $n=0$ it follows from $T_1,T_2\gg 0$.

We cannot have an infinite sequence of lattice
$\bar J=K_1\subset K_2\subset\cdots\subset L$, since we would have an
infinite sequence of ideals
$\vol K_1\subset\vol K_2\subset\cdots\subset\vol L$. It 
follows that there is some lattice $K$, with
$J\subseteq\bar J\subseteq K\subset L$, which is maximal with the
property that $R_1(K)+R_2(K)>R_1+R_2$. If $R_i(K)=S_i$, then
$R_1+R_2<S_1+S_2$ and by the maximality of $K$, we get from Lemma
3.3 that $S_1+S_2\leq R_1+R_2+2$. \qed

\blm If $L\cong{\prec a_1,\ldots,a_n\succ}$ realtively to a (possibly
bad) BONG and $\ord a_i=R_i$, in particular, if $R_i=R_i(L)$ then:

(i) $R_{j+1}-R_j$ cannot be odd and negative for any $1\leq j\leq n-1$.

(ii) If $R_j=R_{j+2}$ for some $1\leq j\leq n-2$ then $R_j\equiv
R_{j+1}\pmod 2$. 
\elm
\pf (i) By [B2, Lemma 2.7(i)] for every $1\leq i\leq n-1$ the lattice
$J={\prec x_j,x_{j+1}\succ}$ exists. Then, by Lemma 1.4, $R_{j+1}-R_j$
cannot be odd and negative.

(ii) We have $(R_{j+1}-R_j)+(R_{j+2}-R_{j+1})=R_{j+2}-R_j=0$ so if
$R_j\not\equiv R_{j+1}\pmod 2$, then one of $R_{j+1}-R_j$ and
$R_{j+2}-R_{j+1}$ is odd and negative. But this contradicts (i). Thus
$R_j\equiv R_{j+1}\pmod 2$. \qed

\blm Let $x_1,x_2,x_3$ be a bad BONG of a ternary quadratic lattice
$L$ and let $R_i=\ord Q(x_i)$. Then $R_2\equiv R_3\pmod 2$ and
$R_i(L)=S_i$, with $(S_1,S_2,S_3)=(R_1,R_2+R_3-R_1,R_1)$.
\elm
\pf We have $S_1=R_1(L)=\ord\nnn L=\ord Q(x_1)=R_1$. Since $L$ has a
bad BONG, by Proposition 1.6(i) it cannot have property A, i.e. we
cannot have $S_1<S_3$. Thus $S_3=S_1=R_1$. For the value of $S_2$ we consider
the volume of $L$. By Lemma 1.3(ii) we have both
$\ord\vol L=R_1+R_2+R_3$ and $\ord\vol L=S_1+S_2+S_3$. Thus
$R_1+R_2+R_3=S_1+S_2+S_3=R_1+S_2+R_1$, which implies
$S_2=R_2+R_3-R_1$, as claimed.

By Lemma 3.5(ii), since $S_1=S_3$, we have $S_1\equiv S_2\pmod 2$,
i.e. $R_1\equiv R_2+R_3-R_1\pmod 2$. Thus $R_2\equiv R_3\pmod 2$. \qed

\blm (i) Let $L$ be a quadratic lattice of rank $n\geq 2$ and let
$R_i(L)=R_i$. Let $x\in L$ be a norm generator. Then
$\ord\nnn pr_{x^\perp}L\geq R_2$, with equality iff $x$ is a first
element in a good BONG of $L$.

(ii) If $n\geq 3$ and $R_2-R_1\geq 2e$, then every norm generator of
$L$ is a first element in a good BONG of $L$.
\elm 
\pf (i) If $x$ is a first element of a good BONG
$x=x_1,x_2,\ldots,x_n$ of $L$, then $\ord Q(x_i)=R_i$. Then $x_2$ is a
norm generator for $pr_{x^\perp}L={\prec x_2,\ldots,x_n\succ}$, so
$\ord\nnn pr_{x^\perp}L=\ord Q(x_2)=R_2$ and we are done.

Suppose now that $x\in L$ is a norm generator, so
$\ord Q(x)=\ord\nnn L=R_1$, but is not a first element of norm
generator. Let $y_2,\ldots,y_n$ be a good BONG for
$L':=pr_{x_1^\perp}L$. Then $x=y_1,y_2,\ldots, y_n$ is a BONG of
$L$. Let $S_i=\ord Q(y_i)$. In particular, $S_1=\ord Q(x)=R_1$ and
$S_2=\ord\nnn{\prec y_2,\ldots,y_n\succ} =\ord\nnn pr_{x^\perp}L$. Since
$x=y_1$ is not a first element in a good BONG, the BONG
$y_1,\ldots,y_n$ is bad, so $S_i>S_{i+2}$ for some $i$. But the BONG
$y_2,\ldots,y_n$ is good, so $S_i\leq S_{i+2}$ for $2\leq i\leq n-2$.
So we must have $S_1>S_3$. By [B2, Lemma 2.7(iii)], the lattice
$M={\prec y_1,y_2,y_3\succ}$ exists and $M\subseteq L$. If
$R_i(M)=T_i$, then, by Lemma 3.1(iii), $R_1+R_2\leq T_1+T_2$. Since
$S_1>S_3$, the BONG $y_1,y_2,y_3$ of $M$ is bad. Hence, by Lemma 3.6,
we have $(T_1,T_2,T_3)=(S_1,S_2+S_3-S_1,S_1)$. Thus
$R_1+R_2\leq T_1+T_2=S_2+S_3$. Together with $R_1=S_1>S_3$, this
implies $R_2<S_2=\ord\nnn pr_{x^\perp}L$, which concludes the proof of
(i).

(ii) Suppose that $R_2-R_1\geq 2e$ and there is $x\in L$ which is a
norm generator, but not a first element in a good BONG. Using the
notations from the proof of (i), we have $R_2<S_2$ and $R_1>S_3$, It
follows that $S_2-S_3>R_2-R_1\geq 2e$. But this is impossible since
$y_2,\ldots,y_n$ is a good BONG so, by Lemma 1.4, $S_3-S_2'\geq -2e$.
\qed 

\blm Let $L$ be a quadratic lattice of rank $n\geq 3$, with
$R_i(L)=R_i$.

If $R_1<R_3$ or $R_3-R_2=-2e$, then every norm generator of $L$ is a
first element in a good BONG.
\elm
\pf If $R_3-R_2=-2e$, then $R_1-R_2\leq R_3-R_2=2e$, so
$R_2-R_1\geq 2e$ and so Lemma 3.7(ii) applies. So we are left with the
case $R_1<R_3$.

Let $x_1,\ldots,x_n$ be a good BONG of $L$ and let
$L=L_1\perp\cdots\perp L_t$ be a Jordan splitting. By [OM, 91.9], $t$,
$\rank L_i$ and $r_i=\ord\sss L_i$ are invariants of $L$, i.e. they
don't depend on the choice of the Jordan splitting.

Let $x\in L$ be a norm generator and let $K=pr_{x^\perp}L$. By Lemma
3.7(i), $\ord\nnn K\geq R_2$, with equality iff $x$ is a first element
in a good BONG. So it is enought to prove that $\ord\nnn K\leq R_2$.

We have two cases:

(1) $R_1<R_2$. By Lemma 1.5(ii), we have
$L={\prec x_1\succ}\perp L'=\oo x_1\perp L'$, where
$L'={\prec x_2,\ldots,x_n\succ}$. By Lemma 1.5(i),
$\nnn L=\p^{R_2}$ and $\sss L=\p^r$, where
$r=\min\{ R_2,(R_2+R_3)/2\}$. Since $R_2,R_3>R_1$ we have both
$2R_2>R_1+R_2$ and $R_2+R_3>R_1+R_2$. Thus
$2r=\min\{ 2R_2,R_2+R_3\} >R_1+R_2>2R_1$ and so $r>R_1$. Then
$\sss\oo x=\p^{R_1}\supset\p^r=\sss L'$. It follows that we have a
Jordan splitting $L=L_1\perp\cdots\perp L_t$ with $L_1=\oo x_1$ and
$L_2\perp\cdots\perp L_t=L'$. Also $r_1=R_1$ and $r_2=r$ and $L_1$ is
unary. We have $L^{\p^{r_2}}=\p^{r_2-r_1}L_1\perp L'$, i.e.
$L^{\p^r}=\p^{r-R_1}x_1\perp L'$. Since
$\nnn\p^{r-R_1}x_1=\p^{2(r-R_1)}Q(x_1)=\p^{2r-R_1}$,
$\nnn L'=\p^{R_2}$ and $2r>R_1+R_2$, so $R_2<2r-R_1$, we have
$\nnn L=\p^{R_2}$.

Since $\ord Q(x)=R_1=r_1$, and $L_1$ is unary there is another Jordan
decomposition $L=K_1\perp\cdots\perp K_t$, with $K_1=\oo x$. Then
$K=pr_{x^\perp}L=K_2\perp\cdots\perp K_t$, with $\sss K=\p^{r_2}=\p^r$.
Same as before, we have $L^{\p^{r_2}}=\p^{r_2-r_1}K_1\perp K$, i.e.
$L^{\p^r}=\p^{r-R_1}x\perp K$, so
$\nnn L^{\p^r}=\nnn\p^{r-R_1}x+\nnn K$. But $\nnn L^{\p^r}=\p^{R_2}$
and $\nnn\p^{r-R_1}x=\p^{2(r-R_1)}Q(x)=\p^{2r-R_1}$, we get
$\p^{R_2}=\p^{2r-R_1}+\nnn K$. Since $R_2<2r-R_1$, this implies
$\nnn K=\p^{R_2}$, so we are done.

(2) $R_1\geq R_2$. By [B2, Corollary 3.4(iii)],
$J:={\prec x_1,x_2\succ}$ is $\p^r$-modular, with
$r=(R_1+R_2)/2$. Since $R_3>R_1\geq R_2$, by Lemma 1.5(ii),
$L=J\perp L'$, where $L'={\prec x_3,\ldots,x_n\succ}$.
By Lemma 1.5(i), $\nnn L'=\p^{R_3}$ and $\sss L=\p^s$, with
$s=\min\{ R_3,(R_3+R_4)/2\}$, with $(R_3+R_4)/2$ ignored if $n=3$. We
have $R_3>R_1\geq R_2$, so $2R_3>R_1+R_2$, and we have $R_3>R_1$ and
$R_4>R_2$, so $R_3+R_4>R_1+R_2$. Thus
$2s=\min\{ 2R_3,R_3+R_4\} >R_1+R_2=2r$ and so $s>r$. Then we have a
Jordan decomposition $L=L_1\perp\cdots\perp L_t$, with $L_1=J$ and
$L_2\perp\cdots\perp L_t=L'$. Also $r_1=r$ and $r_2=\ord\sss L'=s$. We
have two subcases:

(a) $R_1=R_2$. Then $r_1=r=(R_1+R_2)/2=R_1$. Since $x\in L$ and
$\ord Q(x)=R_1=r_1=\ord\sss L$, we have that $\oo x$ splits
$L$. Thince $pr_{x^\perp}L=K$, this means that $L=\oo x\perp K$. Since
the first term of the Jordan decomposition of $L$ is $L_1$, of scale
$\p^{r_1}=\p^{R_1}$ and rank $2$, and $\oo x$ is $\p^{R_1}$-modular of
rank $1$, the first term in a Jordan decompostion of $K$ is
$\p^{R_1}$-modular of rank $1$. Hence $K$ writes as $K=\oo y\perp K'$,
with $\ord Q(y)=R_1$. It follows that
$\ord\nnn K\leq\ord Q(y)=R_1=R_2$ and we are done.

(b) $R_1>R_2$. Then $R_1>(R_1+R_2)/2=r$, so
$\nnn L=\p^{R_1}\subset\p^r=\sss L$. Now $J$ is $\p^r$-modular and,
since $s>r$, we have $\sss L'=\p^s\subseteq\p^{r+1}$. It follows that
$L^{\p^{r+1}}=(J\perp L')^{\p^{r+1}}=\p J\perp L'$. Since
$\ord\nnn\p J=\ord\nnn J+2=R_1+2>R_1$ and $\ord\nnn L'=R_3>R_1$, we
have $\ord\nnn L^{\p^{r+1}}>R_1$.   

We have $B(x,L)\subseteq\sss L=\p^r$. If $B(x,L)\subseteq\p^{r+1}$,
then $x\in L^{\p^{r+1}}$ so
$R_1=\ord Q(x)\geq\ord\nnn L^{\p^{r+1}}>R_1$. Contradiction. So
$B(x,L)=\p^r$. Let $y\in L$ with $B(x,y)\oo =\p^r=\sss L$. Since
$\nnn L\subset\sss L$, the lattice $J':=\oo x+\oo y$ is binary and
$\sss L=\p^r$-modular (see [OM, \S91C]). We have
$x\in J'\subseteq L$  and $x$ is a norm generator of $L$ so it is a
norm generator for $J'$. Then $J'={\prec x,x'\succ}$, for some
$x'\in FJ'$. We have $\ord Q(x)=R_1$ so if $\ord Q(x')=R'_2$
then $\p^{(R_1+R'_2)/2}=\sss J=\sss L=\p^{(R_1+R_2)/2}$ so
$R'_2=R_2$. Since $x'\in pr_{x^\perp}J\subseteq pr_{x^\perp}L=K$, we get
$\ord\nnn K\leq\ord (Q'(x'))=R_2$. \qed

\bco Let $M$ and $N$ be quadratic lattices with $N\subseteq M$,
$\nnn M=\nnn N$, $\rank M=m\geq 3$, $\rank N=n$, $R_i(M)=R_i$ and
$R_i(N)=S_i$. Then in each of the following cases every element in a
good BONG of $N$ is also an element of a good BONG of $M$.

$$(a)~~ R_1<R_3,\quad (b)~~ R_3-R_2=-2e,\quad (c)~~ n\geq 2\text{ and
}R_2=S_2.$$
\eco
\pf Since $\nnn N=\nnn M$ every norm generator of $N$, in particular,
every first element of a good BONG of $N$, is also a norm generator of
$M$. In the cases (a) and (b), by Lemma 3.8, every norm generator of
$M$ is a first element of a good  BONG so we are done.

If $R_2=S_2$ then, by Lemma 3.7(i), for a first element $x$ of a good
BONG of $N$ we have $\ord\nnn pr_{x^\perp}N=S_2=R_2$. Since $x$ is a
norm generator of $M$ we also have
$\ord\nnn pr_{x^\perp}M\geq R_2$. But $N\subseteq M$ so
$\nnn pr_{x^\perp}N\subseteq\nnn pr_{x^\perp}M$ and so
$R_2=\ord\nnn pr_{x^\perp}N\geq\ord\nnn pr_{x^\perp}M\geq R_2$. Hence
$\ord\nnn pr_{x^\perp}M=R_2$ so, by Lemma 3.7(i), $x$ is a first element
in a good BONG of $M$. \qed

\section{The set of all not norm generators of a lattice}

In this section for every quadratic lattice $L$ we denote by
$$L'=\{ x\in L\mid\, x\text{ not a norm generator}\}
=\{ x\in L\mid\, Q(x)\in\p\nnn L\}.$$

We fix a matrix $L\cong{\prec a_1,\ldots,a_n\succ}$, relative to some
good BONG $x_1,\ldots,x_n$ and we denote by $R_i(L)=R_i$.

\blm If $x,x'$ are elements of a quadratic lattice such that
$\ord Q(x+x')=R<\ord Q(x),\ord Q(x')$ then
$\oo x+\oo x'\cong\frac 12A(0,0)$.
\elm
\pf This is [B2, Lemma 3.19].

\blm If $V$ is an anisotropic quadratic space and $\mathfrak a$ is a
fractionary ideal, then $L:=\{ x\in V\mid\, Q(x)\in\mathfrak a\}$ 
is the only $\mathfrak a$-maximal lattice on $V$.

We also have $O(L)=O(V)$.
\elm
\pf The first statement is just [OM, Theorem 91.1]. For the second we
note that if $\sigma\in O(V)$, then $\sigma (L)\cong L$ so, same as
$L$, $\sigma (L)$ is $\mathfrak a$-maximal. By the unicity of the
$\mathfrak a$-maximal lattice over $V$, we have $L=\sigma (L)$ so
$\sigma\in O(L)$. Hence $O(L)=O(V)$. \qed

\blm (i) If $L=J\perp K$, with $\nnn L=\nnn J\supset\nnn K$, then
$L'=J'\perp K$.

(ii) $L'=L'+\p L$. In particular, $\p L\subseteq L'$.

(iii) $L'$ is not a lattice iff $\nnn L=2\sss L$ and $L$ splits a
hyperbolic lattice of scale $\sss L$, . 

(iv) If $\nnn L\supset 2\sss L$ then $L/L'\cong\oo/\p$.
\elm
\pf Let $R=\ord\nnn L$ so $x\in L'$ iff $x\in L$ and $\ord Q(x)>R$.

(i) We have $R=\ord\nnn L=\ord\nnn J<\ord\nnn L$. If $x\in L$ then
$x=y+z$ with $y\in J$, $z\in K$. Then $x\in L'$ iff $\ord Q(x)>R$.
Since $Q(x)=Q(y)+Q(z)$ and $\ord Q(z)\geq\ord\nnn K>R$, this 
is equivalent to $\ord Q(y)>R=\ord\nnn J$, i.e. to $y\in J'$. Hence
$L'=J'\perp K$.

(ii) We must prove that $x+\alpha y\in L'$ for every if $x\in L'$,
$y\in L$ and $\alpha\in\p$. We have
$Q(x+\alpha y)=Q(x)+\alpha^2Q(y)+2\alpha B(x,y)R$. Since
$Q(y)\in\nnn L$, $2B(x,y)\in 2\sss L\subseteq\nnn L$ and
$\alpha\in\p$, we have
$\alpha^2Q(y)+2\alpha B(x,y)\in\p\nnn L=\p^{R+1}$ so
$\ord (\alpha^2Q(y)+2\alpha B(x,y))>R$. Since also $x\in L'$, so
$\ord Q(x)>R$, we get $\ord Q(x+\alpha y)>R$. Thus
$x+\alpha y\in L'$. 

(iii) If $x\in L'$, $\alpha\in\oo$, then
$\ord Q(\alpha x)\geq\ord Q(x)>R$ so $\alpha x\in L'$. So $L'$ is
closed under multiplication with scalars from $\oo$. It is a lattice
iff it is also closed under addition. Suppose the contrary, i.e. there
are $x,y\in L'$ such that $x+y\notin L'$, i.e.
$\ord Q(x),\ord Q(y)>R$, but $\ord Q(x+y)=R$. By Lemma 4.1, this
implies that $H:=\oo x+\oo y\cong\frac 12\pi^RA(0,0)$. We have
$H\subseteq L$, so
$\frac 12\p^R=\sss J\subseteq\sss L\subseteq \frac 12\nnn L=\frac
12\p^R$, so $\sss L=\sss H=\frac 12\nnn L$ and so $H$ splits
$L$. Conversely, if $\p^R=\nnn L=2\sss L$ and $L$ splits a hyperbolic
lattice $H$ with $\sss H=\sss L$, i.e. $H\cong\frac 12\pi^RA(0,0)$
relative to some basis $x,y$, then $Q(x)=Q(y)=0$, so $x,y\in L'$, but
$Q(x+y)=2B(x,y)=\pi^R$ so $\ord Q(x+y)=R$ so $x+y\notin L'$. Thus $L'$
is not a lattice.

(iv) Since $2\sss L\subset\nnn L$, by (iii), $L'$ is a lattice. Let
$x\in L$ be a norm generator. We prove that $L=\oo x+L'$. Let
$x'\in L$. If $x'\in L'$ then $x'\in\oo x+L'$ trivially. If
$x'\in L\setminus L'$ then $\ord Q(x')=\ord Q(x)=R$ so
$Q(x')=\eta Q(x)$, with $\eta\in\ooo$. As seen in the beginning of
[OM1, \S63], we have $-\eta\equiv\delta^2\pmod\p$, so
$\eta +\delta^2\in\p$, for some  $\delta\in\ooo$. If $y=x'-\delta x$,
then $Q(y)=Q(x')+\delta^2Q(x)-2\delta B(x,x')
=(\eta +\delta^2)Q(x)-2\delta B(x,x')$. We have $\eta +\delta^2\in\p$,
so $\ord (\eta +\delta^2)Q(x)>\ord Q(x)=R$ and
$2\delta B(x,x')\in 2\sss L\subset\nnn L=\p^R$ so
$\ord 2\delta B(x,x')>R$. Thus $\ord Q(y)>R$, so $y\in L'$. It follows
that $x'=\delta x+y\in\oo x+L'$.

Consequently,
$L/L'=(\oo x+L')/L'\cong\oo x/(\oo x\cap L')=\oo x/\p x\cong\oo/\p$.
(If $\alpha\in\oo$, we have $\alpha x\in L'$ iff
$\ord Q(\alpha x)>R=\ord Q(x)$, which is equivalent to
$\alpha\in\p$. Thus $\oo x\cap L'=\p x$.) \qed

\blm If $M=\p^{-1}x_1+L$, then
$M\cong{\prec\pi^{-2}a_1,\ldots,a_n\succ}$  relative to the BONG
$\pi^{-1}x_1,\ldots,x_n$ and $M'=L$.
\elm
\pf Let $R=\ord a_1=\ord\nnn L$. We have
$\p^{-1}x_1\subseteq\p^{-1}L$ and $L\subseteq\p^{-1}L$ so
$M\subseteq\p^{-1}M$. It follows that
$\nnn M\subseteq\p^{-2}\nnn L=\p^{R-2}$. But $\pi^{-1}x_1\in M$ and
$\ord Q(\pi^{-1}x_1)=\ord\pi^{-2}a_1=R-2$, so $\nnn M=\p^{R-2}$ and
$\pi^{-1}x_1$ is a norm generator for $M$. Since
$pr_{\pi^{-1}x_1^\perp}M=pr_{x_1^\perp}(\p^{-1}x_1+L)=pr_{x_1^\perp}L
={\prec x_2,\ldots,x_n\succ}$, we have
$M={\prec\pi^{-1}x_1,x_2,\ldots,x_n\succ}$, as claimed.

To prove that $M'=L$, we note first that $\nnn L\subset\nnn M$, so
$L\subseteq M'$. Conversely, let $x\in M$, $x=\alpha x_1+y$, with
$\alpha\in\p^{-1}$ and $y\in K$. We must prove that if $x\in M'$,
i.e. if $\ord Q(x)>\ord\nnn M=R-2$, then $x\in L$. We have
$\alpha\in\p^{-1}$ so $\ord\alpha\geq -1$. Suppose that
$\ord\alpha =-1$. Then $\ord\alpha^2Q(x_1)=R-2$, $y\in K$, so
$\ord Q(y)\geq\ord\nnn K=R$, and
$2B(x_1,y)\in 2\sss L\subseteq\nnn L=\p^R$, so
$\ord 2\alpha B(x_1,y)=-1+\ord B(x_1,y)\geq R-1$. By the domination
principle, we have
$\ord Q(x)=\ord (\alpha^2Q(x_1)+Q(y)+2\alpha B(x_1,y))=R-2$.
Contradiction. Thus $\ord\alpha\geq 0$, i.e. $\alpha\in\oo$. Since
$x_1,y\in L$, we get $x=\alpha x_1+y\in L$. \qed

Upon replacing $L$, $M$ and $x_1$ by $K$, $L$ and $\pi x_1$ and noting
that $\p^{-1}\pi x_1=\oo x_1$, we get the following equivalent
statement.

\bco If there is a lattice
$K\cong{\prec\pi^2a_1,\ldots,a_n\succ}$ relative to the BONG
$\pi x_1,\ldots,x_n$, then $L=\oo x_1+K$ and $L'=K$.
\eco

{\bf Note.} In Lemma 4.4 and Corollary 4.5 the BONGs need not be good.
\medskip

In the next three lemmas $L$ is binary.

\blm $L$ is not $\nnn L$-maximal iff $\pi^{-2}a_2/a_1\in\mathcal A$
or, equivalently, iff $R_2-R_1\geq 2-2e$ and
$R_2-R_1\geq 2-d(-a_1a_2)$.
\elm
\pf First note that $\ord\pi^{-2}a_2/a_1=R_2-R_1-2$ and
$d(-\pi^{-2}a_2/a_1)=d(-a_2/a_1)$ so $\pi^{-2}a_2/a_1\in\mathcal A$
iff $R_2-R_1-2\geq -2e$ and $R_2-R_1-2\geq -d(-a_1a_2)$, as claimed.

For the ``if'' part, we note that
$Q(\pi^{-1}x_2)/Q(x_1)=\pi^{-2}a_2/a_1\in\mathcal A$ so the lattice
$M={\prec x_1,\pi^{-1}x_2\succ}$ exists. Since $x_1$ is a norm generator
for both $M$ and $L$ and
$pr_{x_1^\perp}M=\oo\pi^{-1}x_2\supset\oo x_2=pr_{x_1^\perp}L$, by
[B2, Lemma 2.2], we have $M\supset L$. Since $\nnn M=\nnn L$ and
$M\supset L$, $L$ is not $\nnn L$-maximal.

For the ``only if'' part, if $L$ is not $\nnn L$-maximal, then there
is a lattice $M$ with $FM=FL=:V$, $M\supset L$ and $\nnn M=\nnn L$. Since
$x_1\in L\subset M$ and $\nnn L=Q(x_1)\oo$, $x_1$ is also a norm
generator of $M$. We have $pr_{x_1^\perp}M\subset pr_{x_1^\perp}V=Fx_2$ so
$pr_{x_1^\perp}M=\p^kx_2$ for some $k\in\ZZ$. Then
$M={\prec x_1,\pi^kx_2\succ}$, so
$\pi^{2k}a_2/a_1=Q(\pi^kx_2)/Q(x_1)\in\mathcal A$. We cannot have
$k\geq 0$ since this would imply that
$pr_{x_1^\perp}M\subseteq\oo x_2=pr_{x_1^\perp}L$, so $M\subseteq L$.
Thus $k\leq -1$. Then, by [B2, Lemma 3.8],
$\pi^{2k}a_2/a_1\in\mathcal A$ implies
$\pi^{-2}a_2/a_1\in\mathcal A$. \qed

\bco The lattice $L$ is $\nnn L$-maximal iff one of the following
happens.

(1) $R_2-R_1=-2e$.

(2) $2-2e\leq R_2-R_1\leq 1$ and $d(-a_1a_2)=R_1-R_2+1$. Equivalently,
either $R_2-R_1=1$ or $R_2-R_1$ is even, $2-2e\leq R_2-R_1\leq 0$ and
$d(-a_1a_2)=R_1-R_2+1$.
\eco
\pf Since $L\cong{\prec a_1,a_2\succ}$, we have $a_2/a_1\in\mathcal A$
so $R_2-R_1\geq -2e$ and $R_2-R_1\geq -d(-a_1a_2)$. By Lemma 4.6,
the condition that $L$ is maximal is equivalent to
$\pi^{-2}a_2/a_1\notin\mathcal A$, which means that $R_2-R_1\geq 2-2e$
or $R_1-R_2\geq 2-d(-a_1a_2)$ fails.

If $R_2-R_1\geq 2-2e$ fails, then $-2e\leq R_2-R_1\leq 1-2e$. But, by
Lemma 3.5(i), $R_2-R_1$ cannot be odd and negative, so $R_2-R_1=-2e$. 

If $R_2-R_1\geq 2-2e$ holds, but $R_1-R_2\geq 2-d(-a_1a_2)$ fails,
then $-d(-a_1a_2)\leq R_2-R_1\leq 1-d(-a_1a_2)$. But, since
$R_2-R_1\neq -2e$, so $a_2/a_2\neq -\Delta/4$, by Lemma 1.4, the
inequality $R_2-R_1\geq -d(-a_1a_2)$ must be strict. Hence
$R_2-R_1=1-d(-a_1a_2)$, i.e. $d(-a_1a_2)=R_1-R_2+1$. Since
$d(-a_1a_2)<\infty$, we have $0\leq d(-a_1a_2)\leq 2e$, which implies
that $1-2e\leq R_2-R_1\leq 1$. But $R_2-R_1$ cannot be odd and
negative, so $2-2e\leq R_2-R_1\leq 1$ and if $2-2e\leq R_2-R_1\leq 0$,
then $R_2-R_1$ is even. If $R_2-R_1=1$ then $\ord a_1a_2=R_1+R_2$ is
odd, so $d(-a_1a_2)=0$, and so the condition that
$d(-a_1a_2)=R_1-R_2+1$ is superfluous. \qed

Recall that $R_2-R_1=-2e$ iff $a_2/a_1\in -\frac 14\oos$ or
$-\frac\Delta 4\oos$, in which case $L$ is isomorphic, up to scalling,
to $A(0,0)$ or $A(2,2\rho )$, respectively. When $L$ is hyperbolic, a
case handled by Lemma 4.3(iii), $L'$ is not a lattice. For the
remaining cases we have:

\blm Assume that $L$ is maximal and not hyperbolic. Let $V=FL$.

(i) $g(a_2/a_1)=\oos\cap{\rm N}(-a_1a_2)$.

(ii) $O(L)=O(V)$, so
$G(a_2/a_1)=\theta (O^+(L))=\theta (O^+(V))={\rm N}(-a_1a_2)$. 

(iii) If $R_2-R_1=-2e$, i.e. if $a_2/a_1\in -\frac\Delta 4\oos$, then
$L'=\p L\cong{\prec\pi^2a_1,\pi^2a_2\succ}$ relative to the BONG
$\pi x_1,\pi x_2$.

(iv) If $2-2e\leq R_2-R_1\leq 1$ and $d(-a_1a_2)=R_1-R_2+1$, then
there is $\eta\in\ooo$ such that $(\pi\eta,-a_1a_2)_\p =1$ and for
every such $\eta$ we have
$L'\cong L^{\pi\eta}\cong{\prec\pi\eta a_1,\pi\eta a_2\succ}$. 
\elm
\pf Note that $\det V\neq -1$, so $V$ is anisotropic and we can apply
Lemma 4.2.  

(i) If $\eta\in\ooo$ then $\eta\in g(a_2/a_1)$
iff $L\cong{\prec\eta a_1,\eta a_2\succ}\cong L^\eta$. A necessary
condition is that
$FL\cong FL^\eta$, i.e. $[a_1,a_2]\cong [\eta a_1,\eta a_2]$, which is
equivalent to $\eta\in{\rm N}(-a_1a_2)$. So
$g(a_2/a_1)\subseteq\oos\cap{\rm N}(-a_1a_2)$. Conversely, if
$\eta\in\ooo\cap{\rm N}(-a_1a_2)$ then $FL\cong FL^\eta$, so $V=FL$
contains a lattice $M\cong L^\eta$. Since $L$ is maximal and
$\eta\in\ooo$, $L^\eta$ too is maximal of the same norm. But by Lemma
4.2, there is only one $\nnn L$-maximal lattice on $V$ so
$L=M\cong L^\eta$. Thus
$[a_1,a_2]\cong [\eta a_1,\eta a_2]$ and so $\eta\in g(a_2/a_1)$.

(ii) We have $O(L)=O(V)$, by Lemma 4.2. Then $O^+(L)=O^+(V)$ and so
$\theta (O^+(L))=\theta (O^+(V))={\rm N}(-\det V)={\rm N}(-a_1a_2)$.

For (iii) and (iv) note that, since $L$ is $\nnn L$-maximal, by Lemma
4.2, we have $L=\{ x\in V\mid\, Q(x)\in\nnn L\}$. It follows that
$L'=\{ x\in V\mid\, Q(x)\in\p\nnn L\}$ so $L'$ is the only
$\p\nnn L$-maximal lattice on $V$.

(iii) Since $a_2/a_1\in -\Delta\fff$, we have
$V\cong [a_1,a_2]\cong [a_1,-\D a_1]$. It follows that $V$ only
represents elements $b\in\fff$ with
$\ord b\equiv\ord a_1=\ord\nnn L\pmod 2$. In particular, $V$ doesn't
reperesent elements of the same order as
$\p\nnn L$. Hence $L'=\{ x\in V\mid\, Q(x)\in\p\nnn L\}
=\{ x\in V\mid\, Q(x)\in\p^2\nnn L\}$ so $L'$ is the only
$\p^2\nnn L$-maximal lattice over $V$. But $L$ is $\nnn L$-maximal, so
$\p L$ is $\p^2\nnn L$-maximal. Thus $L'=\p L$.

(iv) Since $R_2-R_1\geq 2-2e$, we have
$d(-a_1a_2)=R_1-R_2+1<2e=d(\Delta )$ so $-a_1a_2\neq\Delta$ in
$\fff/\fff^2$. It follows that
${\rm N}(-a_1a_2)\not\subseteq{\rm N}(\Delta )=\ooo\fff^2$, so
${\rm N}(-a_1a_2)$ contains elements of odd orders. Let $\eta\in\ooo$
such that $\pi\eta\in{\rm N}(-a_1a_2)$, 
i.e. such that $(\pi\eta,-a_1a_2)_\p =1$. It follows that
$[a_1,a_2]\cong [\pi\eta a_1,\pi\eta a_2]$ so, since
$L^{\pi\eta}\cong{\prec\pi\eta a_1,\pi\eta a_2\succ}$, we have
$V=FL\cong FL^{\pi\eta}$. Let $M\cong L^{\pi\eta}$ be a lattice over
$V$. Since $L$ is $\nnn L$-maximal, $M$ is $\p\nnn L$-maximal and so
$L'=M\cong L^{\pi\eta}$. \qed

Note that if $L'$ is a lattice, then Lemma 4.3(iv) describes $L'$ in
all cases except when $\nnn L=2\sss L$, but $L$ doesn't split a
hyperbolic lattice of scale $\sss L$. We now handle this remaining
case.

\blm We have $\nnn L=2\sss L$ iff $R_2-R_1=-2e$.
\elm
\pf The condition $\nnn L=2\sss L$ writes as
$\ord\nnn L=\ord\sss L+e$. Since, by Lemma 1.5(i),
$\ord\nnn L=R_1$ and $\ord\sss L=\min\{ R_1,(R_1+R_2)/2\}$, this is
equivalent to $R_1=(R_1+R_2)/2+e$, i.e. to $R_2-R_1=-2e$. \qed

\blm Suppose that $\nnn L=\p^R$, i.e. $R_1=R$. Then the following are
equivalent.

(1) $\nnn L=2\sss L$ and $L$ doesn't spilt any hyperbolic lattice of
scale $\sss L$.

(2) $L=J\perp K$, where $J\cong\frac 12\pi^RA(2,2\rho )$ and
$\nnn K\subset\nnn L$ (or $K=0$, if $\rank L=2$).

(3) $R_2-R_1=-2e$, $R_1<R_3$ (or $n=2$) and $d(-a_1a_2)=2e$,
i.e. $a_1a_2\in -\Delta\fff$. 
\medskip

Moreover, if (1)-(3) hold, then $L'=\p J\perp K$ and
$L/L'\cong\oo/\p\oplus\oo/\p$. Also we can take
$J\cong{\prec a_1,a_2\succ}$ and
$K\cong{\prec a_3,\ldots,a_n\succ}$.
\elm
\pf Suppose that (2) holds. Then $\nnn L=\nnn J$ and
$\sss L\supseteq\sss J=\frac 12\nnn J=\frac 12\nnn L$. It follows that
$\sss L=\frac 12\nnn L$, i.e. $\nnn L=2\sss L$.

By Lemmas 4.3(i) and 4.8(iii), $L'=J'\perp K=\p J\perp K$. Then, since
$L'$ is a lattice, by Lemma 4.3(iii), $L$ cannot split a hyperbolic
lattice of scale $\sss L$, so (1) holds. So we proved
(2)$\Rightarrow$(1).

Suppose that (3) holds. We have $R=\ord\nnn L=R_1$. Since
$R_2-R_1=-2e$, by Lemma 4.9, we have $\nnn L=2\sss L$. If $n\geq 3$,
then $R_3>R_1>R_2$ and if $n\geq 4$, then $R_4-R_3\geq -2e$, so
$R_4\geq R_3-2e>R_1-2e=R_2$. Then, by Lemma 1.5(ii), we have
a splitting $L=J\perp K$, where $J\cong{\prec a_1,a_2\succ}$ and
$K\cong{\prec a_3,\ldots,a_n\succ}$. Since $R_2-R_1=-2e$ and
$a_1a_2\in -\Delta\fff^2$, we have $a_2/a_1\in -\frac\Delta 4\oos$ so,
by Lemma 1.4, we have
$J\cong\frac 12\pi^{R_1}A(2,2\rho )=\frac 12\pi^RA(2,2\rho )$. We also
have $\ord\nnn K=R_3>R_1=\ord\nnn L$ so (2) holds. Thus
(3)$\Rightarrow$(2).

Suppose now that (1) holds. Since $\sss L=\frac 12\nnn L=\frac 12\p^R$
and $L$  doesn't split $\frac 12\p^RA(0,0)$, $L$ doesn't represent
$\frac 12\p^RA(0,0)$. The condition $\nnn L=2\sss L$ is equivalent, by
Lemmma 4.9, to $R_2-R_1=-2e$. This implies that
${\prec a_1,a_2\succ}\cong\frac 12\pi^RA(0,0)$ or
$\frac 12\pi^RA(2,2\rho )$, corresponding to $a_1a_2\in -\fff^2$ or
$-\Delta\fff^2$, respectively. (Recall that $R=R_1$.) But $L$ doesn't
represent $\frac 12\pi^RA(0,0)$, so the latter holds. Assume that
$n\geq 3$ and $R_1=R_3$. Since $\ord a_3/a_2=R_3-R_2=2e$, by Lemma
2.19(ii), $\Delta\in g(a_3/a_2)$, which, by Lemma 1.5(iv), implies
$L\cong{\prec a_1,\Delta a_2,\Delta a_3,a_4,\ldots,a_n\succ}$ relative
to some good BONG. It follows that $L$ represents
${\prec a_1,\Delta a_2\succ}$. But $R_2-R_1=-2e$ and
$a_1a_2\in -\Delta\fff^2$ so $a_1\Delta a_2\in -\fff^2$ so
${\prec a_1,\Delta a_2\succ}\cong\frac 12\pi^RA(0,0)$. This
contradicts the hypothesis that $L$ doesn't represent
$\frac 12\pi^RA(0,0)$. Hence we must have $R_1<R_3$ (or $n=2$) and so
(3) holds. Thus (1)$\Rightarrow$(3). \qed

We now give an explicit description of $L'$ in terms of good BONGs in
the cases when $\theta (O^+(L))\neq\fff$. In [B2] we defined
properties A and B of a quadratic lattice. By Proposition 1.6(ii) and
(iii), if a quadratic lattice $L$ doesn't have property A then
$\ooo\fff^2\subseteq\theta (O^+(L))$ and if $L$ has property A, but
not property B, then $\theta (O^+(L))=\fff$. So if
$\theta (O^+(L))\neq\fff$ then $L$ has property B or
$\theta (O^+(L))=\ooo\fff^2$. (We have $|\fff/\ooo\fff^2|=2$ so there
is no intermediate group between $\ooo\fff^2$ and $\fff$.)

\blm Assume that $R_k=R_l=:R$ for some $1\leq k<l\leq n$, with
$k\equiv l\pmod 2$. 

(i) $R_k=R_{k+2}=\cdots =R_l=R$,
$R_k\equiv R_{k+1}\equiv\cdots\equiv R_l=R\pmod 2$ and
$a_{k,l}\in\pi^R\ooo\fff^2$.

(ii) If $G(a_{i+1}/a_i)\subseteq\ooo\fff^2$ for
$k\leq i\leq l-1$ and $\varepsilon\in\ooo$, such that
$a_{k,l}\in (-1)^{(l-k)/2}\pi^R\varepsilon\fff^2$ then
$L=H_1\perp\cdots\perp H_{(l-k)/2}\perp K$, where
$H_t\cong\pi^{(R+R_{k+2t-1})/2}A(0,0)$ and $K\cong
{\prec a_1,\ldots,a_{k-1},\pi^R\varepsilon,a_{l+1},\ldots,a_n\succ}$.
\elm
\pf (i) We have $R=R_k\leq R_{k+2}\leq\cdots\leq R_l=R$, so
$R_k=R_{k+2}=\cdots =R_l=R$. In particular, $R_i\equiv R\pmod 2$ for
$k\leq i\leq l$, $i\equiv k\pmod 2$. If $k<i<l$, $i\equiv k+1\pmod 2$,
then $R_{i-1}=R_{i+1}=R$, so, by Lemma 3.5(ii),
$R_i\equiv R_{i-1}=R\pmod 2$. 

We write $a_i=\pi^{R_i}\varepsilon_i$, with
$\varepsilon_i\in\ooo$. Then
$a_{k,l}=\pi^{R_k+\cdots +R_l}\varepsilon_{k,l}$. But
$R_k+\cdots +R_l$ is a sum of odd number of integers which are
$\equiv R\pmod 2$, so $R_k+\cdots +R_l\equiv R\pmod 2$. Hence in
$\fff/\ffs$ we have
$a_{k,l}=\pi^R\varepsilon_{k,l}\in\pi^R\ooo\fff^2$.

(ii) In $\fff/\fff^2$ we have $a_{k,l}=\pi^R\varepsilon_{k,l}$, so
$a_{k,l}=(-1)^{(l-k)/2}\pi^R\varepsilon$ iff
$\varepsilon =(-1)^{(l-k)/2}\varepsilon_{k,l}$ in $\ooo/\oos$.

We prove our result by induction on $l-k$. When $l-k=2$ this is just
[B2, Lemma 7.3]. (In this case $\varepsilon =-\varepsilon_{k,k+2}$ and
$H_1\cong\pi^{(R+R_{k+1})/2}A(0,0)=\pi^{(R_k+R_{k+1})/2}A(0,0)$.)

Suppose now that $l-k\geq 4$. From the case $l-k=2$ we get
$L=H_1\perp M$, with $H_1\cong\pi^{(R+R_{k+1})/2}A(0,0)$ and
$M\cong{\prec a_1,\ldots,a_{k-1},a'_{k+2},a_{k+3},\ldots,a_n\succ}$,
where $a'_{k+2}=\pi^R\xi =-a_{k,k+2}$ in $\fff/\fff^2$, i.e.
$\xi =-\varepsilon_{k,k+2}$ in $\ooo/\oos$. We have $k+2<l$,
$k+2\equiv l\pmod 2$, $\ord a'_{k+2}=R=R_{k+2}=R_l$ and
$G(a_{i+1}/a_i)\subseteq\ooo\fff^2$ for $k+3\leq i\leq l-1$. In order
to use the induction hypothesis, we still have to prove that
$G(a_{k+3}/a'_{k+2})\subseteq\ooo\fff^2$. Since in $\fff/\oos$ we have
$a_{k+3}=\pi^{R_{k+3}}\varepsilon_{k+3}$ and
$a'_{k+2}=-\pi^R\varepsilon_{k,k+2}$, this is equivalent to
$G(-\pi^{R_{k+3}-R}\varepsilon_{k,k+3})\subseteq\ooo\fff^2$. By
[B2, Lemma 7.2(ii)], this follows from $G(a_{k+1}/a_k)=
G(\pi^{R_{k+1}-R_k}\varepsilon_{k,k+1})\subseteq\ooo\fff^2$,
$G(a_{k+3}/a_{k+2})=
G(\pi^{R_{k+3}-R_{k+2}}\varepsilon_{k+2,k+3})\subseteq\ooo\fff^2$,
$R_{k+3}-R=R_{k+3}-R_k\geq R_{k+1}-R_k$ and
$R_{k+3}-R=R_{k+3}-R_{k+2}$.

By the induction hypothesis,
$M=H_2\perp\cdots H_{(l-k)/2}\perp K$. Here $H_2,\ldots,H_{(l-k)/2}$
are hyperbolic lattices with scales of orders
$(R+R_{k+3})/2,(R+R_{k+5})/2,\ldots,(R+R_{l-1})/2$,
i.e. $H_t\cong\pi^{(R+R_{k+2t-1})/2}A(0,0)$ for $2\leq t\leq (l-k)/2$
and $K\cong
{\prec a_1,\ldots,a_{k-1},\pi^R\varepsilon,a_{j+1},\ldots,a_n\succ}$,
with $\varepsilon\in\ooo$ such that in $\fff/\fff^2$ we have
$\pi^R\varepsilon =(-1)^{(l-(k+2))/2}a'_{k+2}a_{k+3}\cdots a_l$. But
in $\fff/\fff^2$ we also have $a'_{k+2}=-a_{k,k+2}$
and so $\pi^R\varepsilon =(-1)^{(k-l)/2}a_{l,k}$, as
claimed. Hence we get the desired relation,
$L=H_1\perp M=H_1\perp\cdots\perp H_{(k-l)/2}\perp K$. \qed

\bco (i) Let $1\leq k<l\leq n$, with $k\equiv l\pmod 2$ and let
$\varepsilon\in\ooo$. Assume that $K\cong{\prec
a_1,\ldots,a_{k-1},\pi^R\varepsilon,a_{l+1},\ldots,a_n\succ}$ relative
to some good BONG, with $\ord a_i=R_i$. Let $R_k,\ldots,R_l\in\ZZ$
with $R_k=R_{k+2}=\cdots =R_l=R$,
$R_k\equiv R_{k+1}\equiv\cdots\equiv R_l\equiv R\pmod 2$,
$R-2e\leq R_{k+1}\leq R_{k+3}\leq\cdots\leq R_{l-1}\leq R+2e$,
$R_{k-1}\leq R_{k+1}$ (if $k>1$) and $R_{l-1}\leq R_{l+1}$ (if
$l<n$).

If $L=H_1\perp\cdots\perp H_t\perp K$, where
$H_t\cong\pi^{(R+R_{k+2t-1})/2}A(0,0)$ for $1\leq t\leq (l-k)/2$, then
$L\cong{\prec a_1,\ldots,a_{k-1},\pi^{R_k}\varepsilon,
-\pi^{R_{k+1}}\varepsilon,\pi^{R_{k+2}}\varepsilon,
\ldots,-\pi^{R_{l-1}}\varepsilon,\pi^{R_l}\varepsilon,
a_{l+1},\ldots,a_n\succ}$.

(ii) If $L$ is the lattice from Lemma 4.11(ii), then

$L\cong{\prec a_1,\ldots,a_{k-1},\pi^{R_k}\varepsilon,
-\pi^{R_{k+1}}\varepsilon,\pi^{R_{k+2}}\varepsilon,\ldots, 
-\pi^{R_{l-1}}\varepsilon,\pi^{R_l}\varepsilon,a_{l+1},\ldots,a_n\succ}$. 
\eco 
\pf (i) We use Lemma 1.4 to prove that there is a lattice
$L\cong{\prec a_1,\ldots,a_n\succ}$, where
$a_i=(-1)^{t-i}\pi^{R_i}\varepsilon$ for $k\leq t\leq l$ and then we
use Lemma 4.11 to prove that
$L\cong H_1\perp\cdots\perp H_{(l-k)/2}\perp K$.

We have $\ord a_i=R_i$ holds for every $1\leq i\leq n$, so we must
verify that $R_i\leq R_{i+2}$ for $1\leq i\leq n-2$ and
$a_{i+1}/a_i\in\mathcal A$ for $1\leq i\leq n-1$. From the
corresponding condition for $K$, $R_i\leq R_{i+2}$ holds for
$1\leq i\leq k-3$ and for $l+1\leq i\leq n-2$ and we also have
$R_{k-2}\leq R=R_k$ (if $k\geq 3$) and $R_l=R\leq R_{l+2}$ (if
$l\leq n-2$). We are left to prove that $R_i\leq R_{i+2}$
when $k-1\leq i\leq l-1$. But this follows from the hypothesis that
$R_{k-1}\leq R_{k+1}\leq R_{k+3}\leq\cdots\leq R_{l-1}\leq R_{l+1}$
and $R_k=R_{k+2}=\cdots =R_l$. From the conditions for $K$, we have
$a_{i+1}/a_i\in\mathcal A$ for $1\leq i\leq k-2$ and for . From the
corresponding condition for $K$, we get that
$a_{i+1}/a_i\in\mathcal A$ holds for $1\leq i\leq k-2$ and for
$l+1\leq i\leq n-1$ and we also have
$a_k/a_{k-1}=\pi^R\varepsilon/a_{k-1}\in\mathcal A$ and
$a_{l+1}/a_l=a_{l+1}/\pi^R\varepsilon\in\mathcal A$. So we are left to
proving that $a_{i+1}/a_i\in\mathcal A$ for $k\leq i\leq l-1$. At
these indices we have $a_{i+1}/a_i=-\pi^{R_{i+1}-R_i}$. Also
$R_{i+1}-R_i$ is even, so $\pi^{R_{i+1}-R_i}\in\fff^2$ and so 
$R_{i+1}-R_i+d(\pi^{R_{i+1}-R_i})=\infty >0$. We still have to prove
that $R_{i+1}-R_i\geq -2e$. By hypothesis, for every, $k\leq j\leq k$,
we have $R_j=R$ if $j\equiv k\pmod 2$ and $R-2e\leq R_j\leq R+2e$ if
$j\equiv k+1\pmod 2$. Thus $R_{i+1}-R_i=R_{i+1}-R\geq -2e$ if
$i\equiv k\pmod 2$ and $R_{i+1}-R_i=R-R_i\geq -2e$ if
$i\equiv k+1\pmod 2$. This concludes the proof.

Hence the lattice $L\cong{\prec a_1,\ldots,a_n\succ}$ exists. We have
$k<l$, $R_k=R_l=R$ and for $k\leq i\leq l-1$ we have that
$R_{i+1}-R_i$ is even and $d(1)=\infty >e-(R_{i+1}-R_i)/2$ so, by
[B2, Lemma 7.2(i)], we have
$G(a_{i+1}/a_i)=G(-\pi^{R_{i+1}-R_i})\subseteq\ooo\fff^2$. If
$k\leq i\leq l$, then $R_i\equiv R\pmod 2$ so in $\fff/\fff^2$ we have
$a_i=(-1)^{i-k}\pi^{R_i}\varepsilon_i$, with
$\varepsilon_i=(-1)^{i-k}\varepsilon$. Same as in the proof of Lemma
4.11, we have $R_k+\cdots +R_l\equiv R\pmod 2$ so
$a_{k,l}=\pi^R\varepsilon_{k,l}$ in $\fff/]\ffs$. Since the sequence
$\varepsilon_k,\ldots,\varepsilon_l$ is
$\varepsilon,-\varepsilon,\ldots,-\varepsilon,\varepsilon$, with
$(k-l)/2+1$ copies of $\varepsilon$ and $(k-l)/2$ copies of
$-\varepsilon$, in $\fff/\ffs$ we have
$\varepsilon_{k,l}=(-1)^{(l-k)/2}\varepsilon$, so
$a_{k,l}=(-1)^{(l-k)/2}\pi^R\varepsilon$. Then, by Lemma 4.11(ii), we
have $L\cong H_1\perp\cdots\perp H_{(l-k)/2}\perp K$.

(ii) If $L$ is like in Lemma 4.11(ii) then we have a splitting
$L=H_1\perp\cdots\perp H_{(l-k)/2}\perp K$, with $K$ and $H_t$ as in
(i). To prove the claimed result, we have to verify that
$R_k,\ldots,R_l$ satisfy the conditions of (i). The conditions
$R_k=R_{k+2}=\cdots =R_l=R$ and
$R_k\equiv R_{k+1}\equiv\cdots\equiv R_l=R\pmod 2$ follow from Lemma
4.11(i). Since $R_i=R_i(L)$ we have
$R_{k-1}\leq R_{k+1}\leq R_{k+3}\leq\cdots\leq R_{l-1}\leq R_{l+1}$
and, since $R_{k+1}-R_k\geq -2e$ and $R_l-R_{l-1}\geq -2e$, we also
have $R-2e=R_k-2e\leq R_{k+1}$ and $R_{l-1}\leq R_l+2e=R+2e$. Hence
all the conditions of (i) are fulfilled. \qed

\blm Assume that $n\geq 2$ and $L$ has property B or
$\theta (O^+(L))=\ooo\fff^2$.

(i) If $1\leq i\leq n-2$ and $R_i<R_{i+2}$ then $R_{i+2}-R_i\geq 2$.

(ii) If $1\leq i\leq n-1$, $R_{i+1}-R_i>-2e$ and
${\prec a_i,a_{i+1}\succ}$ is maximal, then
$\theta (O^+(L))\not\subseteq\ooo\fff^2$, so $L$ has property B.
\elm
\pf (i) If $L$ has property B, this is just [B2, Remark 4.12]. If
$\theta (O^+(L))=\ooo\fff^2$, this follows from the fact that
$R_1,\ldots,R_n$ have the same parity. (See the remark following
Theorem 1.9.)

By Corollary 4.7, $2-2e\leq R_{i+1}-R_i\leq 1$ and
$d(-a_ia_{i+1})=R_i-R_{i+1}+1\leq 2e-1$ and, by Lemma 4.8(ii),
$\theta (O^+(L))\supseteq\theta (O^+({\prec a_i,a_{i+1}\succ}))
={\rm N}(-a_ia_{i+1})$. Since $d(-a_ia_{i+1})<2e=d(\Delta )$, we have
$-a_ia_{i+1}\neq\Delta$ in $\fff/\fff^2$, so
${\rm N}(-a_ia_{i+1})\not\subseteq{\rm N}(\Delta )=\ooo\fff^2$. Thus
$\theta (O^+(L))\not\subseteq\ooo\fff^2$. By hypothesis, this implies
that $L$ has property B. \qed

\blm Let $L$ be like in Lemma 4.13.

If $R_2-R_1=-2e$ then:

\noindent (i) If $a_1a_2\in -\Delta\fff^2$ and $R_1<R_3$ (or $n=2$),
then $L'\cong{\prec\pi^2a_1,\pi^2a_2,a_3,\ldots,a_n\succ}$.

\noindent (ii) If $a_1a_2\in -\fff^2$ or $R_1=R_3$, then $L'$ is not a
lattice.

If $R_2-R_1>-2e$ or $n=1$ then:

\noindent (iii) If ${\prec a_1,a_2\succ}$ is not maximal and $R_1<R_3$
(or $n=2$) then $L'\cong{\prec\pi^2a_1,a_2,a_3,\ldots,a_n\succ}$. Same
holds if $n=1$.

\noindent (iv) If ${\prec a_1,a_2\succ}$ is maximal then there is
$\eta\in\ooo$ such that $(\pi\eta,-a_1a_2)_p =1$ and for every such
$\eta$ we have
$L'\cong{\prec\pi\eta a_1,\pi\eta a_2,a_3,\ldots,a_n\succ}$.

\noindent (v) If $R_1=R_3$ and $1\leq k\leq\frac{n-1}2$ is maximal
with the property that $R_1=R_{2k+1}$, then
$a_{1,2k+1}\in (-1)^k\varepsilon\pi^{R_1}\fff^2$  for some
$\varepsilon\in\ooo$ and we have
$$L\cong{\prec\pi^{R_1}\varepsilon, -\pi^{R_2}\varepsilon,
\pi^{R_3}\varepsilon, \ldots,-\pi^{R_{2k}}\varepsilon,
\pi^{R_{2k+1}}\varepsilon,a_{2k+2},\ldots,a_n\succ}$$
and
$$L'\cong{\prec\pi^{R_1+2}\varepsilon, -\pi^{R_2-2}\varepsilon,
\pi^{R_3+2}\varepsilon, \ldots,-\pi^{R_{2k}-2}\varepsilon,
\pi^{R_{2k+1}+2}\varepsilon,a_{2k+2},\ldots,a_n\succ}.$$
\elm
\pf (i) By Lemma 4.10, we have $L'=\p J\perp K$, with
$J\cong{\prec a_1,a_2\succ}$ and $K\cong{\prec a_3,\ldots,a_n\succ}$. We
have $\p J\cong{\prec\pi^2a_1,\pi^2a_2\succ}$. So to prove that
$L'=\p J\perp K\cong{\prec\pi^2a_1,\pi^2a_2,a_3,\ldots,a_n\succ}$, by
Lemma 1.5(ii), we have to show that $\ord\pi^2a_1\leq\ord a_3$,
$\ord\pi^2a_2\leq\ord a_3$ (if $n\geq 3$) and
$\ord\pi^2a_2\leq\ord a_4$ (if $n\geq 4$), i.e. $R_1+2\leq R_3$,
$R_2+2\leq R_3$ and $R_2+2\leq R_4$. By Lemma 4.13(i) $R_1<R_3$
implies $R_3-R_1\geq 2$, i.e. $R_1+2\leq R_3$. Since
$R_2<R_1$ we also have $R_2+2<R_3$. And since
$R_4-R_3\geq -2e=R_2-R_1$, we have $R_4-R_2\geq R_3-R_1\geq 2$, so
$R_2+2\leq R_4$. 

(ii) follows from Lemmas 4.3(iii), 4.9 and 4.10.

(iii) By Corollary 4.5, it is enough to prove that a lattice
$K\cong{\prec\pi^2a_1,a_2,\ldots,a_n\succ}$ exists. Most of the
conditions of Lemma 1.4 follow from the fact that
$L\cong{\prec a_1,a_2,\ldots,a_n\succ}$ relative to the good BONG. The
remaining conditions, $a_2/(\pi^2a_1)=\pi^{-2}a_2/a_1\in\mathcal A$ and
$\ord\pi^2a_1=R_1+2\leq R_3=\ord a_3$, follow from Lemma 4.6 and Lemma
4.13(i). (We have $R_1<R_3$, so $R_3-R_1\geq 2$.)

(iv) By Lemma 4.13(ii), $L$ must have property B so $R_1+2\leq R_3$ and
$R_2+2\leq R_4$. By Corollary 4.7, $R_2-R_1\leq 1$ so
$R_2+1\leq R_1+2\leq R_3$. By Lemma 1.5(ii), since $R_2<R_3$,
we have $L=J\perp K$, where $J\cong{\prec a_1,a_2\succ}$ and
$K\cong{\prec a_3,\ldots,a_n\succ}$. We have
$\nnn K=\p^{R_3}\subset\p^{R_1}=\nnn L$ so, by Lemma 4.3(i), we have
$L'=J'\perp K$. By Lemma 4.8(iv),
$J'\cong{\prec\pi\eta a_1,\pi\eta a_2\succ}$. Then
$L'=J'\perp K\cong{\prec\pi\eta a_1,\pi\eta a_2,\ldots,a_n\succ}$
follows, by Lemma 1.5(ii), from $\ord\pi\eta a_1=R_1+1<R_3=\ord a_3$,
$\ord\pi\eta a_2=R_2+1\leq R_3=\ord a_3$ and
$\ord\pi\eta a_2=R_2+1<R_4=\ord a_4$. 

(v) We will apply Lemma 4.11 and Corollary 4.12 when the pair $(k,l)$
is $(1,2k+1)$. Since $R_1=R_3$, $L$ cannot have property B so we have
$\theta (O^+(L))=\ooo\fff^2$. By Theorem 1.9,
$G(a_{i+1}/a_i)\subseteq\ooo\fff^2$ for $1\leq i\leq n-1$. Also, since
$R_1=R_3$, we have $(R_2-R_1)/2\equiv e\pmod 2$, i.e.
$R_2-R_1\equiv 2e\pmod 4$. Hence $R_2-R_1>-2e$ implies that
$R_2-R_1\geq 4-2e$.

By Lemma 4.11, since $R_1=R_{2k+1}=:R$ and
$G(a_{i+1}/a_i)\subseteq\ooo\fff^2$ for $1\leq i\leq 2k$, we have
$L=H_1\perp\cdots\perp H_k\perp K$, where
$H_t=\pi^{(R+R_{2t})/2}A(0,0)$ and
$K\cong{\prec\pi^R\varepsilon,a_{2k+2},\ldots,a_n\succ}$, with
$\varepsilon\in\ooo$ such that
$a_{1,2k+1}\in (-1)^k\pi^R\varepsilon\fff^2$. Also, by Corollary
4.12(ii), we get the stated formula for $L$.

We have $\theta (O^+(K))\subseteq\theta (O^+(L))=\ooo\fff^2$. Since
$\theta (O^+(K))\neq\fff$, $K$ has property B or
$\theta (O^+(K))=\ooo\fff^2$. If $n\geq 2k+2$, then
$R_{2k+2}-R\geq R_2-R=R_2-R_1>-2e$, i.e.
$\ord a_{2k+2}-\ord\pi^R\varepsilon >-2e$, so $K$ is not in one of the
cases (i) and (ii) of our lemma. By the maximality of $k$, if
$2k+3\leq n$ then $R_{2k+3}>R_1=R$, i.e.
$\ord a_{2k+3}>\ord\pi^R\varepsilon$, so $K$ is not in the case (v).
Finally, since $\ord a_{2k+2}-\ord\pi^R\varepsilon >-2e$ and 
$\theta (O^+(K))\subseteq\ooo\fff^2$, by Lemma 4.13(ii),
${\prec\pi^R\varepsilon,a_{2k+2}\succ}$ cannot be maximal. Hence $K$ is
in the case (iii) of our lemma and we have
$K'\cong{\prec\pi^{R'}\varepsilon,a_{2k+2},\ldots,a_n\succ}$, where
$R'=R+2$.

We denote by $R'_1=R'_3=\cdots =R'_{2k+1}=R'=R+2$ and by
$R'_{2t}=R_{2t}-2$ for $1\leq t\leq k$. We use Corollary 4.12(i) to
prove that $L'=H_1\perp\cdots\perp H_k\perp K'
\cong{\prec a'_1,\ldots,a'_{2k+1},a_{2k+3}\ldots,a_n\succ}$, where
$a'_i=(-1)^{i-1}\pi^{R'_i}\varepsilon$ for $1\leq i\leq 2k+1$. We have
$R'_1=R'_3=\cdots =R'_{2k+1}=R'$ by definition. The conditions
$R'_1\equiv R'_2\equiv\cdots\equiv R'_{2k+1}\equiv R'\pmod 2$ and
$R'_2\leq R'_4\leq\cdots\leq R'_{2k}$ follow from
$R_1\equiv R_2\equiv\cdots\equiv R_{2k+1}\equiv R\pmod 2$ and
$R_2\leq R_4\leq\cdots\leq R_{2k}$. It remains to prove that
$R'_2\geq R'-2e$, $R'_{2k}\leq R'+2e$ and, if $n\geq 2k+2$,
$R'_{2k}\leq R_{2k+2}$. Since $R_2-R_1\geq 4-2e$, we have
$R'_2=R_2-2\geq R_1+2-2e=R'-2e$. We have $R-R_{2k}=R_{2k+1}-R_{2k}\geq
-2e$ so $R_{2k}\leq R+2e$ and also $R_{2k}\leq R_{2k+2}$. Since
$R'_{2k}=R_{2k}-2<R_{2k}$ and $R=R'-2<R'$, these imply $R_{2k}'<R'+2e$
and $R'_{2k}<R_{2k+2}$. 

Hence, by Corollary 4.12(i), there is a lattice
$M\cong{\prec a'_1,\ldots,a'_{2k+1},a_{2k+2}\ldots,a_n\succ}$ and we
have $M\cong\bar H_1\perp\cdots\perp\bar H_k\perp K'$, where
$\bar H_t=\pi^{(R'+R'_{2t})/2}A(0,0)$. But for $1\leq t\leq k$ we have
$R'+R'_{2t}=(R+2)+(R_{2t}-2)=R+R_{2t}$ so
$\bar H_t=\pi^{(R+R_{2t})/2}A(0,0)\cong H_t$. Thus
$M\cong H_1\perp\cdots\perp H_t\perp K'\cong L'$, as claimed. \qed

\section{BONGs of lattices with $\theta (O^+(L))\neq\fff$}

Throughot this section $L\cong{\prec a_1,\ldots,a_n\succ}$ relative to a
good BONG $x_1,\ldots,x_n$ and $R_i(L)=R_i$. We will assume that
either $L$ has property B (relative to the BONG $x_1,\ldots,x_n$) or
$\theta (O^+(L))\subseteq\ooo\ffs$. By Proposition 1.6, if $L$ has
neither of these two properties, then $\theta (O^+(L))=\fff$.

\blm Let $a\in\mathcal A$ such that $\ord a=R\leq 2e$.

(i) If $G(a)\subseteq\ooo\ffs$, then $g(a)=\upo{R/2+e}$.

(ii) If $G(a)\not\subseteq\ooo\ffs$, then $g(a)=\upon{R+d(-a)}{-a}$.
\elm
\pf By [B2, Lemma 7.2], we have $G(a)\subseteq\ooo\ffs$ iff
$d(-a)>e-R/2$ or $a=-\Delta/4$ (in $\fff/\oos$). Hence, except for the
case $a=-\Delta/4$, both (i) and (ii) follow from the definition of
$g$ ([B2, Definition 6]). In the exceptional case we have $R=-2e$, so
$\upo{R/2+e}=\upo 0=\ooo$, but $d(-a)=d(\Delta )=2e=e-R/2$, so, by
definition,
$g(a)=\upon{R+d(-a)}{-a}=\upon 0\Delta =\ooo\cap\N (\Delta )=\ooo$.
This concludes the proof. \qed

\blm (i) For every $1\leq i\leq n-1$ we have
$$\alpha_i=\alpha_1({\prec a_i,a_{i+1}\succ} )=
\begin{cases}
(R_{i+1}-R_i)/2+e&\text{if }d(-a_{i,i+1})>e-(R_{i+1}-R_i)/2\\
R_{i+1}-R_i+d(-a_{i,i+1})&\text{if }d(-a_{i,i+1})\leq
e-(R_{i+1}-R_i)/2
\end{cases}.$$

(ii) If $\theta (O^+(L))\subseteq\ooo\fff^2$ then
$\alpha_i=(R_{i+1}-R_i)/2+e$ $\forall 1\leq i\leq n-1$.

(iii) For every $1\leq i\leq n-1$ such that
$\alpha_i<(R_{i+1}-R_i)/2+e$ we have $R_i-R_{i-1}>2e$ (if $i>1$) and
$R_{i+2}-R_{i+1}>2e$ (if $i<n-1$).
\elm
\pf (i) We have $\alpha_1({\prec a_i,a_{i+1}\succ} )
=\min\{ (R_{i+1}-R_i)/2+e,R_{i+1}-R_i+d(-a_{i,i+1})\}$, which is equal
to $(R_{i+1}-R_i)/2+e$ or $R_{i+1}-R_i+d(-a_{i,i+1})$, when
$d(-a_{i,i+1})\geq (R_{i+1}-R_i)/2$ or
$d(-a_{i,i+1})\leq (R_{i+1}-R_i)/2$, accordingly.

Note that $(R_{i+1}-R_i)/2+e$ and $R_{i+1}-R_i+d(-a_{i,i+1})$ are
elements of the set whose minimum is $\alpha_i$ so is suffices to
prove that all the other elements of that set are
$\geq (R_{i+1}-R_i)/2+e$.

First we take elements of the type $t=R_{i+1}-R_j+d(-a_{j,j+1})$ with
$j<i$. Since $j+1\leq i$, we have
$R_j+R_{j+1}\leq R_{j+1}+R_{j+2}\leq R_i+R_{i+1}$. If
$d(-a_{j,j+1})\geq e-(R_{j+1}-R_j)/2$ then
$t\geq R_{i+1}-(R_j+R_{j+1})/2+e\geq
R_{i+1}-(R_i+R_{i+1})/2+e=(R_{i+1}-R_i)/2+e$. If
$d(-a_{j,j+1})<e-(R_{j+1}-R_j)/2$ then, by [B2, Lemma 7.2(i)], we have
$G(a_{j+1}/a_j)\not\subseteq\ooo\ffs$, so
$\theta (O^+(L))\not\subseteq\ooo\fff^2$. Thus $L$ must have property
B, so $d(-a_{j,j+1})<e-(R_{j+1}-R_j)/2$ implies
$R_{j+2}-R_{j+1}>2e$. Together with $R_{j+1}-R_j+d(-a_{j,j+1})\geq 0$
and $R_{j+1}+R_{j+2}\leq R_i+R_{i+1}$ this implies that
$t\geq R_{i+1}-R_{j+1}>R_{i+1}-(R_{j+1}+R_{j+2})/2+e\geq
R_{i+1}-(R_i+R_{i+1})/2+e=(R_{i+1}-R_i)/2+e$.

Let now $t=R_{j+1}-R_i+d(-a_{j,j+1})$ with $j>i$. Since $j-1\geq i$,
we have $R_j+R_{j+1}\geq R_{j-1}+R_j\geq R_i+R_{i+1}$. If
$d(-a_{j,j+1})\geq e-(R_{j+1}-R_j)/2$, then
$t\geq (R_j+R_{j+1})/2-R_i+e\geq (R_i+R_{i+1})/2-R_i+e=(R_{i+1}-R_i)/2+e$.
If $d(-a_{j,j+1})<e-(R_{j+1}-R_j)/2$, same as in the previous case,
$L$ has property B, so $d(-a_{j,j+1})<e-(R_{j+1}-R_j)/2$ implies
$R_j-R_{j-1}>2e$. Together with $R_{j+1}-R_j+d(-a_{j,j+1})\geq 0$ and
$R_{j-1}+R_j\geq R_i+R_{i+1}$, this implies that
$t\geq R_j-R_i>(R_{j-1}+R_j)/2-R_i+e\geq
(R_i+R_{i+1})/2-R_i+e=(R_{i+1}-R_i)/2+e$.

(ii) If $\theta (O^+(L))\subseteq\ooo\fff$, then, by
[B2, Lemma 7.2(i)], $d(-a_{i,i+1})\geq e-(R_{i+1}-R_i)/2$ so, by (i),
we have $\alpha_i=(R_{i+1}-R_i)/2+e$.

(iii) If $\alpha_i<(R_{i+1}-R_i)/2+e$, then, by (ii),
$\theta (O^+(L))\not\subseteq\ooo\ffs$, so $L$ must have property
B. By (i), we have $d(-a_{i,i+1})<e-(R_{i+1}-R_i)/2$, so property B
implies that $R_i-R_{i-1}>2e$ and $R_{i+2}-R_{i+1}>2e$. \qed

\blm (i) Property B is independent of the choice of the BONG.

(ii) For $1\leq i\leq n-1$ the group $g(a_{i+1}/a_i)$ is independent
of the choice of the BONG.

(iii) If $L\cong{\prec b_1,\ldots,b_n\succ}$ relative to another good
BONG, then $b_{1,i}/a_{1,i}\in g(a_{i+1}/a_i)$ for all
$1\leq i\leq n-1$.   
\elm
\pf Assume that $L\cong{\prec b_1,\ldots,b_n\succ}$ relative to anuother
good BONG, $y_1,\ldots,y_n$.

(i) We must prove that if $L$ has property B relative to the BONG
$x_1,\ldots,x_n$, it also has property B relative to the BONG
$y_1,\ldots,y_n$. Since the numbers $R_i$ are invariants, it suffices
to prove that $d(-a_{i,i+1})\leq e-(R_{i+1}-R_i)/2$ for every
$1\leq i\leq n-1$ such that $d(-b_{i,i+1})\leq e-(R_{i+1}-R_i)/2$. 

Fist note that if $d(-b_{i,i+1})\leq e-(R_{i+1}-R_i)/2$, then
$0\leq e-(R_{i+1}-R_i)/2$, so $R_{i+1}-R_i\leq 2e$. By Lemma 5.2(iii),
this implies $\alpha_{i-1}=(R_i-R_{i-1})/2+e>e-(R_{i+1}-R_i)/2$ (if
$i>1$) and $\alpha_{i+1}=(R_{i+2}-R_{i+1})/2+e>e-(R_{i+1}-R_i)/2$ (if
$i<n-1$). (Since $L$ has property B, it has property A, so
$R_{i-1}<R_{i+1}$ and $R_i<R_{i+2}$.) Then, by Theorem 1.12(iii),
$d(a_{1,i-1}b_{1,i-1})\geq\alpha_{i-1}>e-(R_{i+1}-R_i)/2$ and 
$d(a_{1,i+1}b_{1,i+1})\geq\alpha_{i+1}>e-(R_{i+1}-R_i)/2$. (The two
inequalities hold also when $i=1$ and $i=n-1$, respectively. We have
$a_{1,0}=b_{1,0}=1$ and $a_{1,n}=b_{1,n}=\det FL$ in $\fff/\ffs$, so
$d(a_{1,0}b_{1,0})=d(a_{1,n}b_{1,n})=\infty$.) Since
$d(a_{1,i-1}b_{1,i-1}),d(a_{1,i+1}b_{1,i+1})>e-(R_{i+1}-R_i)/2\geq
d(-b_{i,i+1})$, by the domination principle, we have
$d(-a_{i,i+1})=d(-b_{i,i+1})\leq e-(R_{i+1}-R_i)/2$.

(ii) and (iii) We must prove that $g(a_{i+1}/a_i)=g(b_{i+1}/b_i)$ and
$b_{1,i}/a_{1,i}\in g(a_{i+1}/a_i)$ for all $1\leq i\leq n-1$. We
start by noting that, by Theorem 1.12(iii),
$d(b_{1,i}/a_{1,i})=d(a_{1,i}b_{1,i})\geq\alpha_i$. Also
$\ord a_{1,i}=\ord b_{1,i}=R_1+\cdots +R_i$, so
$b_{1,i}/a_{1,i}\in\ooo$. It follows that
$b_{1,i}/a_{1,i}\in\upo{\alpha_i}$. 

If $R_{i+1}-R_i>2e$, then $g(a_{i+1}/a_i)=\oos$, regardless of the
BONG, so we have (ii). By [B4, Corollary 2.8(ii)], $\alpha_i>2e$, so
$b_{1,i}/a_{1,i}\in\upo{\alpha_i}=\oos =g(a_{i+1}/a_i)$, so we have
(iii). From now on we will assume that $R_{i+1}-R_i\leq 2e$.

If $\theta (O^+(L))\subseteq\ooo\ffs$, then
$G(a_{i+1}/a_i)\subseteq\ooo\ffs$ so, by Lemma 5.1(i),
$g(a_{i+1}/a_i)=\upo{(R_{i+1}-R_i)/2+e}$, regardless of the BONG. This
proves (ii). Also, by Lemma 5.1(ii), $\alpha_i=(R_{i+1}-R_i)/2+e$, so
$b_{1,i}/a_{1,i}\in\upo{\alpha_i}=\upo{(R_{i+1}-R_i)/2+e}=g(a_{i+1}/a_i)$.
So we have (iii).

We are left with the case when $L$ has property B, which, by (i), is
independent of the BONG. By the proof of (i),
$d(-a_{i,i+1})\leq e-(R_{i+1}-R_i)/2$ iff
$d(-b_{i,i+1})\leq e-(R_{i+1}-R_i)/2$. So we have two cases.

(a) $d(-a_{i,i+1}),\,d(-b_{i,i+1})>e-(R_{i+1}-R_i)/2$. We have
$g(a_{i+1}/a_i)=g(b_{i+1}/b_i)=\upo{(R_{i+1}-R_i)/2+e}$ so (ii)
holds. By Lemma 5.2(i), $\alpha_i=(R_{i+1}-R_i)/2+e$, so
$b_{1,i}/a_{1,i}\in\upo{\alpha_i}=g(a_{i+1}/a_i)$ and so (iii) holds. 

(b) $d(-a_{i,i+1}),\,d(-b_{i,i+1})\leq e-(R_{i+1}-R_i)/2$. We have
$g(a_{i+1}/a_i)=\upon{R_{i+1}-R_i+d(-a_{i,i+1})}{-a_{i,i+1}}$. Also,
by Lemma 5.2(i), $\alpha_i=R_{i+1}-R_i+d(-a_{i,i+1})$, so
$b_{1,i}/a_{1,i}\in\upo{R_{i+1}-R_i+d(-a_{i,i+1})}$. To complete the
proof of (iii), we need to show that $b_{1,i}/a_{1,i}\in\N (-a_{i,i+1})$.
Since $d(-a_{i,i+1})\leq e-(R_{i+1}-R_i)/2$, property B implies that
$R_i-R_{i-1}>2e$ (or $i=1$) and $R_{i+2}-R_{i+1}>2e$ (or
$i=n-1$). The, by Corollary 1.13, we have
$[a_1,\ldots,a_{i-1}]\cong [b_1,\ldots,b_{i-1}]$ and
$[a_1,\ldots,a_{i+1}]\cong [b_1,\ldots,b_{i+1}]$, which imply that
$[a_i,a_{i+1}]\cong [b_i,b_{i+1}]$, and also
$a_{1,i-1}b_{1,i-1}=1$ in $\fff/\ffs$. From
$[a_i,a_{i+1}]\cong [b_i,b_{i+1}]$
we get that $a_{i+1}/a_i=b_{i+1}/b_i$ in $\fff/\ffs$. But
$\ord a_{i+1}/a_i=\ord b_{i+1}/b_i=R_{i+1}-R_i$, so in fact
$a_{i+1}/a_i=b_{i+1}/b_i$ in $\fff/\oos$, which implies that
$g(a_{i+1}/a_i)=g(b_{i+1}/b_i)$, so (ii) holds. We also have
$b_i\rep [a_i,a_{i+1}]$, which implies that
$(a_ib_i,-a_{i,i+1})_\p =1$, i.e. $a_ib_i\in\N (-a_{i,i+1})$. But in
$\fff\ffs$ we have $a_{1,i-1}b_{1,i-1}=1$, so
$b_{1,i}/a_{1,i}=a_ib_i$. Hence $b_{1,i}/a_{1,i}\in\N (-a_{i,i+1})$,
which concludes the proof of (iii). \qed

\bff Let $\mathcal S_L:=\{ (a_1,\ldots,a_n)\in (\fff/\oos )^n\,
\mid\, L\cong{\prec a_1,\ldots,a_n\succ}\text{ relative to a good
BONG}\}$.

Since $g(a_{i+1}/a_i)$ are invariants of the lattice $L$ we may denote
$g_i(L)=g(a_{i+1}/a_i)$. If $\eta\in g_i(L)=g(a_{i+1}/a_i)$, then, by
Lemma 1.5(iv) $L\cong{\prec a_1,\ldots,a_{i-1},\eta a_i,\eta a_{i+1},
a_{i+2},\ldots,a_n\succ}$, so
$(a_1,\ldots,\eta a_i,\eta a_{i+1},\ldots,a_n)\in\mathcal S_L$. Hence
for every $\eta\in g_i(L)$ we have a map
$t_i(\eta ):\mathcal S_L\to\mathcal S_L$, given by
$(a_1,\ldots,a_n)\mapsto (a_1,\ldots,\eta a_i,\eta a_{i+1},\ldots,a_n)$.

Obviously, every two such maps, $t_i(\eta )$ and $t_j(\eta')$,
commute. Mereover, we have $t_i(\eta )t_i(\eta')=t_i(\eta\eta')$.
\eff

\blm (i) The map $t_i(\eta )$ changes the product $a_{1,i}$ to
$\eta a_{1,i}$ and leaves all products $a_{1,j}$ with $j\neq i$
unchanged (in $\fff/\oos$). 

(ii) If $1\leq i\leq n-1$ and $(a_1,\ldots,a_n)\in\mathcal S_L$, we
denote by $l_i(a_1,\ldots,a_n)$ the class of
${\prec a_1,\ldots,a_i\succ}$ and by $r_i(a_1,\ldots,a_n)$ the class of
${\prec a_{i+1},\ldots,a_n\succ}$. Then for every $j\neq i$ and
$\eta\in g_j(L)$ we have $l_it_j(\eta )=l_i$ and $r_it_j(\eta )=r_i$.
\elm
\pf (i) If $(a'_1,\ldots,a'_n)=t_i(\eta)(a_1,\ldots,a_n)$, then
$a'_i=\eta a_i$, $a'_{i+1}=\eta a_{i+1}$ and $a'_j=a_j$ for
$j\neq i,i+1$. Then for $j<i$ we have $a'_{1,j}=a_{1,j}$,
$a'_{1,i}=\eta a_{1,i}$ and for $j>i$ we have
$a'_{1,j}=\eta^2a_{1,j}=a_{1,j}$ in $\fff/\oos$.

(ii) Let $a=(a_1,\ldots,a_n)$. If $j<i$, then
$r_it_j(\eta )(a)={\prec a_{i+1},\ldots a_n\succ} =r_i(a)$ and
$l_it_j(\eta )(a)={\prec a_1,\ldots,\eta a_j,\eta
a_{j+1},\ldots,a_i\succ}$. But $\eta\in g_j(L)=g(a_{j+1}/a_j)$ so, by
Lemma 1.5(iv),
${\prec a_1,\ldots,\eta a_j,\eta a_{j+1},\ldots,a_i\succ}
={\prec a_1,\ldots,a_i\succ} =l_i(a)$. Similarly, 
if $j>i$, then $l_it_j(\eta )(a)={\prec a_1,\ldots a_i\succ} =l_i(a)$
and $r_it_j(\eta )(a)={\prec a_{i+1},\ldots,\eta a_j,\eta a_{j+1},
\ldots a_n\succ}$. Then, again by Lemma 1.5(iv), we have 
${\prec a_{i+1},\ldots,\eta a_j,\eta a_{j+1},\ldots a_n\succ}
={\prec a_{i+1},\ldots a_n\succ} =r_i(a)$. \qed

\blm (i) If
$L\cong{\prec a_1,\ldots,a_n\succ}\cong{\prec b_1,\ldots,b_n\succ}$
relative to some good BONGs, then
$(b_1,\ldots,b_n)=t_1(\eta_1)\cdots t_{n-1}(\eta_{n-1})(a_1,\ldots,_n)$
for some $\eta_i\in g_i(L)$, Namely, $\eta_i=b_{1,i}/a_{1,i}$.

(ii) If $a_{1,i}=b_{1,i}$ for some $1<i<n$, then
${\prec a_1,\ldots,a_i\succ}\cong{\prec b_1,\ldots,b_i\succ}$ and
${\prec a_{i+1},\ldots,a_n\succ}\cong{\prec b_{i+1},\ldots,b_n\succ}$.
\elm
\pf (i) By Lemma 5.3(iii),
$\eta_i:=b_{1,i}/a_{1,i}\in g(a_{i+1}/a_i)=g_i(L)$. Let
$t_1(\eta_1)\cdots t_{n-1}(\eta_{n-1})(a_1,\ldots,a_n)
=(c_1,\ldots,c_n)$. By Lemma 5.5(i), for $1\leq i\leq n-1$ the product
$a_{1,i}$ is multiplied by $\eta_i$ when we apply $t_i(\eta_i)$ and it
is left unchanged by $t_j(\eta_j)$ for $j\neq i$. Therefore
$c_{1,i}=\eta_ia_{1,i}=b_{1,i}$ for $1\leq i\leq n-1$. Same happens
for $i=n$, as $a_{1,n}=b_{1,n}=\det L$. It follows that $c_i=b_i$ for
$1\leq i\leq n$.

(ii) Let $a=(a_1,\ldots,a_n)$ and $b=(b_1,\ldots,b_n)$. Since
$l_i(a)$, $l_i(b)$, $r_i(a)$ and $r_i(b)$ are the classes of
${\prec a_1,\ldots,a_i\succ}$, ${\prec b_1,\ldots,b_i\succ}$,
${\prec a_{i+1},\ldots,a_n\succ}$ and ${\prec b_{i+1},\ldots,b_i\succ}$,
we must prove that $l_i(a)=l_i(b)$ and $r_i(a)=r_i(b)$.

By (i) we have $b=t(a)$ where
$t=t_1(\eta_1)\cdots t_{n-1}(\eta_{n-1})$. But
$\eta_i=b_{1,i}/a_{1,i}=1$, so $t_i(\eta_i)=t_i(1)=1$. Hence
$t=t_1(\eta_1)\cdots\widehat{t_i(\eta_i)}\cdots t_{n-1}(\eta_{n-1})$.
Since $l_it_j(\eta_j)=l_i$ and $r_it_j(\eta_j)=r_i$ for every
$j\neq i$, we have $l_it=l_i$ and $r_it=r_i$. We conclude that
$l_i(b)=l_it(a)=l_i(a)$ and $r_i(b)=r_it(a)=r_i(a)$. \qed

\bco If $x$ and $y$ are first elements in good BONGs of $L$ and
$Q(x)=Q(y)$, then $y=\sigma (x)$ for some $\sigma\in O^+(L)$.
\eco
\pf We have $L\cong{\prec a_1,\ldots,a_n\succ}$ and
$L\cong{\prec b_1,\ldots,b_n\succ}$ relative to some good BONGs
$x=x_1,\ldots,x_n$ and $y=y_1,\ldots,y_n$. Since $Q(x)=Q(y)$,
i.e. $a_1=b_1$, by Lemma 5.6(ii), we have that
${\prec y_2,\ldots,y_n\succ}\cong{\prec b_2,\ldots,b_n\succ}
\cong{\prec a_2,\ldots,a_n\succ}$ relative to another good BONG
$z_2,\ldots,z_n$. Then, by Lemma 1.5(iv),
$L\cong{\prec b_1,a_2,\ldots,a_n\succ} ={\prec a_1,\ldots,a_n\succ}$
relative to the BONG $y_1,z_2,\ldots,z_n$. Then the linear map
$\sigma$, given by $x_1\mapsto y_1$, $x_2\mapsto z_2$, ...,
$x_n\mapsto z_n$, is an isometry sending $x_1=x$ to $y_1=y$. If
$\sigma\in O^+(L)$, then we are done. If $\sigma\in O^-(L)$, then by
Lemma 1.5(v), we also have $\tau_{x_2}\in O^-(L)$, so
$\sigma':=\sigma\tau_{x_2}\in O^+(L)$. Since $x=x_1$ and $x_2$ are
mutually orthogonal, we have $\tau_{x_2}x=x$, so
$\sigma'(x)=\sigma (x)=y$ and so $\sigma'$ has the required
property. \qed

\section{Two reduction steps}

In this section $M\cong{\prec a_1,\ldots,a_m\succ}$ and
$N\cong{\prec b_1,\ldots,b_n\succ}$ relative to some good BONGs
$x_1,\ldots,x_m$ and $y_1,\ldots,y_n$, $N\subseteq M$ and $M$ has
property B or $\theta (O^+(M))=\ooo\fff^2$. We denote by $R_i=R_i(M)$
for $1\leq i\leq m$ and $S_i=R_i(N)$ for $1\leq i\leq n$.

The first reduction step applies when $\nnn M\supset\nnn N$ and
consists of decresing the norm of $M$.

Same as in \S4, we denote by
$M'=\{ x\in M\mid x\text{ not a norm generator}\}$.

\blm If $\theta (X(M/N))\subseteq\ooo\ffs$, then
$R_1\equiv\cdots\equiv R_m\equiv S_1\equiv\cdots\equiv S_n\pmod 2$.
\elm
\pf Since $\theta (O^+(M)),\,\theta O^+(N))\subseteq\theta (X(M/N)
\subseteq\ooo\ffs$, by the Remark following Theorem 1.9, we have
$R_1\equiv\cdots\equiv R_m\pmod 2$ and
$S_1\equiv\cdots\equiv S_n\pmod 2$. To complete the proof we will show
that $R_1\equiv S_1\pmod 2$, i.e. $R_1+S_1$ is even. By
Lemma 1.5(v), we have $a_1\in\theta (O^-(M))$ and
$b_1\in\theta (O^-(N))$. Since $O^-(N)O^-(M)\subseteq X(M/N)$, this
implies $b_1a_1\in\theta (X(M/N))\subseteq\ooo\ffs$, so
$R_1+S_1=\ord a_1b_1$ is even. \qed

\blm If $N\subset M$ and $\nnn N\subset\nnn M$, then $X(M/N)=X(M'/N)$
if $M'$ is a lattice and $\theta (X(M/N))=\fff$ otherwise.

In particular, when $N=M'$ we get
$O^+(M)\subseteq X(M/M')=X(M'/M')=O^+(M')$.
\elm
\pf Since $Q(N)\subseteq\nnn N\subset\nnn M$, we have
$N\subseteq M'$.

If $M'$ is a lattice, then for every $\sigma\in X(M/N)$ we have
$\sigma^{-1}N\subseteq M$ and $\nnn\sigma^{-1}N=\nnn N\subset\nnn M$,
so $\sigma^{-1}N\subseteq M'$ and so $\sigma\in X(M'/N)$. Thus
$X(M/N)=X(M'/N)$.

If $M'$ is not a lattice then, by Lemma 4.3(ii), $M=H\perp K$, where
$H\cong\frac 12\pi^{R_1}A(0,0)$ relative to some basis $x,y$. We have
$Q(x)=Q(y)=0$. Also $\ooo\ffs=\theta (O^+(H))\subseteq\theta (O^+(M))
\subseteq\theta (X(M/N))$ and so $\theta (X(M/N))=\ooo\ffs$ or $\fff$.

Since for every intermediate lattice $N\subseteq\bar N\subseteq M$ we
have $X(M/\bar N)\subseteq X(M/N)$, it is enough to prove that
$X(M/N)=\fff$ when $N$ is a maximal sublatice of $M$ with the property
that $N\subseteq M'$, i.e. $\nnn N\subset\nnn M$, so $S_1>R_1$. We
prove that $S_1=R_1+1$. Then, by Lemma 6.2, since
$R_1\not\equiv S_1\pmod 2$, we cannot have $\theta (X(M/N))=\ooo\ffs$
and so $\theta (X(M/N))=\fff$.

By Lemma 4.3(ii), $\p M+N\subseteq\p M+M'=M'$, so, by the maximality
of $N$, we must have $\p M+N=N$, i.e. $\p M\subseteq N$. Since
$\nnn H=\nnn M\supset\nnn N$, we cannot have $H\subseteq N$, so $x$ or
$y\notin N$. Say, $x\notin N$. Since $N\subset\oo x+N$, by the
maximality of $N$, we cannot have $\oo x+N\subseteq M'$. Thus
$\nnn (\oo x+N)=\nnn M$. Then there are $\alpha\in\oo$ and $z\in N$
such that $\ord Q(\alpha x+z)=\ord\nnn M=R_1$, i.e.
$\ord (Q(z)+2\alpha B(x,z))=R_1$. Also, since $z\in N$, we have
$\ord Q(z)\geq\ord\nnn N=S_1>R_1$. It follows that
$\ord 2\alpha B(x,z)=R_1$. If $z'=\pi\alpha x+z$, then $z'\in\p M+N=N$
and we have $Q(z')=Q(z)+2\pi\alpha B(x,z)$. Since $z,z'\in N$, we have
$2\pi\alpha B(x,z)=Q(z')-Q(z)\in\nnn N=\p^{S_1}$. Then
$R_1<S_1=\ord\nnn N\leq\ord 2\pi\alpha B(x,z)=\ord 2\alpha
B(x,z)+1=R_1+1$ and so $S_1=R_1+1$, as claimed. \qed

The second reduction step applies when $\nnn M=\nnn N$ and consists of
decreasing the ranks of the two lattices.

\blm Assume that $R_1=S_1$. Let $M^*=pr_{y_1^\perp}M$ and
$N^*=pr_{y_1^\perp}N$. Then in each of the following cases: 
$$\text{(a) }R_1<R_3\text{ (if }m\geq 3),\text{ (b) }R_3-R_2=-2e,
\text{ (c) }R_2=S_2$$
we have:

(i) $y_1$ is a first element in a good BONG of $M$, $b_1=\eta a_1$ for
some $\eta\in g(a_2/a_1)$,
$M^*\cong{\prec\eta a_2,a_3,\ldots,a_m\succ}$ 
and $N^*\cong{\prec b_2,\ldots,b_n\succ}$ relative to good BONGs.

(ii) $X(M/N)=X(M^*/N^*)O^+(M)$, $O^+(M^*)\subseteq O^+(M)$ and
$O^+(N^*)\subseteq O^+(N)$.
\elm
\pf (i) By Corollary 3.9, in each of the cases (a), (b) and (c) we
have that $y_1$ is a first element in a BONG of $M$. Then, by the case
$i=1$ of Lemma 5.3(iii), we have $\eta :=b_1/a_1\in g(a_2/a_1)$. By
Lemma 1.5(iv), we have
$M\cong{\prec\eta a_1,\eta a_2,a_3,\ldots.a_m\succ}$
relative to some BONG $x'_1,x'_2,x_3,\ldots,x_m$. Since
$Q(x'_1)=\eta a_1=b_1=Q(y_1)$, by Corollary 5.7, there is
$\sigma\in O(M)$ with $\sigma x'_1=y_1$. It follows that
$M=\sigma (M)\cong{\prec\eta a_1,\eta a_2,a_3,\ldots.a_m\succ}$
relative to the good BONG
$y_1=\sigma x'_1,\sigma x'_2,\sigma x_3,\ldots,\sigma x_m$. Then
$M^*\cong{\prec\eta a_2,a_3,\ldots.a_m\succ}$ relative to the BONG
$\sigma x'_2,\sigma x_3,\ldots,\sigma x_m$. We also have
$N^*\cong{\prec b_2,\ldots,b_n\succ}$ relative to the BONG
$y_2,\ldots,y_n$.

(ii) As seen in the proof of (i), we have
$M\cong{\prec b_1,\eta a_2,a_3,\ldots,a_m\succ}$ relative to a BONG
$y_1,z_2,\ldots,z_m$. By [B2, 2.5] we have the embedding
$X(M^*/N^*)\subseteq X(M/N)$. Thus
$X(M^*/N^*)O^+(M)\subseteq X(M/N)O^+(M)=X(M/N)$. Conversely, assume
that $\sigma\in X(M/N)$. We have $\sigma^{-1}N\subseteq M$ and
$\sigma ^{-1}y_1$ is a first element in a good BONG of
$\sigma^{-1}N$. By applying Corollary 3.9, to $M$ and
$\sigma^{-1}N\cong N$, we have that $\sigma^{-1}y_1$ is a first
element in a good BONG of $M$. Since $Q(\sigma^{-1}y_1)=Q(y_1)$, by
Corollary 5.7, there is some $\tau\in O^+(M)$ with
$\tau y_1=\sigma^{-1}y_1$. Hence $y_1=\psi y_1$, with
$\psi :=\sigma\tau\in X(M/N)O^+(M)=X(M/N)$. By [B2, 2.3], this implies
that $\psi\in X(M^*/N^*)$, so
$\sigma =\psi\tau^{-1}\in X(M^*/N^*)O^+(M)$.

In particular, if we take $N=M$, in which case we have $R_1=S_1$ and
$R_2=S_2$, so (c) holds, we get $X(M/M)=X(M^*/M^*)O^+(M)$,
i.e. $O^+(M)=O^+(M^*)O^+(M)$, so $O^+(M^*)\subseteq O^+(M)$. Similarly
$O^+(N^*)\subseteq O^+(N)$. \qed

\blm If $R_1=S_1$ and none of the conditions (a), (b) or (c) of Lemma
6.3 is satisfied, then $\theta (X(M/N))=\fff$.
\elm
\pf We use induction on $m$. If $m=2$ then (a) holds by default, so we
have nothing to prove. Suppose now that $m\geq 3$. We have $R_1=R_3$,
$R_3-R_2>-2e$ and $S_2\neq R_2$. But, by Lemma 3.1(iii),
$R_1+R_2\leq S_1+S_2=R_1+S_2$ so $S_2\geq R_2$. Hence $S_2>R_2$ (or
$n=1$).

Assume that $\theta (X(M/N))\neq\fff$. Since $R_1=R_3$, so $M$
doesn't have property A, by Proposition 1.6(iii), we have
$\ooo\fff^2\subseteq\theta (O^+(M))\subseteq\theta (X(M/N))$. Hence we
must have $\theta (O^+(M))=\theta (X(M/N))=\ooo\fff^2$. By Theorem
1.9, $R_1=R_3$ implies that $(R_2-R_1)/2\equiv e\pmod 2$. We also have
$R_1\equiv\cdots\equiv R_m\pmod 2$

First we reduce to the case when $S_2=R_2+2$. We have $n=1$ or
$S_1+S_2=R_1+S_2>R_1+R_2$ so, by Corollary 3.4, there is $K$, with
$N\subseteq K\subset M$ and $FK=FM$, such that $R_1+R_2<T_1+T_2\leq
R_1+R_2+2$, where $T_i=R_i(K)$. By taking the orders in
$\nnn N\subseteq\nnn K\subseteq\nnn M$, we get
$S_1\geq T_1\geq R_1=S_1$, so $R_1=T_1$. Then we have
$R_1<T_1\leq R_1+2$. But
$\theta (X(M/K))\subseteq\theta (X(M/N))=\ooo\fff^2$ so, by Lemma 6.1, 
$R_2\equiv T_2\pmod 2$. Thus $T_2=R_2+2$. We have $R_1=T_1$, $R_2<T_2$
and if we prove that $\theta (X(M/K))=\fff$ then also
$\theta (X(M/N))=\fff$ and we are done. Hence, by replcing $N$ by $K$,
we reduce to the case $S_2=R_2+2$.

Note that if $S_1=S_3$ then, by Theorem 1.9, we have
$(R_2+2-R_1)/2=(S_2-S_1)/2\equiv e\pmod 2$. But this contradicts
$(R_2-R_1)/2\equiv e\pmod 2$. Hence $S_3>S_1=R_1$.

Next we reduce to the case when $x_1=y_1$. To do this, we must prove
that $N$ contains a first element in a good BONG of $M$.

Since $R_1=R_3$, by Lemma 4.11(ii), we have $M=H\perp K$, with
$H\cong\pi^{(R_1+R_2)/2}A(0,0)$ relative to some basis $x,y$ and
$K\cong{\prec\pi^{R_1}\varepsilon,a_4,\ldots,a_m\succ}$, where
$\varepsilon\in\ooo$. We cannot have $N\subseteq\p H\perp K$, for this
would imply $X(M/N)\supseteq X(H\perp K/\p H\perp K)\supseteq X(H/\p
H)$. But, since $H'$ is not a lattice and $\nnn\p H\subset\nnn H$, we
have by Lemma 6.2 that $\theta (X(H/\p H))=\fff$ and so so we get
$\theta (X(M/N))=\fff$.

Hence $N$ contains elements $v=z+t$, with $z\in K$ and $t\in
H\setminus\p H$. We claim that we have such elements that are norm
generators. We have two elements of $N=H\perp K$, $v'=z'+t'$ and
$v''=z''+t''$, with $z',z''\in K$ and $t,t'\in H$ such that
$t'\notin\p H$ and $v''$ is a norm generator. (For $v'$ see above and
for $v''$ we may take $v''=y_1$.) If $v'$ is a norm generator or
$t''\notin\p H$, then $v=v'$ or $v=v''$, respectively, 
has the desired properties. So we may assume that $v'$ is not a norm
generator and $t''\in\p H$. Then we claim that $v=v''-v'$ has the
desired properties. We have $v=z+t$, where $z=z''-z'\in K$ and
$t=t''-t'\in H$. But $t'\in H\setminus\p H$ and $t''\in\p H$, so
$t\in H\setminus\p H$. We must prove that also $v$ is a norm generator
of $N$. Assume the contrary. Then $v,v'\in N$ are not norm generators,
but $v''=v+v'$ is. Hence $N'=\{ u\in N\mid x\text{ not a norm
generator}\}$ is not a lattice. By Lemmas 4.3(iii) and 4.9, this
implies that $S_2-S_1=-2e$. Since $R_1=S_1$ and $R_2<S_2$, we get
$R_2-R_1<S_2-S_1=-2e$, which is impossible. Hence $v$ is a norm
generator of $N$. 

Since $N\subseteq M$ and $\nnn M=\nnn N=\p^{R_1}$, $v$ is also a norm
generator of $M$, with $\ord Q(v)=R_1$. We claim that $v$ is a first
element in a good BONG of $M$. By Lemma 3.7(i), it suffices to prove
that $pr_{v^\perp}M$ contains an element $u$ with $\ord Q(u)=R_2$. We
write $t\in H$ in terms of the basis $x,y$ as $t=\alpha x+\beta y$,
with $\alpha,\beta\in\oo$. Since $t\notin\p H$ we may assume that, say
$\alpha\in\ooo$. We prove that $u:=pr_{v^\perp}y$ satisfies
$\ord Q(u)=R_2$. We have $Q(u)=(Q(v)Q(y)-B(v,y)^2)Q(v)^{-1}$. But
$Q(x)=Q(y)=0$, $B(x,y)=\pi^{(R_1+R_2)/2}$ and $z\in K$, so
$B(z,x)=B(y,z)=0$. Hence
$B(v,y)=B(z+\alpha x+\beta y,y)=\pi^{(R_1+R_2)/2}\alpha$. It follows
that
$Q(v)Q(y)-B(v,y)^2=0-(\pi^{(R_1+R_2)/2}\alpha )^2=-\pi^{R_1+R_2}\alpha^2$
and so $Q(u)=-\pi^{R_1+R_2}\alpha^2Q(v)^{-1}$. Since $\alpha\in\ooo$
and $\ord Q(v)=R_1$, we have $\ord Q(u)=R_2$, as claimed. So $v$ is a
first element in a good BONG of $M$. Since $v$ is a norm generator for
$N$ and $S_1<S_3$, by Lemma 3.8, $v$ s a first element of a good BONG
of $N$. Then, by replacing the good BONGs $x_1,\ldots,x_m$ and
$y_1,\ldots,y_n$ of $M$ and $N$ with
good BONGs that begin with $v$, we may assume that $x_1=y_1=v$.

If $M^*=pr_{v^\perp}M$ and $N^*=pr_{v^\perp}N$ then 
$M^*\cong{\prec a_2,\ldots,a_m\succ}$ and
$N^*\cong{\prec b_2,\ldots,b_n\succ}$ relative
to the good BONGs $x_2,\ldots,x_m$ and $y_2,\ldots,y_n$. We have
$N^*\subseteq M^*$ and, by [B2, Corollary 2.4],
$X(M^*/N^*)\subseteq X(M/N)$. Then $\theta (O^+(M^*))\subseteq
\theta (X(M^*/N^*))\subseteq\theta (X(M/N))=\ooo\fff^2$. The
invariants $R_1,R_2,\ldots$ corresponding to $M^*$ and $N^*$ are
$R_2,\ldots,R_m$ and $S_2,\ldots,S_n$. Since $S_2>R_2$, we have
$N^*\subseteq K\subset M^*$, where
$K=(M^*)'=\{ u\in M^*\mid u\text{ not a norm generator}\}$. Since
$R_3-R_2\neq -2e$ and $\theta (O^+(M^*))\subseteq\ooo\fff^2$ so $M^*$
is in one of the cases (iii) and (v) of Lemma 4.14. (By Lemma
4.13(ii), ${\prec a_1,a_2\succ}$ cannot be maximal, so we are not in
the case (iv).) So we have two cases:

a. $R_2<R_4$ or $m=3$. By Lemma 4.14, case (iii) we have
$K={\prec\pi x_2,\ldots,x_m\succ}$. Then $N^*\subseteq K\subset M^*$
writes as ${\prec y_2,\ldots,y_n\succ}\subseteq{\prec\pi
x_2,x_3,\ldots,x_m\succ}\subseteq{\prec x_2,\ldots,x_m\succ}$. Since
$M={\prec v,x_2,\ldots,x_m\succ}$, by [B2, Lemma 2.2], the lattice
$L={\prec v,\pi x_2,x_3,\ldots,x_m\succ}$ exists and we have
${\prec v,y_1,\ldots,y_n\succ}\subseteq{\prec v,\pi
x_2,x_3,\ldots,x_m\succ}\subseteq{\prec v,x_2,\ldots,x_m\succ}$,
i.e. $N\subseteq L\subset M$. It follows that $\theta (O^+(L))
\subseteq\theta (X(M/L))\subseteq\theta (X(M/N))=\ooo\fff^2$. Since
the BONG $\pi  x_2,x_3,\ldots,x_m$ of $K$ is good and
$\ord Q(v)=R_1=R_3=\ord Q(x_3)$, the BONG $v,\pi x_2,x_3,\ldots,x_m$
of $K$ is good as well. Hence if $T_i=R_i(L)$, then the sequence
$T_1,\ldots,T_m$ is $R_1,R_2+2,R_3,\ldots,R_m$. Since
$\theta (O^+(L))\subseteq\ooo\fff^2$ and $T_1=T_3$, by Theorem 1.9, we
have $(T_2-T_1)/2\equiv e\pmod 2$, i.e.
$(R_2+2-R_1)/2\equiv e\pmod 2$. But this contradicts
$(R_2-R_1)/2\equiv e\pmod 2$.

b. $R_2=R_4$. By the case (v) of Lemma 4.14, if the sequence of the
invariants $R_1,R_2,\ldots$ of $K=(M^*)'$ is $T_2,\ldots,T_m$, then
$T_2=T_4=R_2+2=S_2$ and $T_3=R_3-2=R_1-2$. It follows that
$T_4-T_3=R_2-R_1+4\geq -2e+4>-2e$, and
$S_3>S_1=R_1=R_3>T_3$. Thus $T_2=S_2$ and $K$ and $N^*$ don't satisfy
any of the conditions (a), (b) and (c) of our lemma. (We have
$T_2=T_4$, $T_2<S_2$ and $T_4-T_3>-2e$.) By the induction
hypothesis, $\theta (X(K/N^*))=\fff$. But $N^*\subseteq K\subset M^*$
so $\theta (X(K/N^*))\subseteq\theta (X(M^*/N^*))\subseteq\theta
(X(M/N))=\ooo\fff^2$. Contradiction. \qed

\section{Main results}

In this section we calculate $\theta (X(M/N))$ for two lattices
$N\subseteq M$. This is done in Theorems 7.3 and 7.5, which are
analogues of Theorems 1.7 and 1.9] for $\theta (O^+(L))$. 

We write $M\cong{\prec a_1,\ldots,a_m\succ}$ and
$N\cong{\prec b_1,\ldots,b_n\succ}$ relative to some good BONGs
$x_1,\ldots,x_m$ and  $y_1,\ldots,y_n$. Let $R_i=R_i(M)$ and
$S_i=S_i(N)$, $\alpha_i=\alpha_i(M)$ and $\beta_i=\alpha_i(N)$. 

Let $V=FM$ and $W=FN$ and let $U$ be the orthogonal complement of $W$
in $V$. Then $O^+(U)\subset X(M/N)$. If $m-n=\dim U\geq 3$ then
$\theta (O^+(U))=\fff$ and so $\theta (X(M/N))=\fff$.

Therefore we may restrict ourselves to the case when $m-n\leq 2$. If
$m-n$, we make the convention that $S_{n+1}\gg 0$.

Recall that for a lattice $L$ we have the notion of property A, which
menas $R_i<R_{i+2}$ $\forall i$. If $L$ has property A, we have a
formula for $\theta (O^+(L))$, given by Theorem 1.7. If $L$ doesn't
have property A, then $\theta (O^+(L))\supseteq\ooo\fff^2$, i.e.
$\theta (O^+(L))$ can only be $\ooo\fff^2$ or $\fff$ and Theorem
1.9. decides which of the two $\theta (O^+(L))$ is.

We have a similar situation when we want to calculate
$\theta X(M/N)$. For this purpose we extend the definition of property
A to a pair of lattices.

\bdf We say that the pair of lattices $M,N$, with $N\subseteq M$ and
$m-n\leq 2$, has property A if $R_{i+2}>S_i$ for $1\leq i\leq m-2$.
\edf 

Same as for $\theta (O^+(L))$, our result is given in two theorems. In
Theorem 7.3 we treat the case when the property A holds, which is
similar to Theorem 1.7. If property A doesn't hold, then
$\ooo\ffs\subseteq\theta (X(M/N))$, and Theorem 7.5, which is an
analogue of Theorem 1.9, decides if $\theta (X(M/N))=\ooo\ffs$ or
$\fff$.

For convenience, we denote by R1-R4 the contions (i)-(iv) from the
representation theorem, [B3, Theorem 4.5] or [B5, Theorem 2.1]. Since
so far this result has not been approved to be published, we will not
use them here. In particular, unlike in [B1], we will not use the
sufficiency of R1. Instead, we will prove that if R1 fails, then
$\theta (X(M/N))=\fff$.

\blm If $M,N$ satisfy R1, then:

(i) $R_i+R_{i+1}\leq S_i+S_{i+1}$ for all $1\leq i<n$.

(ii) For every $1<i<m$ such that $R_{i+1}\geq S_{i-1}$
we have $R_i\leq S_i$.
\elm
\pf (i) If $R_i\leq S_i$ and $R_{i+1}\leq S_{i+1}$, then
$R_i+R_{i+1}\leq S_i+S_{i+1}$ trivially. If $R_i>S_i$, then R1 implies
$R_i+R_{i+1}\leq S_{i-1}+S_i\leq S_i+S_{i+1}$. And if
$R_{i+1}\leq S_{i+1}$, then R1 implies
$R_i+R_{i+1}\leq R_{i+1}+R_{i+2}\leq S_i+S_{i+1}$.

(ii) Suppose that $R_i>S_i$. Then R1 implies
$R_i+R_{i+1}\leq S_{i-1}+S_i$, which, together with $R_i>S_i$, implies
$R_{i+1}<S_{i-1}$. Contradiction. \qed

\bco If the pair $M,N$ has property A then condition R1 is equivalent
to $R_i\leq S_i$ for all $1\leq i\leq n$. Consequently, it implies
that $M$ has property A.
\eco
\pf If $R_i\leq S_i$ for all $i$, then R1 holds by
definition. Conversely, if R1 holds, then, when $1<i<m$, $R_i\leq S_i$
follows, by Lemma 7.1(ii), from $R_{i+1}>S_{i-1}$. If $i=1$ or $m$, this
follows dirrectly from R1.

If the pair $M,N$ has property A and satifies R1, then for
$1\leq i\leq m-2$ we have $R_{i+2}>S_i\geq R_i$, so $M$ has property
A. \qed

\btm (i) If $N\subseteq M$, $m-n\leq 2$, and the pair of lattices
$M,N$ has property A and satisfies condition R1, then
$$\theta (X(M/N))=\ups\beta\prod_{i=1}^{m-2}
\bar G(a_{1,i+1}b_{1,i-1},R_{i+1}-S_i)\theta (O^+(M)),$$
where $\beta =\min\{ [(R_{i+2}-S_i)/2]\mid\, 1\leq i\leq m-2\}$.

(ii) If $R_1+\cdots +R_i\not\equiv S_1+\cdots +S_i\pmod 2$, then the
factor $\bar G(a_{1,i+1}b_{1,i-1},R_{i+1}-S_i)$ in the formula from
(i) can be replaced by ${\rm N}(-a_{1,i+1}b_{1,i-1})$.

(iii) If the pair $M,N$ doesn't have property A, then
$\theta (X(M/N))\supseteq\ooo\ffs$. If condition R1 fails, then
$\theta (X(M/N))=\fff$.  

If $m-n=2$, by our convention $S_{m-1}\gg 0$, so $R_m-S_{m-1}\ll 0$.
Then, by Corrolary 2.13(iii), the last factor of the product from (i),
$\bar G(a_{1,m}b_{1,m-2},R_m-S_{m-1})$, is equal to
${\rm N}(-a_{1,m}b_{1,m-2})$.

And if $m\leq 2$, then $\beta$ is the  minimum of an empty set so we
put $\beta =\infty$. Then the factor $\ups\beta =\ffs$ can be ignored
in the formula for $\theta (X(M/N))$.
\etm 

\bff {\bf Remarks.} (1) Theorem 7.1(i) states that
$$\theta (X(M/N))=\tilde G(M/N):=G(M/N)\theta (O^+(M)),$$
where
$G(M/N)=\ups\beta\prod_{i=1}^{m-1}\bar G(a_{1,i+1}b_{1,i-1},R_{i+1}-S_i)$.

By Corollary 7.2, since the pair $M,N$ has property A and satisfies
R1, we have $R_i\leq S_i$ for $1\leq i\leq n$. Then, since
$G(a_{i+1}/a_i)=\bar G(a_{i,i+1},R_{i+1}-R_i)$, Theorem 1.7
writes as
$\theta (O^+(M))=\ups\alpha\prod_{i=1}^{m-1}\bar G(a_{i,i+1},R_{i+1}-R_i)$,
where $\alpha =\min\{ [(R_{i+2}-R_i)/2]\mid 1\leq i\leq m-2\}$.

But for every $1\leq i\leq n-2$ we have $R_i\leq S_i$, so
$[(R_{i+2}-R_i)/2]\geq [(R_{i+2}-S_i)/2]$. By takeng minima, we get
$\alpha\geq\beta$, so $\ups\alpha\subseteq\ups\beta\subseteq G(M/N)$.
We also have $R_2-R_1\geq R_2-S_1$ so, by Lemma 2.12(i), 
$\bar G(a_{1,2},R_2-R_1)\subseteq\bar G(a_{1,2},R_2-R_1)\subseteq G(M/N)$.

Thus both factors $\ups\alpha$ and $\bar G(a_{1,2},R_2-R_1)$ of
$\theta (O^+(M))$ are superfluous in the product
$\tilde G(M/N)=G(M/N)\theta (O^+(M))$. It follows that
$$\tilde G(M/N)=G(M/N)\prod_{i=2}^{n-1}\bar G(a_{i,i+1},R_{i+1}-R_i).$$

(2) Since
$\bar G(a_{1,i+1}b_{1,i-1},R_{i+1}-S_i)\subseteq\N (-a_{1,i+1}b_{1,i-1})$,
(ii) is equivalent to
$\N (-a_{1,i+1}b_{1,i-1})\subseteq\theta (X(M/N))$ for every $i$ such
that $R_1+\cdots +R_i\not\equiv S_1+\cdots +S_i\pmod 2$.
\eff

\btm Suppose that $N\subseteq M$, $m-n\leq 2$ and the pair of
lattices $M,N$ satisfies R1. For $1\leq i\leq m-1$ we put
$T_i=\max\{ R_{i+1},S_{i-1}\} -\min\{ S_i,R_{i+2}\}$. (If $i=1$, we
ignore $S_{i-1}$ in $\max\{ R_{i+1},S_{i-1}\}$; if $i=m-1$, we ignore
$R_{i+2}$ in $\min\{ S_i,R_{i+2}\}$.)

Then $\theta (M/N)\subseteq\ooo\fff^2$ iff the following conditions
hold.

(1)  $R_1\equiv\cdots\equiv R_m\equiv S_1\equiv\cdots\equiv S_n\pmod 2$.

(2) If $R_{i+2}\leq S_i$ for some $1\leq i\leq m-2$, then one of the
following happens:

$\qquad$ (a) $R_{i+1}+R_{i+2}=S_i+S_{i+1}$ and
$(R_{i+2}-R_i)/2\equiv (S_{i+1}-S_i)/2\equiv e\pmod 2$.

$\qquad$ (b) $R_{i+2}=S_i$ and $-2e\in\{ R_{i+2}-R_{i+1},S_{i+1}-S_i\}$.

(3) For $1\leq i<m$ we have 
$\bar G(a_{1,i+1}b_{1,i-1},T_i)\subseteq\ooo\fff^2$.

(4) $\theta (O^+(M)),\theta (O^+(M))\subseteq\ooo\fff^2$.

If $m-n=2$, by our convention, $S_{m-1}\gg 0$, so we cannot have
$R_{m-1}+R_m=S_{m-2}+S_{m-1}$ or $S_{m-1}-S_{m-2}=-2e$. Thus the
condition (2) at $i=m-1$ states that if $R_m\leq S_{m-2}$, then
$R_m=S_{m-2}$ and $R_m-R_{m-2}=-2e$. Also in this case
$T_{m-1}=\max\{ R_m,S_{m-2}\} -S_{m-1}\ll 0$. So, by Lemma 2.17(ii),
condition (3) at $i=m-1$,
$\bar G(a_{1,m}b_{1,m-2},T_{m-1})\subseteq\ooo\ffs$, is equivalent to
$d(-a_{1,m}b_{1,m-2})=2e$. 
\etm

\bff {\bf Remarks.} (1) Assuming that (1) holds and
$R_{i+1}+R_{i+2}=S_i+S_{i+1}$, we note that
$(R_{i+2}-R_{i+1})/2=(R_{i+1}+R_{i+2})/2-R_{i+1}$ and
$(S_{i+1}-S_i)/2=(S_i+S_{i+1})/2-S_i$. Since
$R_{i+1}+R_{i+2}=S_i+S_{i+1}$ and $R_{i+1}\equiv S_i\pmod 2$, we have
$(R_{i+2}-R_{i+1})/2=(S_{i+1}-S_i)/2\pmod 2$. Thus the condition
$(R_{i+2}-R_{i+1})/2\equiv (S_{i+1}-S_i)/2\equiv e\pmod 2$ from (2)(a)
may be replaced by $(R_{i+2}-R_{i+1})/2\equiv e\pmod 2$ or
$(S_{i+1}-S_i)/2\equiv e\pmod 2$.

(2) An equivalent condition for (3) is the following:

{\it (3') For every $1\leq i\leq m-1$ we have either
$d(-a_{1,i+1}b_{1,i-1})=2e$ or $i\leq n$, $T_i\geq -2e$ and
$d(-a_{1,i+1}b_{1,i-1})>e-T_i/2$.}

Indeed, if $i=n+1=m-1$, then (3) and (3') are the same. And if
$i\leq n$, they are equivalent by Lemma 2.17(i). The missing
condition, that $\ord a_{1,i+1}b_{1,i-1}$ is even, is superfluous
by (1). We have that
$\ord a_{1,i+1}b_{1,i-1}=R_1+\cdots +R_{i+1}+S_1+\cdots +S_{i-1}$ is a
sum of $2i$ integers of the same parity, so it is even.
\eff

As a consequence of Theorems 7.3 and 7.5, when
$\theta (X(M/N))\neq\fff$ we have a stronger version of R1.

\bco If $\theta (X(M/N))\neq\fff$, then for all $1\leq i\leq n$ we
have $R_i\leq S_i$ or $1<i<m$ and $R_{i+1}+R_{i+2}=S_i+S_{i+1}$.
\eco
\pf Since $\theta (X(M/N))\neq\fff$, by Theorem 7.3(iii), the pair
$M,N$ satifies R1. If the pair $M,N$ has property A, then, by
Corollary 7.2, $R_i\leq S_i$ for all $1\leq i\leq n$ and we are
done. If the pair $M,N$ doesn't have property A, by Theorem 7.3(iii),
we have $\theta (X(M/N))\supseteq\ooo\ffs$. Thus we must have
$\theta (X(M/N))=\ooo\ffs$, so Theorem 7.5 applies. If $R_i>S_i$ for
some $1\leq i\leq n$, then, by R1, $R_i+R_{i+1}\leq S_{i-1}+S_i$.
Together with, $R_i>S_i$, this implies $R_{i+1}<S_{i-1}$, which, by
condition (2) of Theorem 7.5, implies $R_i+R_{i+1}=S_{i-1}+S_i$.
\qed

\blm Suppose that the pair $M,N$ satisfy condition R1 and condition
(2) of Theorem 7.5.

(i) $R_{i+1}+R_{i+2}\geq S_{i-1}+S_i$ for every $1<i\leq m-2$.

(ii) $T_i=\max\{ R_{i+1}-S_i,R_{i+1}-R_{i+2},S_{i-1}-S_i\}$.
Explicitly, $T_i=R_{i+1}-S_i$ if $R_{i+1}\geq S_{i-1}$ and
$R_{i+2}\geq S_i$, $T_i=R_{i+1}-R_{i+2}$ if $R_{i+2}\leq S_i$ and
$T_i=S_{i-1}-S_i$ if $R_{i+1}\leq S_{i-1}$.
\elm
\pf (i) If $R_{i+1}\geq S_{i-1}$ and $R_{i+2}\geq S_i$ the stement is
trivial. If $R_{i+1}<S_{i-1}$ or $R_{i+2}<S_i$, then, by the condition
(2) of Theorem 7.5, $R_{i+1}+R_{i+2}\geq R_i+R_{i+1}=S_{i-1}+S_i$ or
$R_{i+1}+R_{i+2}=S_i+S_{i+1}\geq S_{i-1}+S_i$, respectively.

(ii) If $R_{i+1}\geq S_{i-1}$ and $R_{i+2}\geq S_i$, then
$T_i=R_{i+1}-S_i$ trivially. By (i), we have
$R_{i+1}+R_{i+2}\geq S_{i-1}+S_i$. Thus if $R_{i+2}\leq S_i$, then
also $R_{i+1}\geq S_{i-1}$ and we get $T_i=R_{i+1}-R_{i+2}$. And if
$R_{i+1}\leq S_{i-1}$, then also $R_{i+2}\geq S_i$ and we get
$T_i=S_{i-1}-S_i$.

Note that $T_i=\max\{ R_{i+1},S_{i-1}\} +\max\{ -S_i,-R_{i+2}\}
=\max\{ R_{i+1}-S_i,R_{i+1}-R_{i+2},S_{i-1}-S_i,S_{i-1}-R_{i+2}\}$.
But $S_{i-1}-R_{i+2}$ can be removed, since we proved that $T_i$ is
one of the remaining terms of the maximum. \qed

\blm If $\theta (O^+(M))=\fff$, then Theorems 7.3 and 7.5 hold.
\elm
\pf Since $O^+(M)\subseteq X(M/N)$, $\theta (O^+(M))=\fff$ implies
$\theta (X(M/N))=\fff$. Now the product from Theorem 7.3(i) contains
$\theta (O^+(M))=\fff$ so it is equal to $\fff$. Thus Theorem
7.3 states that $\theta (X(M/N))=\fff$, which is true. And since
$\theta (O^+(M))=\fff\not\subseteq\ooo\ffs$, Theorem 7.5 states that
$\theta (X(M/N))\not\subseteq\ooo\ffs$, which is true. \qed

\blm Theorem 7.5 is independent of the choice of the BONG of $M$.
\elm
\pf All conditions of Theorem 7.5 except (3) are in terms of the
invariants $R_1,\ldots,R_m$, so they are independent of the BONG. So
we are left to proving that if conditions (1), (2) and (4) are
satisfied, then (3) is independent of the BONG. If $i=m-1=n+1$, then
(3) at $i$ states that
$-\det FM\det FN=-a_{1,i+1}b_{1,i-1}\in\Delta\ffs$, which
is independent of the BONG. For the remaining indices it states that
$\bar G(a_{1,i+1}b_{1,i-1},T_i)\subseteq\ooo\ffs$. Let
$M\cong{\prec a_1'\ldots,a_m'\succ}$ relative to another good BONG. For
convenience, let $a,b\in\fff/\ffs$, $a=a_{1,i+1}b_{1,i-1}$,
$b=a_{1,i+1}'b_{1,i-1}$. We must prove that if
$\bar G(a,T_i)\subseteq\ooo\ffs$ then also
$\bar G(b,T_i)\subseteq\ooo\ffs$. If $i=m-1$, then
$a=b=\det FM\, b_{1,m-2}$, so we are done. We assume that $i\leq m-2$.
As a consequence of (1),
$\ord a=\ord b=R_1+\cdots +R_{i+1}+S_1+\cdots +S_{i-1}$ is even. Then
$\ord a,\ord b>0$, so also $d(ab)>0$.

If $T_i\geq 2e$, since $\ord a$ and $\ord b$ are even, by Lemma
17.(ii), $\bar G(a,T_i),\bar G(T_i,b)\subseteq\ooo\ffs$ and we are
done. So we assume that $T_i<2e$.

If $T_i=-2e$ and $\{ a,b\} =\{ -1,-\Delta\}$, then
$\bar G(a,T_i)=\bar G(b,T_i)=\ooo\ffs$ and we are done. In the
remaining cases, by Lemma 2.16, we have
$\bar G(a,T_i)\bar G(b,T_i)=\ups{T_i/2-e+d(ab)}\bar G(a,T_i)$, so if
$\bar G(a,T_i)\subseteq\ooo\ffs$, then
$\bar G(b,T_i)\subseteq\ooo\ffs$ is equivalent to
$\ups{T_i/2-e+d(ab)}\subseteq\ooo\ffs$, i.e. to $T_i/2-e+d(ab)>0$.

In $\fff/\ffs$ we have $ab=a_{1,i+1}a_{1,i+1}'$. Then, by Theorem
1.12(iii) and Lemma 5.2(ii), we have
$d(ab)\geq\alpha_{i+1}=(R_{i+2}-R_{i+1})/2+e$. Since also
$T_i\geq R_{i+1}-R_{i+2}$, we have $T_i/2-e+d(ab)\geq 0$. If one of
the two equalities is strict, then $T_i/2-e+d(ab)>0$ and we are
done. So we may assume that $T_i=R_{i+1}-R_{i+2}$ and
$d(ab)=(R_{i+2}-R_{i+1})/2+e=e-T_i/2$.

Suppose that $T_i\leq -2e$. Then $d(ab)=e-T_i/2\geq 2e$, which can
only happen if $d(ab)=2e$, so $T_i=-2e$. Since $T_i=-2e$, by Lemma
2.17(ii), $\bar G(a,T_i)\subseteq\ooo\ffs$ implies
$a\in\{ 1,\Delta\}$. On the other hand, $d(ab)=2e$, so in $\fff/\ffs$
we have $ab=\Delta$, i.e. $b=\Delta a$. So we are back to the case
$T_i=-2e$ and $\{ a,b\} =\{ -1,-\Delta\}$, already discussed. So from
now on we assume that $-2e<T_i<2e$.

We have $R_{i+1}-R_{i+2}=T_i\geq R_{i+1}-S_i$ so $R_{i+2}\leq S_i$ and
so one of the cases (a) or (b) of (2) holds at $i+1$. Suppose fist
that (b) holds, i.e. $R_{i+2}=S_i$ and
$-2e\in\{ R_{i+2}-R_{i+1},S_{i+1}-S_i\}$. If
$R_{i+2}-R_{i+1}=-2e$, then $T_i=R_{i+1}-R_{i+2}=2e$, which
contradicts the assumption that $T_i<2e$. If $S_{i+1}-S_i=-2e$, we
note that, by Lemma 7.1(ii), $R_{i+2}=S_i$ implies
$R_{i+1}\leq S_{i+1}$. Then $T_i=R_{i+1}-R_{i+2}\leq S_{i+1}-S_i=-2e$,
which contradicts the assumption that $T_i>-2e$. 

If (a) holds, then $-T_i/2=(R_{i+2}-R_{i+1})/2\equiv e\pmod 2$, which
implies that $e-T_i/2$ is even. Since $T_i\in (-2e,2e)$, we also have
$e-T_i/2\in (0,2e)$. Thus $d(ab)=e-T_i/2$ is an even integer in the
interval $(0,2e)$, which is impossible. This concludes the proof. \qed

\blm Suppose that the pair of lattices $M,N$ satisfies R1 and the
condition (1) of Theorem 7.5.

(i) If $R_2-R_1=-2e$, then the condition (2) of Theorem 7.5 holds at
$i=1$ iff $R_3>S_1$ or $R_3=S_1$ and $S_2-S_1=-2e$.

(ii) Suppose that $R_1=R_{2k+1}$ for some $k\geq 1$,
$(R_1-R_{2j})/2\equiv e\pmod 2$ for $1\leq j\leq k$ and we have either
$R_1<S_1$ or $R_1-R_{2k}>-2e$.

Then the condition (2) of Theorem 7.5 holds at all indices
$1\leq i\leq\min\{ 2k,m-2\}$ iff $S_{2j-1}+S_{2j}=R_1+R_{2j}$ holds
for all $1\leq j\leq k$.

Moreover, if this happens, then for every $1\leq j\leq k$ we have
$S_{2j}\leq R_{2j}$, $T_{2j-1}=R_{2j}-R_1$ and $T_{2j}=S_{2j-1}-S_{2j}$.
\elm
\pf (i) The ``if'' part is trivial. Suppose now that condtion (2) of
Theorem 7.5 holds at $i=1$. If $R_3>S_1$, then we are done. If
$R_3\leq S_1$ then one of the cases (a) and (b) of (2) holds. If (b)
holds, then $R_3=S_1$ and $-2e\in\min\{ R_3-R_2,S_2-S_1\}$. Since
$R_3-R_2\geq R_1-R_2=2e$, we must have $S_2-S_1=-2e$, so we are
done. If (a) holds, then $R_2+R_3=S_1+S_2$, which, together with
$R_3\leq S_1$, implies $R_2\geq S_2$. Since also $R_1\leq S_1$
$-2e\leq S_2-S_1\leq R_2-R_1=-2e$, so we must have the equalities
$S_2-S_1=-2e$, $R_1=S_1$ and $R_2=S_2$. But from
$R_3\leq S_1\leq R_1\leq R_3$ we also have $R_3=S_1$, so we are
done. 

(ii) We have $R_1\leq R_3\leq\cdots\leq R_{2k+1}$ so $R_1=R_{2k+1}$
implies $R_1=R_3=\cdots =R_{2k+1}$. Also not that
$S_{2j-1}\geq S_1\geq R_1=R_{2j+1}$ for $1\leq j\leq k$.

Suppose Theorem 7.5(2) holds at $1\leq i\leq 2k-1$. For
$1\leq j\leq k$ we have $S_{2j-1}\geq R_{2j+1}$ so the condition (2)
of Theorem 7.5 applies at $i=2j-1$. If
$S_{2j-1}+S_{2j}=R_{2j}+R_{2j+1}=R_1+R_{2j}$, then we are
done. Suppose now that $S_{2j-1}=R_{2j+1}=R_1$ and
$-2e\in\{ R_{2j+1}-R_{2j},S_{2j}-S_{2j-1}\}$. Then 
$S_{2j-1}\geq S_1\geq R_1$ implies $S_1=R_1$, so by hypothesis, we
must have $R_1-R_{2k}>-2e$. Since $R_{2j+1}=R_1$ and
$R_{2j}\leq R_{2k}$, we have $R_{2j+1}-R_{2j}\geq R_1-R_{2k}>-2e$, so
we are left with the case when $S_{2j}-S_{2j-1}=-2e$. By Lemma 7.1(i),
property R1 implies that $S_{2j-1}+S_{2j}\geq R_{2j-1}+R_{2j}$. We
also have $R_{2j}-R_{2j-1}\geq -2e=S_{2j}-S_{2j-1}$. Since
$S_{2j-1}=R_1=R_{2j-1}$, these inequlities are equivalent to
$S_{2j}\geq R_{2j}$ and $R_{2j}\geq S_{2j}$. Thus $S_{2j}=R_{2j}$ and
so $S_{2j-1}+S_{2j}=R_1+R_{2j}$.

Conversely, assume that $S_{2j-1}+S_{2j}=R_1+R_{2j}$ for
$1\leq j\leq k$. Then for $1\leq j\leq k$ we have
$S_{2j-1}+S_{2j}=R_1+R_{2j}=R_{2j}+R_{2j+1}$ and
$(R_{2j+1}-R_{2j})/2=(R_1-R_{2j})/2\equiv e\pmod 2$, so condition (2)
of Theorem 7.5 holds at $i=2j-1$. (See also Remark 7.6(1).) If
$i=2j$, with $1\leq j\leq k$, we note that
$S_{2j-1}\geq S_1\geq R_1$, so
$R_1+S_{2j}\leq S_{2j-1}+S_{2j}=R_1+R_{2j}$ and so
$S_{2j}\leq R_{2j}\leq R_{2j+2}$. If $S_{2j}<R_{2j+2}$, condition (2)
of Theorem 7.5 holds trivially at $i=2j$. Otherwise
$S_{2j}=R_{2j}=R_{2j+2}$ and  we have
$S_{2j-1}+S_{2j}=R_1+R_{2j}=R_{2j+1}+R_{2j+2}$. Since
also $(R_{2j+2}-R_{2j+1})/2=(R_{2j}-R_1)/2\equiv e\pmod 2$, the
condition (2) of Theorem 7.5 holds at $i=2j$.

If $1\leq j\leq k$, then from $S_{2j-1}\geq R_1=R_{2j+1}$ we get, by
Lemma 7.8(ii), both $T_{2j-1}=R_{2j}-R_{2j+1}=R_{2j}-R_1$ and
$T_{2j}=S_{2j-1}-S_{2j}$. And from $S_{2j-1}+S_{2j}=R_1+R_{2j}$ and
$S_{2j-1}\geq R_1$ we get $S_{2j}\leq R_{2j}$. \qed

\blm Suppose that the pair $M,N$ satisfies the condition R1 and the
conditions (1) and (2) of Theorem 7.5. Then for every
$1\leq i\leq m-2$ such that $R_i=R_{i+2}$ we have
$(R_{i+1}-R_i)/2\equiv e\pmod 2$.
\elm
\pf Note that $R_i=R_{i+2}$ implies $R_{i+2}-R_{i+1}=-(R_{i+1}-R_i)$,
so $(R_{i+1}-R_i)/2\equiv e\pmod 2$ is equivalent to
$(R_{i+2}-R_{i+1})/2\equiv e\pmod 2$.

If $S_i\geq R_i=R_{i+2}$, then we are in one of the two cases of the
condition (2) of Theorem 7.5. Suppose first that we are in the case
(b), i.e. $S_i=R_{i+2}=R_i$ and
$-2e\in\{ R_{i+2}-R_{i+1},S_{i+1}-S_i\}$. If $R_{i+2}-R_{i+1}=-2e$
then $(R_{i+2}-R_{i+1})/2=-e\equiv e\pmod 2$ and we are done. If
$S_{i+1}-S_i=-2e$, then, by Lemma 7.1(i), we have
$R_i+R_{i+1}\leq S_i+S_{i+1}$. By subtracting $2R_i=2S_i$, we get
$R_{i+1}-R_i\leq S_{i+1}-S_i=-2e$. This implies $R_{i+1}-R_i=-2e$, so
$(R_{i+1}-R_i)/2=-e\equiv e\pmod 2$. If we are in the case (a), then
$(R_{i+2}-R_{i+1})/2\equiv e\pmod 2$ and we are done. 

If $S_i<R_i$ then, by Lemma 7.1(ii), $R_{i+1}<S_{i-1}$, so the index
$i-1$ is in the case (a) of condition (2). Thus  
$(R_{i+1}-R_i)/2=-e\equiv e\pmod 2$. \qed

We now start our proof. It will be done in two steps.
\medskip

{\bf Step 1.} Reducing to the case $\nnn M=\nnn N$.

We assume that $\nnn M\supset\nnn N$, i.e. $R_1<S_1$. By Lemma 6.2, we
have $\theta (X(M/N))=\theta (X(M'/N))$ if $M'$ is a lattice
and $\theta (X(M/N))=\fff$ otherwise.

\blm If $M'$ is not a lattice, then Theorems 7.3 and 7.5 are true.
\elm
\pf Since $M'$ is not a lattice, we have $\theta (X(M/N))=\fff$, so
Theorem 7.3(ii) and (iii) are trivial, (For (ii) see the equivalent
statement from Remark 7.4(2).) So we are left with proving Theorem
7.3(i) and Theorem 7.5. 

If $M'$ is not a lattice, then, by Lemma 4.14, $R_2-R_1=-2e$ and
either $d(-a_{1,2})=\infty$ or $R_1=R_3$. We have two cases.

{\bf 1.} $R_3\geq S_1$. Then $R_3\geq S_1>R_1$, so we musty have
$d(-a_{1,2})=\infty$, i.e. $a_{1,2}=-1$ in $\fff/\ffs$.

For Theorem 7.3, note that the product in (i) contains the factor
$\bar G(a_{1,2},R_2-S_1)=\bar G(-1,R_2-S_1)$. But $R_2-S_1<-2e$ so
$\bar G(-1,R_2-S_1)=\N (1)=\fff$. Thus the product is $\fff$ and
Theorem 7.3(i) states that $\theta (X(M/N))=\fff$, which is true. 

For Theorem 7.5, note that $T_1=R_2-\min\{ S_1,R_3\} =R_2-S_1$. Then,
same as above,
$\bar G(a_{1,2},T_1)=\bar G(a_{1,2},R_2-S_1)=\fff\not\subseteq\ooo\ffs$.
Thus condition (3) of Theorem 7.5 fails, so Theorem 7.5 states that
$\theta (X(M/N))\not\subseteq\ooo\ffs$, which is true.

{\bf 2.} $R_3<S_1$. Then the pair $M,N$ doesn't have property A and so
Theorem 7.3(i) is vacuous. Suppose that the pair $M,N$ satisfies
condition (2) of Theorem 7.5. Then $R_3<S_1$ implies that
$R_2+R_3=S_1+S_2$. Together with $R_3<S_1$, this implies
$R_2>S_2$. Since also $R_1\leq R_3<S_1$, we have
$-2e=R_2-R_1>S_2-S_1$, which is impossible. Thus condition (ii) fails
and so Theorem 7.5 states that $\theta (X(M/N))\not\subseteq\ooo\ffs$,
which is true. \qed

From now on we assune that $M'$ is a lattice. We write
$M'\cong{\prec a_1',\ldots,a_m'\succ}$ relative to a good BONG, and
$R_i(M')=R_i'$. By Lemma 6.2, we have
$\theta (X(M/N))=\theta (X(M'/N))$.

For convenience, we denote by 7.3(i)-(iii) the statements of Theorem
7.3 and by 7.5(1)-(4) the conditions of Theorem 7.5 for the pair of
lattices $M,N$. For the pair $M',N$ we denote them by 7.3(i')-(iii')
and 7.5(1')-(4'), respectively.
\blm (i) The pair $M,N$ satisfies R1 iff the pair $M',N$ does so.

(ii) The pair $M,N$ has property A iff the pair $M',N$ does so.

If the pair $M,N$ satisfies condition R1 and
$\theta (O^+(M))\subseteq\ooo\ffs$ then:

(iii) If $M'$ is in the case (v) of Lemma 4.14, then the hypothesis of
Lemma 7.11(ii) holds for both pairs $M,N$ and $M',N$.

(iv) 7.5(1) and (2) hold iff 7.5(1') and (2') hold.
\elm
\pf (i) If $M$ is not in the case (v) of Lemma 4.14, then $R_i=R'_i$ for
$i\geq 3$ so at every $i\geq 3$ the condition R1 for the pairs $M,N$
and $M',N$ is the same. For $i=1,2$, the condition R1 holds for both
$M,N$ and $M',N$ by Lemma 3.1(i) and (ii). If $M$ is in the case (v) 
of Lemma then $R_i=R_i'$ for $i>2k+1$, so at every $i>2k+1$ the
condition R1 for the pairs $M,N$ and $M',N$ is the same. If
$1\leq i\leq 2k+1$ is odd, then $S_i\geq S_1\geq R_1'=R_i'>R_i$, so
condition R1 at $i$ holds for both pairs $M,N$ and $M',N$. If
$1<i<2k+1$ is even, then $R_{i-1}'=R_{i+1}'=R_{i-1}+2=R_{i+1}+2$ and
$R_i'=R_i-2$, so
$R_{i-1}+R_i=R_{i-1}'+R_i'=R_i+R_{i+1}=R_i'+R_{i+1}'$. If the pair
$M,N$ satisfies the condition R1, then, by Lemma 7.1(i),
$R_i'+R_{i+1}'=R_{i-1}+R_i\leq S_{i-1}+S_i$, so the condition R1 at
index $i$ holds for the pair $M',N$. Similarly, if the pair $M',N$
satisfies the condition R1, then
$R_i+R_{i+1}=R_{i-1}'+R_i'\leq S_{i-1}+S_i$, so the condition R1 at
index $i$ holds for the pair $M,N$. 

(ii) If $M'$ is not in the case (v) of Lemma 4.14, then
$R_i'=R_i$ for $i\geq 3$, so property A for $M,N$ and for $M',N$ is
the same. If $M'$ is in the case (v) of Lemma 4.14, then $S_1>R_1=R_3$
and $S_1\geq R'_1=R'_3$, so condition A fails for both $M,N$ and
$M',N$.

(iii) Since $\theta (O^+(M))\subseteq\ooo\ffs$, $M'$ cannot be in the
case (iv) of Lemma 4.14. Then in each of the remaining cases,
(i), (iii) and (v) have $R_i\equiv R'_i\pmod 2$ for $1\leq i\leq m$.
Thus 7.5(1) and 7.5(1') are equivalent. Assume now that 7.5(1) and
7.5(1') hold. 

If $M'$ is in the case (i) of Lemma 4.14, then $R_i=R'_i$ for
$i\geq 3$, so 7.5(2) and 7.5(2') are the same for $i\geq 2$. And at
$i=1$, since $R_2-R_1=R'_2-R'_1=-2e$, by Lemma 7.11(i), both 7.5(2)
and 7.5(2') are equivalent to $S_2-S_1=-2e$.

If $M'$ is in the case (iii) of Lemma 4.14, then $R_i=R'_i$ for
$i\geq 2$, so 7.5(2) and 7.5(2') are the same. 

If $M'$ is in the case (v) of Lemma 4.14, then for $i\geq 2k+2$ we
have $R_i=R'_i$. Hence conditions 7.5(2) and 7.5(2') are the same at
indices $i\geq 2k+1$. For $i\leq 2k$ we note that $R_1=R_{2k+1}$,
$R'_1=R'_{2k+1}$, $R_1<S_1$ and
$R'_1-R'_{2k}=(R_1+2)-(R'_{2k}-2)=R_{2k+1}-R_{2k}+4\geq -2e+4>-2e$, so
Lemma 7.11(ii) applies to both pairs $M,N$ and $M',N$. Then 7.5(2) and
7.5(2') for indices $i\leq 2k$ are equivalent to
$R_1+R_{2j}=S_{2j-1}+S_{2j}$ and $R'_1+R'_{2j}=S_{2j-1}+S_{2j}$,
respectively, for $1\leq j\leq k$. But
$R'_1+R'_{2j}=(R_1+2)+(R_{2j}-2)=R_1+R_{2j}$, so the two conditions
are the same. \qed

\blm If the pair $M,N$ satisfies property A and condition R1, then $M$
is in one of the cases (i), (iii) or (iv) of Lemma 4.14 and we have:

(i) For every $1\leq i\leq m$, except when $M$ is in the case (iv) of
Lemma 4.14 and $i=1$, we have
$R_1+\cdots +R_i\equiv R'_1+\cdots +R'_i\pmod 2$ and
$a_{1,i}=a'_{1,i}$ in $\fff/\ffs$.

(ii) $\tilde G(M/N)=\tilde G(M',N)$.
\elm
\pf By Corollary 7.2, the lattice $M$ satisfies property A, so we
cannot have $R_1=R_3$. Therefore we cannot be in the case (v) of Lemma
4.14.

(i) If $M$ is in the case (i) or (iii) of Lemma 4.14, then $a_i=a'_i$
in $\fff/\ffs$ for all $1\leq i\leq m$, so $a_{1,i}=a'_{1,i}$ in
$\fff/\ffs$ for $1\leq i\leq m$. In the case (iv) note that
$a'_{1,2}=\pi^2\eta^2a_{1,2}=a_{1,2}$ in $\fff/\ffs$ and
$a_i=a'_i$ for $i\geq 3$. So in this case $a_{1,i}=a'_{1,i}$ in
$\fff/\ffs$ holds for $2\leq i\leq m$.

Since $a_{1,i}$ and $a'_{1,i}$ differ by a square factor, we have
$\ord a_{1,i}\equiv\ord a'_{1,i}\pmod 2$, i.e.
$R_1+\cdots +R_i\equiv R'_1+\cdots +R'_i\pmod 2$. 

(ii) We have $R_i=R'_i$ for $i\geq 3$, which implies that
$[(R_{i+2}-S_i)/2]=[(R'_{i+2}-S_i)/2]$ for $1\leq i\leq m-2$. By
taking minima, we get $\beta =\beta'$.

By (i), for $1\leq i\leq m-2$ we have $a_{1,i+1}=a'_{1,i+1}$ in
$\fff/\ffs$. If moreover $i\geq 2$ or $M$ is in the case (iii) of
Lemma 4.14, then also $R_{i+1}=R'_{i+1}$ and we get
$\bar G (a_{1,i+1}b_{1,i-1},R_{i+1}-S_i)
=\bar G (a'_{1,i+1}b_{1,i-1},R'_{i+1}-S_i)$. The 
same relation also holds in the remaining cases, when $i=1$ and $M$
is in the case (i) or (iv) of Lemma 4.14, when it writes as
$\bar G (a_{1,2},R_2-S_1)=\bar G (a'_{1,2},R'_2-S_1)$. We still have
$a_{1,2}=a'_{1,2}$ in $\fff/\ffs$, but this time $R_2<R'_2$. In both
cases (i) and (iv), by Lemma 4.8(ii), we have
$\bar G(a_{1,2},R_2-R_1)=G(a_2/a_1)=\N (-a_{1,2})$. Then, by Lemma
2.12(i) and (iv), for every $R_2-R_1\geq R$ we have
$\N (-a_{1,2})=\bar G(a_{1,2},R_2-R_1)\subseteq\bar G(a_{1,2},R)
\subseteq\N (-a_{1,2})$, so $\bar G(a_{1,2},R)=\N (-a_{1,2})$. But in
both cases we have $R_2-R_1=R'_2-R'_1\geq R_2'-S_1>R_2-S_1$. It
follows that $\bar G(a_{1,2},R_2-S_1)=\N (-a_{1,2})$ and 
$\bar G(a'_{1,2},R_2-S_1)=\bar G(a_{1,2},R_2-S_1)=\N (-a_{1,2})$, so
we have the claim equality.

In conclusion,
$\ups\beta\prod_{i=1}^{m-2}\bar G(a_{1,i+1}b_{1,i-1},R_{i+2}-S_i)=
\ups{\beta'}\prod_{i=1}^{m-2}\bar G(a'_{1,i+1}b_{1,i-1},R'_{i+2}-S_i)$, 
i.e. $G(M/N)=G(M'/N)=:G$. Moreover, if $M$ is in one of the cases (i)
and (iv) of Lemma 4.14, then
$G\supseteq\bar G (a_{1,2},R_2-S_1)=\N (-a_{1,2})$.

By Remark 7.2(1), we have
$\tilde G(M/N)=G\prod_{i=2}^{m-1}\bar G(a_{i,i+1},R_{i+1}-R_i)$ and
$\tilde G(M'/N)=G\prod_{i=2}^{m-1}\bar G(a'_{i,i+1},R'_{i+1}-R'_i)$.
To prove that $\tilde G(M/N)=\tilde G(M'/N)$ we will show that
$G\bar G(a_{i,i+1}R_{i+1}-R_i)=G\bar G(a'_{i,i+1},R'_{i+1}-R'_i)$ for
$2\leq i\leq m-1$. If $M$ is in the case (iii) of Lemma 4.14, then
$a_i=a'_i$ for $2\leq i\leq m$, so for $2\leq i\leq m-1$ we have
$G(a_{i+1}/a_i)=G(a'_{i+1}/a'_i)$, i.e.
$\bar G(a_{i,i+1},R_{i+1}-R_i)=\bar G(a'_{i,i+1},R'_{i+1}-R'_i)$, so
our statement follows trivially. If $M$ is in one of the cases (i) and
(iv), then $a_i=a'_i$ only holds for $i\geq 3$ so
$\bar G(a_{i,i+1},R_{i+1}-R_i)=\bar G(a'_{i,i+1},R'_{i+1}-R'_i)$ holds
for $3\leq i\leq m-1$.

So we are left to prove that
$G\bar G(a_{2,3},R_3-R_2)=G\bar G(a'_{2,3},R'_3-R'_2)$ in the cases (i)
and (iv) of Lemma 4.14. Recall that $G\supseteq\N (c)$, where $c$ is
the value of $-a_{1,2}$ in $\fff/\ffs$. If $G=\fff$ the statement we
want to prove is trivial. So we may assume that $G=\N (c)$. Then
$G\bar G(a_{2,3},R_3-R_2)=\N (c)$ if
$\bar G(a_{2,3},R_3-R_2)\subseteq\ N(c)$ and $=\fff$
otherwise. Similarly for
$\bar G(a_{2,3},R'_3-R'_2)=\bar G(a'_{2,3},R_3-R'_2)$ so we must prove
that $\bar G(a_{2,3},R_3-R_2)\subseteq\ N(c)$ iff
$\bar G(a'_{2,3},R_3-R'_2)\subseteq\ N(c)$.

If we are in the case (i) of Lemma 4.14, then
$\N (c)=\N (\Delta )=\ooo\ffs$. We have $R_3>S_1\geq R'_1=R'_2+2e$, so
$R_3-R_2>R_3-R'_2>2e$. Then, by Lemma 2.17,
$\bar G(a_{2,3},R_3-R_2)\subseteq\ooo\ffs$ and
$\bar G(a'_{2,3},R_3-R'_2)\subseteq\ooo\ffs$ are equivalent to
$a_{2,3}\in\ooo\ffs$ and $a'_{2,3}\in\ooo\ffs$, respectively. But
$a'_{2,3}=\pi^2a_{2,3}$, so the two are equivalent.

If we are in the case (iv) of Lemma 4.14, then
$R'_2-R'_1=R_2-R_1=1-d(-a_{1,2})=1-d(c)$. We have
$\ups\beta\subseteq G=\N (c)$, so, by Lemma 1.2(i),
$\beta +d(c)>2e$. Thus $[(R_3-S_1)/2]\geq\beta\geq 2e+1-d(c)$, which
implies $R_3-S_1\geq 4e+2-2d(c)$. Since 
$R'_2=R'_1+1-d(c)\leq S_1+1-d(c)$, we have
$R_3-R_2>R_3-R'_2\geq R_3-S_1-1+d(c)\geq 4e+2-2d(c)-1+d(c)=4e+1-d(c)$.
It follows that $R_3-R_2-2e>R_3-R'_2-2e\geq 2e+1-d(c)$, which, by
Lemma 1.2(i), implies
$\ups{R_3-R_2-2e}\subseteq\ups{R_3-R'_2-2e}\subseteq\N (c)$. By
Corollary 2.13(iv), we have
$a_{1,2}\in\bar G(a_{2,3},R_3-R_2)\subseteq 
\langle a_{2,3}\rangle\ups{R_3-R_2-2e}$. Since
$\ups{R_3-R_2-2e}\subseteq\N (c)$, this implies that
$\bar G(a_{2,3},R_3-R_2)\subseteq\N (c)$ iff $a_{2,3}\in\N (c)$.
Similarly, $\bar G(a'_{2,3},R_3-R'_2)\subseteq\N (c)$ iff
$a'_{2,3}\in\N (c)$. But, using the notations of Lemma 4.14(iii),
$a'_{2,3}=\pi\eta a_{2,3}$ and $\pi\eta\in\N (c)$. Hence
$a_{2,3}\in\N (c)$ iff $a'_{2,3}\in\N (c)$, which concludes the
proof. \qed

\bpr If $M'$ is a lattice, then Theorem 7.3 holds for the pair of
lattices $M,N$ iff it holds for the pair $M',N$.
\epr
\pf By Lemma 7.14(i) and (ii), the pair $M,N$ satisfies the condition
R1 and has property A iff the pair $M,N$ does so. If this happens,
then 7.3(i) and 7.3(i') state that $\theta (X(M/N))=\tilde G(M/N)$ and
$\theta (X(M'/N))=\tilde G(M'/N)$, respectively. Since
$X(M/N)=X(M'/N)$ and, by Lemma 7.15(ii),
$\tilde G(M/N)=\tilde G(M'/N)$, the two statements are equivalent.

By Remark 7.4(2), 7.3(ii) writes as 
$\N (-a_{1,i+1}b_{1,i-1})\subseteq\theta (X(M/N))$ for every
$1\leq i\leq\min\{ n,m-1\}$ such that
$R_1+\cdots +R_i\not\equiv S_1+\cdots +S_i\pmod 2$. Similarly for
7.3(ii'). By Lemma 7.15(i), with exception of the case when $i=1$ and
$M$ is in the case (iv) of Lemma 4.14, we have
$R_1+\cdots +R_i\equiv R'_1+\cdots +R'_i\pmod 2$ and
$-a_{1,i+1}b_{1,i-1}=-a'_{1,i+1}b_{1,i-1}$ in $\fff/\ffs$. Since also
$X(M/N)=X(M'/N)$, the conditions for 7.3(ii) and 7.3(ii') are
equivalent. In the exceptional case we have
$\N (-a_{1,2})=\N (a'_{1,2})\subseteq\theta (X(M/N))=\theta (X(M'/N)$
regardless of the parity of $R_1$ and $R'_1$. So at $i=1$ the
condition for both 7.3(ii) and 7.3(ii') hold unconditionally.

For 7.3(iii) we have three cases: If the the pair $M,N$ satisfies
condition R1 and has property A, then 7.3(iii) is vacuous. If it
doesn't satisfy condition R1, then 7.3(iii) states that
$\theta (X(M/N))=\fff$. (This includes the case when both R1 and
property A fail, since $\theta (X(M/N))=\fff$ implies
$\theta (X(M/N))\supseteq\ooo\ffs$.) And if R1 holds, but property A
not, it states that $\theta (X(M/N))\supseteq\ooo\ffs$.

By Lemma 7.14(i) and (ii), the pair $M',N$ is in the same of the
three cases above as the pair $M,N$. Since also
$\theta (X(M/N))=\theta (X(M'/N))$, conditions 7.3(iii) and 7.3(iii')
are equivalent. \qed

\blm Suppose that $M$ is in the case (v) of Lemma 4.14, the pair of
lattices $M,N$ satisfies condition R1, 7.5(i) and 7.5(ii) and
$\theta (O^+(N))\subseteq\ooo\ffs$. Then for every $1\leq j\leq k$ we
have $\bar G((-1)^{j-1}b_{1,2j},R'_{2j}-R'_1)\subseteq\ooo\ffs$.
\elm
\pf For every $1\leq h\leq j$ we have
$\bar G(b_{2h-1,2h},S_{2h}-S_{2h-1})=G(b_{2h}/b_{2h-1})\subseteq\theta
(O^+(N))\subseteq\ooo\ffs$. We also have $S_{2h-1}\geq S_1\geq R'_1$
and, by Lemma 7.11(ii), $S_{2h-1}+S_{2h}=R_1+R_{2h}=R'_1+R'_{2h}$. By
subtracting the inequality $2S_{2h-1}\geq 2R'_1$, we get
$S_{2h}-S_{2h-1}\leq R'_{2h}-R'_1\leq R'_{2j}-R'_1$. Since also
$R'_{2j}-R'_1=R'_{2j}-R'_{2j-1}\geq -2e$, by Corollary 2.18, we get
$\bar G((-1)^{j-1}b_{1,2j},R'_{2j}-R'_1)\subseteq\ooo\ffs$. \qed

\bpr Theorem 7.5 holds for the pair of lattices $M,N$ iff it holds for
the pair $M',N$. 
\epr
\pf By Lemma 7.14(i), conditions R1 for the pairs $M,N$ and $M',N$ are
equivalent. If they both fail, then Theorem 7.5 is vacuous for both
pairs. So we will assume that R1 holds for both $M,N$ and $M',N$. 

By Lemma 6.2, $O^+(M)\subseteq O^+(M')$. Hence if
$\theta (O^+(M))\not\subseteq\ooo\ffs$ then also
$\theta (O^+(M'))\not\subseteq\ooo\ffs$ and so both 7.5(iv) and
7.5(iv') fail and we are done. Therefore we will assume that 
$\theta (O^+(M))\subseteq\ooo\ffs$. We will also assume that
$\theta (O^+(N))\subseteq\ooo\ffs$, since this condition is common to
7.5(iv) and 7.5(iv'). Consequently, Lemmas 7.14(iii) and 7.17 apply.

By Lemma 7.14(iii), 7.5(i) and (ii) are equivalent to 7.5(i') and
(ii'), so we may assume that these conditions hold. We must prove that
7.5(iii) and (iv) hold iff 7.5(iii') and (iv') hold.

Since $\theta (O^+(M))\subseteq\ooo\ffs$, $M'$ cannot be in the case
(iv) of Lemma 4.14. If we are in the case (i) or (iii), then
$a_i=a'_i$ in $\fff/\ffs$ for all $1\leq i\leq m$. The same may be
assumed in the case (v), by choosing as BONGs for $M$ and $M'$ those
from Lemma 4.14. (Recall that, by Lemma 7.10, Theorem 7.5 is
independent of the choice of the good BONG of $M$.) Consequently,
$a_{1,i}=a'_{1,i}$ in $\fff/\ffs$ for $1\leq i\leq m$.

For the equivalence between 7.5(iii) and 7.5(iii') we prove that for
every $1\leq i\leq m-1$ we have
$\bar G(a_{1,i+1}b_{1,i-1},T_i)\subseteq\ooo\ffs$ iff 
$\bar G(a_{1,i+1}b_{1,i-1},T'_i)=\bar G(a'_{1,i+1}b_{1,i-1},T'_i)
\subseteq\ooo\ffs$. In most cases $T_i=T'_i$, so the statement is
trivial.

If we are in the case (iii) of Lemma 4.14, then $R_i=R'_i$ for
$i\geq 2$, so $T_i=T'_i$ for $i\geq 1$ and we are done.

If we are in the case (i), then $R_i=R'_i$ for $i\geq 3$,
so $T_i=T'_i$ for $i\geq 2$. In the remaining case, $i=1$, we note
that $a_{1,2}=-\Delta$ in $\fff/\ffs$, so, by Lemma 2.12(iv), both
$\bar G(a_{1,2},T_1)$ and $\bar G(a_{1,2},T'_1)$ are contained in
$\N (-a_{1,2})=\N (\Delta )=\ooo\ffs$.

If we are in the case (v), then for every $1\leq h\leq k$ we have
$a_{2h-1}=\pi^{R_1}\varepsilon$ and
$a_{2h}=-\pi^{R_{2h}}\varepsilon =\pi^{R_1}\varepsilon$ in
$\fff/\ffs$. (We have $R_1\equiv R_{2h}\pmod 2$.) Hence for
$1\leq\ell\leq 2k$ in $\fff/\ffs$ we have $a_\ell =-a_{\ell+1}$, i.e.
$a_{\ell,\ell +1}=-1$. Consequently, $a_{1,2j}=(-1)^j$ for
$1\leq j\leq k$. 

We have $R_i=R'_i$ for $i\geq 2k+2$, so $T_i=T'_i$ for $i\geq 2k+1$.
For $i=2j$, with $1\leq j\leq k$, we have, by Lemma 7.14(iii),
$T_{2j}=T'_{2j}=S_{2j-1}-S_{2j}$ and we are done.

So we are left with the case when $i=2j-1$, with $1\leq j\leq k$. Note
that in $\fff/\ffs$ we have
$a_{1,2j}b_{1,2j-2}=(-1)^jb_{1,2j-2}=(-1)^{j-2}b_{2j-2}$ and, by Lemma
7.14(iii), $T_{2j-1}>T'_{2j-1}=R'_{2j}-R'_1\geq R'_{2j-2}-R'_1$.
Hence, by Lemmas 2.6(i) and 7.17, we get
$\bar G((-1)^{j-2}b_{2j-2},T_{2j-1},\subseteq
\bar G((-1)^{j-2}b_{2j-2},T'_{2j-1})\subseteq
\bar G((-1)^{j-2}b_{2j-2},R'_{2j-2}-R'_i)\ooo\ffs$. Hence both
7.5(3) and 7.5 (3') hold unconditionally at $i=2j-1$.

We now assume that 7.5(i)-(iii) and 7.5(i')-(iii') hold and we prove
the equivalence between 7.5(iv) and 7.5(iv'). Recall that we are
under the assumption that
$\theta (O^+(M)),\theta (O^+(N))\subseteq\ooo\ffs$ i.e. 7.5(iv)
holds. So we must prove that $\theta (O^+(M'))\subseteq\ooo\ffs$, so
that 7.5(iv') also holds. Since we assume that conditions R1, 7.5(i')
and 7.5(ii') hold, by Lemma 7.12, the condition that
$(R'_{i+1}-R'_i)/2\equiv e\pmod 2$ whenever $R'_i=R'_{i+2}$ from
Theorem 1.9 is satified. So we are left to prove that
$G(a'_{i+1}/a'_i)\subseteq\ooo\ffs$ for $1\leq i\leq m-1$. Since
$\theta (O^+(M))\subseteq\ooo\ffs$, we have
$G(a_{i+1}/a_i)\subseteq\ooo\ffs$. In most cases
$a'_{i+1}/a'_i=a_{i+1}/a_i$, so
$G(a'_{i+1}/a'_i)=G(a_{i+1}/a_i)\subseteq\ooo\ffs$. Note that
in $\fff/\ffs$ we have $a'_{i,i+1}=a_{i,i+1}$, so
$G(a'_{i+1}/a'_i)=\bar G(a'_{i,i+1},R'_{i+1}-R'_i,)
=\bar G(,a_{i,i+1},R'_{i+1}-R'_i)$.

If $M'$ is in the case (iii) of Lemma 3,14, then for $i\geq 2$ we have
$a'_{i+1}/a'_i=a_{i+1}/a_i$, so
$G(a'_{i+1}/a'_i)\subseteq\ooo\ffs$. For $i=1$, by 7.5(3'), we have
$\bar G(a'_{1,2},T'_1)\subseteq\ooo\ffs$. By
Lemma 2.6(i), in order to 
prove that $G(a'_2/a'_1)=\bar G(a'_{1,2},R'_2-R'_1)\subseteq\ooo\ffs$,
it is enough to show that $R'_2-R'_1\geq T'_1=R'_2-\min\{ S_1,R'_3\}$.
But this follows from $R'_1\leq S_1$ and $R'_1\leq R'_3$, so
$R'_1\leq\min\{ S_1,R'_3\}$. 

If we are in the case (i) of Lemma 4.14, then
$a'_{i+1}/a'_i=a_{i+1}/a_i$ for $i\geq 3$ and same happens at
$i=1$. (We have $a'_2/a'_1=(\pi^2a_2)/(\pi^2a_1)$.) For $i=2$ we note
that $R'_3-R'_2\geq R'_1-R'_2=2e$ and, by 7.5(1'),
$\ord a'_{2,3}=R'_2+R'_3$ is even. Then, by Lemma 2.17(i),
$G(a'_3/a'_2)=\bar G(a'_{2,3},R'_3-R'_2)\subseteq\ooo\ffs$.

If $M'$ is in the case (v) of Lemma 4.14, then
$a'_{i+1}/a'_i=a_{i+1}/a_i$ for $i\geq 2k+2$. For $i\leq 2k$ we have
$a'_{i+1}/a'_i=-\pi^{R'_{i+1}-R'_i}=-1$ in $\fff/\ffs$. (By 7.5(1'),
$R'_{i+1}-R'_i$ is even.) We also have $R'_{i+1}-R'_i\geq -2e$, so, by
Lemma 2.6(i), we have
$G(a'_{i+1}/a'_i)=\bar G(-1,R'_{i+1}-R'_i)\subseteq
\bar G(-1,-2e)=\ooo\ffs$. So we are left with the case $i=2k+1$. 

So we are left to prove that $G(a'_{2k+2}/a'_{2k+1})
=\bar G(a_{2k+1,2k+2},R'_{2k+2}-R'_{2k+1})\subseteq\ooo\ffs$.  In
$\fff/\ffs$  we have $a_{2k+1,2k+2}=-ab$, where
$a=a_{1,2k+2}b_{1,2k}$ and $b=-a_{1,2k}b_{1,k}=(-1)^{k-1}b_{1,k}$. By
7.5(3), we have $\bar G(a,T_{2k+1})\subseteq\ooo\ffs$ and, by Lemma
7.17, $\bar G(b,R'_{2k}-R'_1)\subseteq\ooo\ffs$. We have
$R'_{2k+2}-R'_{2k+1}\geq -2e$ so, by the case $s=2$ of Corollary 2.18,
in order to prove that
$G(a'_{2k+2}/a'_{2k+1})=\bar G(-ab,R'_{2k+2}-R'_{2k+1})\subseteq\ooo\ffs$,
it suffices to show that
$R'_{2k+2}-R'_{2k+1}=R_{2k+2}-R'_1\geq T'_{2k+1},R'_{2k}-R'_1$. By
Lemma 7.11(ii), we have $R_{2k+2}\geq R_{2k}>S_{2k}$, so
$R_{2k+2}=\max\{ R_{2k+2},S_{2k}\}$. We also have 
$R'_1\leq S_1\leq S_{2k+1}$ and $R'_1\leq R'_{2k+3}=R_{2k+3}$ so
$R'_1\leq\min\{ S_{2k+1},R_{2k+3}\}$. It follows that
$R_{2k+2}-R'_1\geq\max\{ R_{2k+2},S_{2k}\} -\min\{ S_{2k+1},R_{2k+3}\}
=T_{2k+1}$. And, since $R_{2k+2}=R'_{2k+2}\geq R'_{2k}$, we also have
$R_{2k+2}-R'_1\geq R'_{2k}-R'_1$. This concludes the proof. \qed

{\bf Step 2. Induction on $n$.} In Step 1. we reduced to the case when
$R_1=S_1$. We now assume that $R_1=S_1$ and we use Lemma 6.3 to reduce
the ranks of $M$ and $N$. First we prove that if the condition
of Lemma is not fulfilled, then Theorems 7.3 and 7.5 hold.

\blm If $M$ and $N$ are not in one of the cases (a), (b) or (c) of
Lemma 6.3, then Theorems 7.3 and 7.5 hold.
\elm
\pf By Lemma 6.4, $\theta (X(M/N))=\fff$. Since (a) doesn't hold, we
have $R_3=R_1=S_1$ so the pair $M,N$ doesn't have property A. Hence
the statements (i) and (ii) of Theorem 7.3 are vacuous. Since
$\theta (X(M/N))=\fff$, statement (iii) holds trivially.

Suppose that condition (2) of Theorem 7.5 holds at $i=1$. We have
$R_3=S_1$. So one of the two cases of (2) must hold. If we are in the
case 2(a), then $R_2+R_3=S_1+S_2$. Since also $R_3=S_1$, we get
$R_2=S_2$, i.e. we are in the case (b) of Lemma 6.3. Contradiction.
Suppose now that we are in the case (2)(b), i.e.
$-2e\in\{ R_3-R_2,S_2-S_1\}$. We cannot have $R_3-R_2=-2e$, 
since this is case (c) of Lemma 6.3. So we are left with the case when
$S_2-S_1=-2e$. By Lemma 3.1(iii), we have $R_1+R_2\leq S_1+S_2$. After
subtracting $2R_1=2S_1$ we get $-2e\leq R_2-R_1\leq S_2-S_1=-2e$, so
we must have equality. Thus $R_2-R_1=S_2-S_1=S_2-R_1$, so
$R_2=S_2$. Again, we get a contradiction.

In conclusion, condition (ii) of Theorem 7.5 fails and so Theorem 7.5
states that $\theta (X(M/N))\not\subseteq\ooo\ffs$, which is true. (We
have $\theta (X(M/N))=\fff$.) \qed

For now on we assume that the pair $M,N$ is in one of the cases of
Lemma 6.3. By Lemma 6.3(i), we have $b_1=\eta a_1$ for some
$\eta\in g(a_2/a_1)$ and
$M\cong{\prec\eta a_1,\eta a_2,a_3,\ldots,a_m\succ}$ relative to a
good BONG beginning with $y_1$. If $M^*=pr_{y_1^\perp}M$ and
$N^*=pr_{y_1^\perp}N$, then $M^*\cong{\prec a_2^*,\ldots,a_m^*\succ}$,
where $a_2^*=\eta a_2$ and $a_i^*=a_i$ for $i\geq 3$, and
$N^*\cong{\prec b_2,\ldots,b_n\succ}$. By Lemma 6.3(ii),
$\theta (X(M/N))=\theta (X(M^*/N^*))\theta (O^+(M))$,
$\theta (O^+(M^*))\subseteq\theta (O^+(M))$ and
$\theta (O^+(N^*))\subseteq\theta (O^+(N))$. 

We have $R_i(M^*)=R_{i+1}$ and $R_i(N^*)=S_{i+1}$, i.e. the analogues
of the sequences $R_1,\ldots,R_m$ and $S_1\ldots,S_n$ corresponding to
$M^*$ and $N^*$ are $R_2,\ldots,R_m$ and $S_2\ldots,S_n$. The
analogues of the products $a_{1,i+1}b_{1,i-1}$ for $1\leq i\leq m-1$
are $a_{2,i+1}^*b_{2,i-1}=\eta a_{2,i+1}b_{1,i-1}$. (Recall,
$a_2^*=\eta a_2$ and $a_i^*=a_i$ for $i\geq 3$.) But $\eta =b_1/a_1$,
so in $\fff/\ffs$ we have $\eta =a_1b_1$, and so
$a_{2,i+1}^*b_{2,i-1}=a_1b_1a_{2,i+1}b_{2,i-1}=a_{1,i+1}b_{1,i-1}$.

We denote by $T_i^*$, with $2\leq i\leq\min\{ m-1,n\}$, the analogues
of $T_i$, with $1\leq i\leq\min\{ m-1,n\}$, corresponding to
$M^*,N^*$. We have
$T_i^*=\max\{ R_{i+1},S_{i-1}\} -\min\{ S_i,R_{i+2}\}$, where
$S_{i-1}$ is ignored if $i=2$ and $R_{i+2}$ is ignored if $i=m-1$.

We first assume that $n\geq 2$ and we prove the induction step
$n-1\to n$ by showing that if Theorems 7.3 and 7.5 are hold for
$M^*,N^*$, then they also hold for $M,N$. Same as for the proof of
Step 1, we denote by 7.3(i)-(iii) the statements of Theorem 7.3 and by
7.5(1)-(4) the condtions of Theorem 7.5 corresponding to the pair of
lattices $M,N$. For $M^*,N^*$ we denote them by 7.3(i$^*$)-(iii$^*$) and
7.5(1$^*$)-(4$^*$).

\blm (i) The pair $M,N$ satisfies condition R1 iff the pair $M^*,N^*$
does so or $n=1$.

(ii)  Assuming that the pair $M,N$ satisfies condition R1, it has
property A iff both $M$ and the pair $M^*,N^*$ have property A. (When
$n=1$ the condition on $M^*,N^*$ is ignored.)

(iii) If the pair $M,N$ satisfies condition R1 and it has property A,
then $\beta =\min\{\alpha,\beta^*\}$. (If $n=1$, then $\beta^*$ is
ignored, so $\beta =\alpha$.)

(iv) If the pair $M,N$ satisfies conditions R1 and 7.5(4), then 
7.5(1) holds unconditionally and, if $n>1$, so does 7.5(1$^*$). And
7.5(2) is equivalent to 7.5(2$^*$) or $n=1$. Also,
$T_1=R_2-R_1$ and $T_i=T_i^*$ for $2\leq i\leq\min\{ m-1,n\}$.
\elm
\pf (i) If $3\leq i\leq n$, then condition R1 at index $i$ for $M,N$
is the same as the condition R1 at index $i-1$ for $M^*,N^*$, namely
either $R_i\leq S_i$ or $i<m$ and $R_i+R_{i+1}\leq S_{i-1}+S_i$. We
are left with the condition R1 for $M,N$ at $i=1$ and $2$ and for
$M^*,N^*$ at $i=1$. But these always hold, by Lemma 3.1.

(ii) If the pair $M,N$ has property A  and $n>1$, then so does the
pair $M^*,N^*$. (If $R_{i+2}>S_i$ holds for $1\leq i\leq m-2$, it also
holds for $2\leq i\leq m-2$.) By Remark 7.6(1), $M$ too has property
A. Conversely, if $M^*,N^*$ has property A, then $R_{i+2}>S_i$ holds
for $2\leq i\leq m-2$. If also $M$ has property A, then $R_{i+2}>S_i$
also holds at $i=1$. (We have $R_3>R_1=S_1$.)

(iii) We may assume that $m\geq 3$ since otherwise
$\beta =\alpha =\beta^*=\infty$ and our statement becomes trivial. We
have $\beta =\min\{ [(R_{i+1}-S_i)/2]\mid 1\leq i\leq m-2\}$ and 
$\beta^*=\min\{ [(R_{i+1}-S_i)/2]\mid 2\leq i\leq m-2\}$, so
$\beta =\min\{ [(R_2-S_1)/2],\beta^*\} =\min\{ [(R_2-R_1)/2],\beta^*\}$. 
Together with $[(R_2-R_1)/2]\leq\alpha\leq\beta$, this implies
$\beta =\min\{\alpha,\beta^*\}$. (The inequalities
$[(R_2-R_1)/2]\leq\alpha\leq\beta$ follow from the definition of
$\alpha$ and Remark 7.4(1).)

(iv) By the remark following Theorem 1.9, since
$\theta (O^+(M)),\theta (O^+(N))\subseteq\ooo\ffs$, we have
$R_1\equiv\cdots\equiv R_2\pmod 2$ and
$S_1\equiv\cdots\equiv S_n\pmod 2$. Since also $R_1=S_1$, we get
7.5(1). If $n\geq 1$, then 7.5(1) implies
$R_2\equiv\cdots\equiv R_m\equiv S_2\equiv\cdots\equiv S_n\pmod 2$,
i.e. 7.5(1$^*$) holds.

If $n>1$, since $R_i(M)=R_{i-1}(M^*)=R_i$ for $2\leq i\leq m$ and
$R_i(N)=R_{i-1}(N^*)=S_i$ for $2\leq i\leq n$, condition 7.5(2$^*$) is
equivalent to 7.5(2) at indices $2\leq i\leq m-2$. (Explicitly, for
$2\leq i\leq m-2$ condition 7.5(2) at index $i$ is equivalent to
7.5(2$^*$) at index $i-1$.) To complete the proof, we show that 7.5(2)
holds unconditionally at $i=1$. If we are in the case (a) of Lemma
6.3, then $R_3>R_1=S_1$ so 7.5(2) at $i=1$ holds trivially. Suppose
that we are not in the case (a) of Lemma 6.3, i.e. that
$R_3=R_1=S_1$. If we are in the case (c), i.e. if $R_3-R_2=-2e$ then
$i=1$ is in the case 7.5(2)(b) and we are done. If we are in case
(b) of Lemma 6.3, i.e. $R_2=S_2$, then $R_2+R_3=S_1+S_2$. Since
$\theta (O^+(M))\subseteq\ooo\ffs$ and
$R_1=R_3$, by Theorem 1.9, we have $(R_2-R_1)/2\equiv 2\pmod 2$. But
$R_2-R_1=S_2-S_1=-(R_3-R_2)$. It follows that
$(R_3-R_2)/2\equiv (S_2-S_1)/2\equiv e\pmod 2$ and so $i=1$ is in the
second case of 7.5(2)(a). This concludes the proof.

We have $R_3\geq R_1=S_1$, so $T_1=R_2-\min\{ S_1,R_3\} =R_2-R_1$ and
$T_2=\max\{ R_3,S_1\}-\min\{ S_2,R_4\} =R_3-\min\{ S_2,R_4\} =T^*_2$.
And if $i\geq 3$, then $T_i=T^*_i$ by definition. \qed

\blm If $n=1$ and $m\leq 2$, then $X(M/N)=O^+(M)$. In particular, if
$m=n=1$, then $X(M/N)=\{1\}$, so $\theta (X(M/N))=\ffs$.

If $n=1$ and $m=3$, then
$\theta (X(M/N))=\N (-a_{1,3}b_1)\theta (O^+(M))$.
\elm
\pf If $m=n=1$, then $FM$ is unary, so
$O^+(M)\subseteq X(M/N)\subseteq O^+(FM)=\{ 1\}$ and our statement is
trivial.

If $m>n=1$, then , by Lemma 6.3(ii), $X(M/N)=O^+(FM^*)O^+(M)$. If
$m=2$, then $FM^*$ is unary, so $O^+(FM^*)=\{ 1\}$ and so
$X(M/N)=O^+(M)$. And if $m=3$, then
$\theta (X(M/N))=\theta (O^+(FM^*))\theta (O^+(M))$. But $FM^*$ is
binary of determinant $a^*_{2,3}=a_{1,3}b_1$, so
$\theta (O^+(FM^*))=\N (-a_{1,3}b_1)$. Hence the conclusion. \qed

\bpr Theorem 7.3 holds for the pair $M,N$ if it holds for the pair
$M^*,N^*$ or $n=1$. 
\epr
\pf If the pair $M,N$ doesn't satisfy R1 then 7.3(i) and 7.3(ii) are
vacuous. By Lemma 7.20(i), we have $n>2$ and the pair does $M^*,N^*$
doesn't satisfy R1 either so, by 7.3(iii$^*$),
$\theta (X(M^*/N^*))=\fff$. Since
$\theta (X(M^*/N^*))\subseteq\theta (X(M/N))$, we have
$\theta (X(M/N))=\fff$, so 7.3(iii) holds. Also 7.3(i) and (ii) are
vacuous when R1 is not satisfied, so we are done.

Assume now that the pair $M,N$ satisfies R1, but doesn't have property
A. Then again 7.3(i) and (ii) are vacuous. By Lemma 7.20(ii), either
$M$ doesn't have property A or $n>1$ and the pair $M^*,N^*$ doesn't
have property A. In the first case, by Proposition 1.6(iii), we have
$\ooo\ffs\subseteq\theta (O^+(M))\subseteq\theta (X(M/N))$ and in the
second case, by 7.3(iii$^*$), we have
$\ooo\ffs\subseteq\theta (X(M^*/N^*))\subseteq\theta (X(M/N))$. Hence
7.3(iii) holds.

Suppose now that the pair $M,N$ satisfies R1 and has property A. By
Lemma 7.20(i) and (ii), if $n>1$, same holds for $M^*,N^*$. Then, by
7.a(i$^*$), written as in Remark 7.4(1),
$\theta (X(M^*/N^*))=\tilde G(M^*/N^*)=G(M/N)\theta (O^+(M^*))$. Since
$\theta (O^+(M^*))\subseteq\theta (O^+(M))$, this implies that
$\theta (X(M/N))=\theta (X(M^*/N^*))\theta (O^+(M))=
G(M^*/N^*)\theta (O^+(M))$. Hence we must prove that
$\tilde G(M/N)=G(M^*/N^*)\theta (O^+(M))$. 

By Theorem 1.7, both $\ups\alpha$ and
$G(a_2/a_1)=\bar G(a_{1,2},R_2-R_1)=\bar G(a_{1,2},R_2-S_1)$ are
contained in $\theta (O^+(M))$. It follows that
$G(M^*/N^*)\theta (O^+(M))=G\theta (O^+(M))$, where
$G=\ups\alpha\bar G(a_{1,2},R_2-S_1)G(M^*/N^*)$.
We have
$$G(M^*/N^*)=\ups{\beta^*}
\prod_{i=2}^{m-1}\bar G(a^*_{2,i+1}b_{2,i-1},R_{i+1}-S_i)
=\ups{\beta^*}\prod_{i=2}^{m-1}\bar G(a_{1,i+1}b_{1,i-1},R_{i+1}-S_i).$$
But $\ups\alpha\ups{\beta^*}=\ups\beta$ (by Lemma 7.20(iii),
$\min\{\alpha,\beta^*\} =\beta$) and
$\bar G(a_{1,2},R_2-S_1 )
\prod_{i=2}^{m-1}\bar G(a_{1,i+1}b_{1,i-1},R_{i+1}-S_i)
=\prod_{i=1}^{m-1}\bar G(a_{1,i+1}b_{1,i-1},R_{i+1}-S_i)$.
It follows that
$G=\ups\beta\prod_{i=1}^{m-1}\bar G(a_{1,i+1}b_{1,i-1},R_{i+1}-S_i)
=G(M/N)$, which implies that
$G(M^*/N^*)\theta (O^+(M))=G(M/N)\theta (O^+(M))=\tilde G(M/N)$, so
7.3(i) holds.

Suppose now that $n=1$. We use Lemma 7.21. If $m=1$, 
$\theta (X(M/N))=\tilde G(M/N)=\ffs$. If $m=2$, then
$\theta (X(M/N))=\theta (O^+(M))$. We also 
have $\tilde G(M/N)=\bar G(a_{1,2},R_2-S_1)\theta (O^+(M))$. But
$\bar G(a_{1,2},R_2-S_1)=\bar G(a_{1,2},R_2-R_1)=G(a_2/a_1)=\theta
(O^+(M))$. Thus $\tilde G(M/N)=\theta (O^+(M))=\theta (X(M/N))$. And
if $m=3$, then $\theta (X(M/N)=\N (-a_{1,3}b_1)\theta (O^+(M))$ and
$G(M/N)=\ups\beta\bar G(a_{1,2},R_2-S_1)\N (-a_{1,3}b_1)$. As seen
above, $\bar G(a_{1,2},R_2-S_1)=G(a_2/a_1)$ and, by Lemma 7.20(iii)
$\beta =\alpha$, so $\ups\beta=\ups\alpha$. But, by Theorem 1.7,
$\ups\alpha,G(a_2/a_1)\subseteq\theta (O^+(M))$ so both these factors
are superfluous in the product $\tilde G(M/N)=G(M/N)\theta (O^+(M))$.
Thus $\tilde G(M/N)=\N (-a_{1,3}b_1)\theta (O^+(M))=\theta (X(M/N))$. 

Since $R_1=S_1$, if
$R_1+\cdots +R_i\not\equiv S_1+\cdots +S_i\pmod 2$, for some
$1\leq i\leq n$, then $i\geq 2$ (in particular, $n\geq 2$)
and $R_2+\cdots +R_i\not\equiv S_2+\cdots +S_i\pmod 2$. Then, by
7.3.(ii$^*$), written in the equivalent form from Remark 7.4(2),
$\N (-a_{1,i+1}b_{1,i-1})=\N (-a^*_{2,i+1}b_{2,i-1})\subseteq
\theta (X(M^*/N^*))\subseteq\theta (X(M/N))$. Hence 7.3(ii)
holds. \qed

\bpr Theorem 7.5 holds for the pair $M,N$ if it
holds for the pair $M^*,N^*$ or $n=1$.
\epr
\pf  If 7.5(4) fails, then Theorem 7.5 states that
$\theta (X(M/N))\not\subseteq\ooo\ffs$, which is true: Otherwise
$\theta (O^+(M)),\theta (O^+(N))\subseteq\theta
(X(M/N))\subseteq\ooo\ffs$, so 7.5(4) holds.

So we may assume that 7.5(4) holds. In particular,
$\theta (O^+(M))\subseteq\ooo\ffs$. By Lemma 7.20(4), $T_1=R_2-R_1$
so, by Theorem 1.9, 
$\bar G(a_{1,2},T_1)=\bar G(a_{1,2},R_2-R_1)=G(a_2/a_1)\subseteq\ooo\ffs$.
Hence 7.5(3) holds unconditionally at $i=1$.

Suppose that $n\geq 2$. Since the pair $M,N$ satisfies condition R1,
by Lemma 7.20(i), so does $M^*,N^*$. We also have
$\theta (O^+(M^*))\subseteq\theta (O^+(M))$ and
$\theta (O^+(N^*))\subseteq\theta (O^+(N))$, so 7.5(4) implies that
7.5(4$^*$) holds, as well. 

By Lemma 6.3(i),
$\theta (X(M/N))=\theta (O^+(M))\theta (X(M^*/N^*))$ and, by 7.5(iv),
$\theta (O^+(M))\subseteq\ooo\ffs$. Hence
$\theta (X(M/N))\subseteq\ooo\ffs$ is equivalent to
$\theta (X(M^*/N^*))\subseteq\ooo\ffs$. Hence, in order to prove our
statement, we must show that conditions 7.5(1)-(4) are equivalent to
7.5(1$^*$)-(4$^*$). 

We already have that 7.5(4) and 7.5(4$^*$) hold. By Lemma 7.20(iv), we
also have the equivalence between 7.5(i) and (ii) and 7.5(i$^*$) and
7.5(ii$^*$). And, as seen above, 7.5(iii) holds unconditionally at
$i=1$. So we are left to prove that 7.5(iii) at indices $\geq 2$ is
equivalent to 7.5(iii*). To do this, we show that for every
$2\leq i\leq m-2$  condition 7.5(iii) at index $i$ is equivalent to
7.5(iii$^*$) at index $i-1$. Indeed, 7.5(iii) at index $i$ states that
either $d(-a_{1,i+1}b_{1,i-1})=2e$ or $i\leq n$ and
$\bar G(a_{1,i+1}b_{1,i-1},T_i)\subseteq\ooo\ffs$, while 7.5(iii$^*$) at
index $i-1$ states that either $d(-a_{2,i+1}^*b_{2,i-1})=2e$ or
$i\leq n$ and $\bar G(a_{2,i+1}^*b_{2,i-1},T_i^*)\subseteq\ooo\ffs$. But
$a_{2,i+1}^*b_{2,i-1}=a_{1,i+1}b_{1,i-1}$ in $\fff/\ffs$ and, by Lemma
7.20(iv), if $i\leq n$, then $T_i^*=T_i$. Hence the two statements are
equivalent.

If $m=n=1$, then
$\theta (O^+(M))=\theta (O^+(N))=\theta (X(M/N))=\ffs$. Since
$\ffs\subset\ooo\ffs$, we have $\theta (X(M/N))\subset\ooo\ffs$ and
7.5(4) holds. To complete the proof, we note that 7.5(1)-(3) hold
trivially. (For 7.5(1) we have $R_1=S_1$, so $R_1\equiv S_1\pmod 2$ and
(2) and (3) are vacuous.)

Suppose now that $n=1$. Then, since R1 and 7.5(4) hold, by Lemma
7.20(iv), 7.5(1) and (2) hold uncondtionally. So we must prove that
$\theta (X(M/N))\subseteq\ooo\ffs$ is equivalent to 7.5(3). If
$m\leq 2$, by Lemma 7.21, we have
$\theta (X(M/N))=\theta (O^+(M))\subseteq\ooo\ffs$
unconditionally. But 7.5(3), too, holds unconditionally: If $m=1$, it
is vacuous and if $n=2$, it involves only the index $i=1$, where it
was already proven. If $m=3$, by Lemma 7.21,
$\theta (X(M/N))=\N (-a_{1,3}b_1)\theta (O^+(M))$. But
$\theta (O^+(M))\subseteq\ooo\ffs$, so
$\theta (X(M/N))\subseteq\ooo\ffs$ iff
$\N (-a_{1,3}b_1)\subseteq\ooo\ffs$, which is equivalent to
$d(-a_{1,3}b_1)=2e$. On the other hand, 7.5(3) holds unconditionally at
$i=1$ and at $i=2$ it states that $d(-a_{1,3}b_1)=2e$. Hence again we
have equivalence. \qed

{\bf Comparisson with the results from [B1].} Condition R1 is Theorem
II.1 from [B1] Lemma. Therefore all references to R1 from Theorems 7.3
and 7.5 are missing from their counterparts, Theorems II.2 and II.3
from [B1].

In [B1] we have $a_i=\pi^{R_i}\varepsilon_i$ and $b_i=\pi^{S_i}\eta_i$
and we put $\xi_i:=\varepsilon_{1,i+1}\eta_{1,i-1}$ for
$1\leq i\leq m-1$. Hence
$a_{1,i+1}b_{1,i-1}=\pi^{R_1+\cdots +R_{i+1}+S_1+\cdots +S_{i-1}}\xi_i$. 

If $M,N$ satisfy the condition (i) of Theorem 7.5, then both
$R_1+\cdots +R_{i+1}+S_1+\cdots +S_{i-1}$ and $T_i$ are even for all
$i$ so in $\fff/\ffs$ we have $a_{1,i+1}b_{1,i-1}=\pi^{T_i}\xi_i=\xi_i$.
Thus Theorem 7.5 and [B1, Theorem II.3] are the same, except for
condition R1 form the hypothesis of Theorem 7.5.

The general formula for $\theta (X(M/N))$ from [B1, Theorem II.2] is
$\theta (X(M/N))=\ups\beta G_1\cdots G_{m-1}\theta (O^+(M))$, where
$G_i=G(\pi^{R_{i+1}-S_i}\xi_i)$ if $i\leq n$ and
$R_1+\cdots +R_i\equiv S_1+\cdots +S_i\pmod 2$ and
$G_i=\N (-a_{1,i+1}b_{1,i-1})$ otherwise. But if
$R_1+\cdots +R_i\equiv S_1+\cdots +S_i\pmod 2$, then
$R_1+\cdots +R_{i+1}+S_1+\cdots +S_{i-1}\equiv R_{i+1}-S_i\pmod 2$, so
$a_{1,i+1}b_{1,i-1}=\pi^{R_{i+1}-S_i}\xi_i$ in $\fff/\ffs$ and so
$G(\pi^{R_{i+1}-S_i}\xi_i)=\bar G(a_{1,i+1}b_{1,i-1},R_{i+1}-S_i)$. So
[B1, Theorem II.2(ii)] follows from Theorem 7.3(i) and (ii). By
[B1, Theorem II.2(i)], in certain conditions, namely if
$\ups\beta G(a_{i+1}/a_i)\neq\fff$ for all $1\leq i\leq m-1$, the
factor $\theta (O^+(M))$ can be dropped from the formula for
$\theta (X(M/N))$.

In the next section we will show that, assuming the representation
theorem [B3, Theorem 4.5] the factor $\theta (O^+(M))$ can always be
dropped from the formula for $\theta (X(M/N))$. We will also prove an
analogue of Remark 1.8, namely that the factor $\ups\beta$ can be
dropped if $\ups\beta =\ups{[(R_{i+2}-S_i)/2]}$ for some
$1\leq i\leq m-2$ such that $R_{i+2}-S_i$ is odd.

\section{Appendix}

We now prove a stronger version of Theorem 1.7, assuming the
unpublished results from [B5]. In particular, we assume R1, so the
consequences of R1 from \S7 hold.

Same as in the previous section,
$M\cong{\prec a_1,\ldots,a_m\succ}$ and
$N\cong{\prec b_1,\ldots,b_n\succ}$ relative to good BONGs,
$R_i(M)=R_i$ and $R_i(N)=S_i$, $m-n\leq 2$ and if $m-n=2$, we make
the convention $S_{m-1}=S_{n+1}\gg 0$. We also put
$\alpha_i(M)=\alpha_i$ and $\alpha_i(N)=\beta_i$.

For every $\varepsilon\in\fff/\ffs$ and every $0\leq i\leq m$,
$0\leq j\leq n$ we define
$$d[\varepsilon a_{1,i}b_{1,j}]
:=\min\{ d(\varepsilon a_{1,i}b_{1,j},\alpha_i,\beta_j\}.$$
(If $i=0$ or $m$, we ignore $\alpha_i$. If $j=0$ or
$n$, we ignore $\beta_j$.) When $M=N$ we have similar expressions. If
$1\leq i-1\leq j\leq m$, then $d[\varepsilon a_{i,j}]:=
d[\varepsilon a_{1,i-1}a_{1,j}]
=\min\{ d(\varepsilon a_{i,j},\alpha_{i-1},\alpha_j\}$. (If $i=1$ or
$m+1$, we ignore $\alpha_i$. If $i=0$ or $m$, we ignore $\alpha_j$.)

The expressions $d[\varepsilon a_{1,i}b_{1,j}]$ and
$d[\varepsilon a_{i,j}]$ are invariants in the sense that they are
independent of the choice of the BONGs. They satisfy the same
domination principles as those satisfied by
$d(\varepsilon a_{1,i}b_{1,j})$ and $d(\varepsilon a_{i,j})$. If
$K\cong{\prec c_1,\ldots,c_k\succ}$ is a third lattice, then
$$d[\varepsilon\varepsilon'a_{1,i}c_{1,l}]\geq\min
\{ d[\varepsilon a_{1,i}b_{1,j}],d[\varepsilon b_{1,j}c_{1,l}]\}.$$
When two or more of the lattices $M,N,K$ coincide we get similar
inequalities, such as those inviolving the triplets
$d[\varepsilon\varepsilon'a_{1,i}b_{1,l}]$,
$d[\varepsilon a_{1,i}b_{1,j}]$, $d[\varepsilon b_{j+1,l}]$ or
$d[\varepsilon\varepsilon'a_{i,l}]$, $d[\varepsilon a_{i,j}]$,
$d[\varepsilon a_{j+1,l}]$.

We also have the relations
$$\alpha_i=\min\{ (R_{i+1}-R_i)/2+e,\, R_{i+1}-R_i+d[-a_{i,i+1}]\}.$$

\btm (i) If $N\subseteq M$ and the pair $M,N$ has property A, then
$$\theta (X(M/N))=(1+\p^\beta )\ffs\prod_{i=1}^{m-2}
\bar G(a_{1,i+1}b_{1,i-1},R_{i+1}-S_i).$$

(iii) If the pair $M,N$ doesn't have property A, then
$\theta (X(M/N))\supseteq\ooo\ffs$.

Same as for Theorem 7.3, if $m-n=2$, since $S_{m-1}\gg 0$, the last
factor in (i), $\bar G(a_{1,m}b_{1,m-2},R_m-S_{m-1})$, is equal to
$\N (-a_{1,m}b_{1,m-2})$. 

And if $m\leq 2$, then $\beta :=\infty$, so the
factor $\ups\beta =\ffs$ may be ignored in (i).
\etm

Note that we are under the assumption that R1 always hold, so, unlike
in Theorem 7.3(i), this condition is missing from the hypothesis of
Theorem 8.1(i). Also Theorem 8.1(ii) is the same as Theorem 7.3(iii).

Using the notation from Remark 7.4(1), Theorem 7.3(i) writes as
$\theta (X(M/N))=\tilde G(M/N)$, while Theorem 8.1(i) writes as
$\theta (X(M/N))=G(M/N)$. So Theorem 8.1(i) follows from Theorem
7.3(i), provided that $G(M/N)=\tilde G(M/N)$. But
$\tilde G(M/N)=G(M/N)\prod_{i=2}^{m-1}\bar G(a_{i,i+1},R_{i+1}-R_i)$,
so all we need is the following result.

\blm For every $2\leq i\leq m-1$ we have
$\bar G(a_{i,i+1},R_{i+1}-R_i)\subseteq G(M/N)$
\elm

For the rest of the Appendix, we assume that $M,N$ satisfy the
hypothesis of Theorem 8.1, i.e. $N\subseteq M$, $m-n\leq 2$ and the
pair $M,N$ has property $A$, and that

We first prove some preliminary results.

\blm (i) $R_i\leq S_i$ for every $1\leq i\leq n$.

(ii) $d[a_{1,i}b_{1,i}]\geq\bar A_i=\min\{ (R_{i+1}-S_i)/2+e,
R_{i+1}-S_i+d[-a_{1,i+1}b_{1,i-1}]\}$ for
$1\leq i\leq\min\{ m-1,n\}$. Explicitely,
$$d[a_{1,i}b_{1,i}]\geq\bar A_i=\begin{cases}
(R_{i+1}-S_i)/2+e&\text{if }d[-a_{1,i+1}b_{1,i-1}]>e+(S_i-R_{i+1})/2\\
R_{i+1}-S_i+d[-a_{1,i+1}b_{1,i-1}]&\text{if }d[-a_{1,i+1}b_{1,i-1}]\leq
e+(S_i-R_{i+1})/2\end{cases}.$$
In particular, if $R_{i+1}-S_i>2e$, then $a_{1,i}b_{1,i}=1$ in
$\fff\ffs$. 

(iii) For every $1\leq i\leq\min\{ n+1,m-1\}$ such that
$d[-a_{1,i}b_{1,i-2}]+d[-a_{1,i+1}b_{1,i-1}]>2e+S_{i-1}-R_{i+1}$ we
have $[b_1,\ldots,b_{i-1}]\rep [a_1,\ldots,a_i]$.
\elm
\pf (i) This was proved in Corollary 7.2, as a consequence of R1.

(ii) We have $R_{i+1}>S_{i-1}$ and $R_{i+2}>S_i$, so
$R_{i+1}+R_{i+2}>S_{i-1}+S_i$. Then, by [B5, Definition 6], $\bar A_i$
is defined and, by the following remark [B5, Remark 2.8], since we
have both $R_{i+1}\geq S_{i-1}$ and $R_{i+2}\geq S_i$, the terms
$2R_{i+1}-S_{i-1}-S_i+\alpha_{i+1}$ and
$R_{i+1}+R_{i+2}-2S_i+\beta_{i-1}$ can be ignored from the definition
of $\bar A_i$. Thus
$\bar A_i=\min\{ (R_{i+1}-S_i)/2+e,R_{i+1}-S_i+d[-a_{1,i+1}b_{1,i-1}]\}$.
Then, by [B5, Lemma 2.9], the condition [B2, Theorem 2.1(ii)],
$d[a_{1,i},b_{1,i}]\geq A_i$, is equivalent to
$d[a_{1,i},b_{1,i}]\geq\bar A_i$. For the explicit formulas, we note
that $(R_{i+1}-S_i)/2+e>$ or $\leq R_{i+1}-S_i+d[-a_{1,i+1}b_{1,i-1}]$
iff $d[-a_{1,i+1}b_{1,i-1}]>$ or $\leq e+(S_i-R_{i+1})/2$,
respectively.

If $R_{i+1}-S_i>2e$, then
$d[-a_{1,i+1}b_{1,i-1}]\geq 0>e+(S_i-R_{i+1})/2$, so
$d(a_{1,i}b_{1,i})\geq d[a_{1,i}b_{1,i}]\geq (R_{i+1}-S_i)/2+e>2e$,
which implies that $a_{1,i}b_{1,i}=1$ in $\fff/\ffs$.

(iii) We have $R_{i+1}>S_{i-1}$, by the property A, and
$d[-a_{1,i}b_{1,i-2}]+d[-a_{1,i+1}b_{1,i-1}]>2e+S_{i-1}-R_{i+1}$, so
$[b_1,\ldots,b_{i-1}]\rep [a_1,\ldots,a_i]$, by the condition
[B5, Theorem 2.1(iii')], which replaces [B5, Theorem 2.1(iii)]. (See
the comments following [B5, Theorem 2.1] and [B5, Lemma 2.16].) \qed

\blm If $G(M/N)\neq\fff$ and for some $1\leq i\leq m-1$ we have
$d(-a_{1,i+1}b_{1,i-1})\leq e+(S_i-R_{i+1})/2$ or $i=m-1=n+1$, then
$G(M/N)=\N (c)$, where $c=-a_{1,i+1}b_{1,i-1}$, and for every
$1\leq j\leq m-2$ we have $R_{j+2}-S_j+2d(c)>4e$.
\elm
\pf Since $d(-a_{1,i+1}b_{1,i-1})\leq e-(R_{i+1}-S_i)/2$, by Corollary
2.3(iii), we have
$\bar G(a_{1,i+1}b_{1,i-1},R_{i+1}-S_i)=\N (-a_{1,i+1}b_{1,i-1})=\N (c)$. 
Same happens by default if $i=n+1=m-1$. It follows that
$\N (c)\subseteq G(M/N)$. Since $[\fff :\N (c)]=2$ and
$G(M/N)\neq\fff$, we must have $G(M/N)=\N (c)$.

Since $(1+\p^\beta )\ffs\subseteq G(M/N)=\N (c)$, by Lemma ,
we have $\beta +d(c)>2e$. Since for every $j$ we have
$(R_{j+2}-S_j)/2\geq\beta$, we get $(R_{j+2}-S_)/2+d(c)>2e$,
i.e. $R_{j+2}-S_j+2d(c)>4e$. \qed

We denote by $I$ the set of all $1\leq i\leq n-1$ such that
$d[-a_{1,i+1}b_{1,i-1}]\leq e+(S_i-R_{i+1})/2$ and
$d[-a_{1,i+1}b_{1,i-1}]<d(-a_{1,i+1}b_{1,i-1})$. Since
$d[-a_{1,i+1}b_{1,i-1}]=\min\{ d(-a_{1,i+1}b_{1,i-1}),
\alpha_{i+1},\beta_{i-1}\}$, we have $I=I_\alpha\cup I_\beta$, where
$I_\alpha=\{ i\in I\,\mid\, d[-a_{1,i+1}b_{1,i-1}]=\alpha_{i+1}\}$ and
$I_\beta=\{ i\in I\,\mid\, d[-a_{1,i+1}b_{1,i-1}]=\beta_{i-1}\}$.

\blm Suppose that $G(M/N)\neq\fff$.

(i) If $i\in I_\alpha$, then
$d[-a_{1,i+2}b_{1,i}]=R_{i+1}-R_{i+2}+d[-a_{1,i+1}b_{1,i-1}]$ and
$i+1\in I$.

(ii) If $i\in I_\beta$, then 
$d[-a_{1,i}b_{1,i-2}]=S_{i-1}-S_i+d[-a_{1,i+1}b_{1,i-1}]$ and
$i-1\in I$. 
\elm
\pf First note that, both for (i) and (ii), we have $i\in I$, so
$d[-a_{1,i+1}b_{1,i-1}]\leq e+(S_i-R_{i+1})/2$. By Lemma
8.3(ii), we have
$d[a_{1,i}b_{1,i}]\geq R_{i+1}-S_i+d[-a_{1,i+1}b_{1,i-1}]$. 

(i) Since $R_{i+2}>S_i$, we have
$\alpha_{i+1}=d[-a_{1,i+1}b_{1,i-1}]\leq
e+(S_i-R_{i+1})/2<(R_{i+2}-R_{i+1})/2+e$. Since $\alpha_{i+1}
=\min\{ (R_{i+2}-R_{i+1})/2+e,R_{i+2}-R_{i+1}+d[-a_{i+1,i+2}]\}$, we have
$\alpha_{i+1}=R_{i+2}-R_{i+1}+d[-a_{i+1,i+2}]$, so
$d[-a_{i+1,i+2}]=R_{i+1}-R_{i+2}+\alpha_{i+1}
=R_{i+1}-R_{i+2}+d[-a_{1,i+1}b_{1,i-1}]$. On the other hand,
$d[a_{1,i}b_{1,i}]\geq R_{i+1}-S_i+d[-a_{1,i+1}b_{1,i-1}]$. Since
$R_{i+2}>S_i$, this implies
$d[a_{1,i}b_{1,i}]>R_{i+1}-R_{i+2}+d[-a_{1,i+1}b_{1,i-1}]=d[-a_{i+1,i+2}]$. By
the domination principle, it follows that
$d[-a_{1,i+2}b_{1,i}]=d[-a_{i+1,i+2}]=R_{i+1}-R_{i+2}+d[-a_{1,i+1}b_{1,i-1}]$.

Since $d[-a_{1,i+1}b_{1,i-1}]\leq e+(S_i-R_{i+1})/2$, we have
$d[-a_{1,i+2}b_{1,i}]\leq R_{i+1}-R_{i+2}+e+(S_i-R_{i+1})/2
=(R_{i+1}+S_i)/2-R_{i+2}+e$. But $R_{i+2}>S_i$ and, by Lemma 8.3(i),
$R_{i+1}\leq S_{i+1}$. It follows that
$(R_{i+1}+S_i)/2-R_{i+2}+e<(R_{i+1}-R_{i+2})/2+e\leq
e+(S_{i+1}-R_{i+2})/2$. Hence
$d[-a_{1,i+2}b_{1,i}]<e+(S_{i+1}-R_{i+2})/2$. To complete the proof of
$i+1\in I$, we must show that
$d[-a_{1,i+2}b_{1,i}]<d(-a_{1,i+2}b_{1,i})$. Suppose the contrary,
i.e. that $d[-a_{1,i+2}b_{1,i}]=d(-a_{1,i+2}b_{1,i})$. Then
$d(-a_{1,i+2}b_{1,i})<e+(S_{i+1}-R_{i+2})/2$, so if
$c=-a_{1,i+2}b_{1,i}$, then, by Lemma 8.4 we have
$R_{i+2}-S_i+2d(c)>4e$. But
$d(c)=d[-a_{1,i+2}b_{1,i}]\leq (R_{i+1}+S_i)/2-R_{i+2}+e$, so
$4e<R_{i+2}-S_i+2d(c)\leq
R_{i+2}-S_i+2((R_{i+1}+S_i)/2-R_{i+2}+e)=R_{i+1}-R_{i+2}+2e$, which
contradicts $R_{i+2}-R_{i+1}\geq -2e$. Hence $i+1\in I$.

(ii) Since $R_{i+1}>S_{i-1}$, we have
$\beta_{i-1}=d[-a_{1,i+1}b_{1,i-1}]\leq
e+(S_i-R_{i+1})/2<(S_i-S_{i-1})/2+e$. Since
$\beta_{i-1}=\min\{ (S_i-S_{i-1})/2+e,S_i-S_{i-1}+d[-b_{i-1,i}]\}$, we
have $\beta_{i-1}=S_i-S_{i-1}+d[-b_{i-1,i}]$, so
$d[-b_{i-1,i}]=S_{i-1}-S_i+\beta_{i-1}
=S_{i-1}-S_i+d[-a_{1,i+1}b_{1,i-1}]$. On the other hand,
$d[a_{1,i}b_{1,i}]\geq R_{i+1}-S_i+d[-a_{1,i+1}b_{1,i-1}]$. Since
$R_{i+1}>S_{i-1}$, we have
$d[a_{1,i}b_{1,i}]>S_{i-1}-S_i+d[-a_{1,i+1}b_{1,i-1}]=d[-b_{i-1,i}]$. By
the domination principle, it follows that
$d[-a_{1,i}b_{1,i-2}]=d[-b_{i-1,i}]=S_{i-1}-S_i+d[-a_{1,i+1}b_{1,i-1}]$.

Since $d[-a_{1,i+1}b_{1,i-1}]\leq e+(S_i-R_{i+1})/2$, we have
$d[-a_{1,i}b_{1,i-2}]\leq S_{i-1}-S_i+e+(S_i-R_{i+1})/2
=S_{i-1}-(R_{i+1}+S_i)/2+e$. But $R_{i+1}>S_{i-1}$ and, by Lemma 8.3(i),
$R_i\leq S_i$. It follows that
$S_{i-1}-(R_{i+1}+S_i)/2+e<(S_{i-1}-S_i)/2+e\leq
e+(S_{i-1}-R_i)/2$. Hence
$d[-a_{1,i+2}b_{1,i}]<e+(S_{i-1}-R_i)/2$. To complete the proof of
$i-1\in I$, we must show that
$d[-a_{1,i}b_{1,i-2}]<d(-a_{1,i}b_{1,i-2})$. Suppose the contrary,
i.e. that $d[-a_{1,i}b_{1,i-2}]=d(-a_{1,i}b_{1,i-2})$. Then
$d(-a_{1,i}b_{1,i-2})<e+(S_{i-1}-R_i)/2$, so if
$c=-a_{1,i}b_{1,i-2}$, then, by Lemma 8.4, we have
$R_{i+1}-S_{i-1}+2d(c)>4e$. But
$d(c)=d[-a_{1,i}b_{1,i-2}]\leq S_{i-1}-(R_{i+1}+S_i)/2+e$, so
$4e<R_{i+1}-S_{i-1}+2d(c)\leq
R_{i+1}-S_{i-1}+2(S_{i-1}-(R_{i+1}+S_i)/2+e)=S_{i-1}-S_i+2e$, which
contradicts $S_i-S_{i-1}\geq -2e$. Hence $i-1\in I$. \qed

\blm If $G(M/N)\neq\fff$ then $I=\emptyset$. 
\elm
\pf Suppose that $I=I_\alpha\cup I_\beta\neq\emptyset$. First note
that $\min I\in I_\alpha$, so $I_\alpha\neq\emptyset$. Indeed, if
$i:=\min I\in I_\beta$, then, by Lemma 8.5(ii), $i-1\in I$, which
contradicts the minimality of $i$.

Let now $i=\max I_\alpha$. Then, by Lemma 8.5(i), $i+1\in I$. We
cannot have $i+1\in I_\alpha$ because of the maximality of $i$. Thus
$i+1\in I_\beta$. Since $i\in I_\alpha$ and $i+1\in I_\beta$, by Lemma
8.5(i)  and (ii), we have
$d[-a_{1,i+2}b_{1,i}]=R_{i+1}-R_{i+2}+d[a_{1,i+1}b_{1,i-1}]$ and
$d[-a_{1,i+1}b_{1,i-1}]=S_i-S_{i+1}+d[a_{1,i+2}b_{1,i}]$. By adding
these relations, we get $0=R_{i+1}-R_{i+2}+S_i-S_{i+1}$. But
$R_{i+2}>S_i$ and, by Lemma 8.3(i), $R_{i+1}\leq S_{i+1}$. Thus
$R_{i+1}-R_{i+2}+S_i-S_{i+1}<0$. Contradiction. Hence
$I=\emptyset$. \qed

\bco If $G(M/N)\neq\fff$ and
$d[-a_{1,i+1}b_{1,i-1}]\leq e+(S_i-R_{i+1})/2$, then
$d[-a_{1,i+1}b_{1,i-1}]=d(-a_{1,i+1}b_{1,i-1})$. Also if
$c=-a_{1,i+1}b_{1,i-1}$, then $G(M/N)=\N (c)$, $R_{j+2}-S_j>4e-2d(c)$
for all $1\leq j\leq m-2$ and 
$d[a_{1,i}b_{1,i}]\geq R_{i+1}-S_i+d(c)$. 
\eco
\pf If $d[-a_{1,i+1}b_{1,i-1}]<d(-a_{1,i+1}b_{1,i-1})$, then, by
definition, $i\in I$. But, by Lemma 8.6 $I=\emptyset$. So we must have 
$d[-a_{1,i+1}b_{1,i-1}]=d(-a_{1,i+1}b_{1,i-1})$. This implies
$d(-a_{1,i+1}b_{1,i-1})\leq e+(S_i-R_{i+1})/2$ so, by Lemma 8.4,
$G(M/N)=\N (c)$ and $R_{j+2}-S_j+2d(c)>4e$
for all $1\leq j\leq m-2$. And, by Lemma 8.3(ii),
$d[a_{1,i}b_{1,i}\geq
R_{i+1}-S_i+d[-a_{1,i+1}b_{1,i-1}]=R_{i+1}-S_i+d(c)$. \qed

\blm Let $a_1,\ldots,a_{i+1},b_1,\ldots,b_{i-1}\in\fff$. If two of the
following statements hold, then so does the third:
(a) $[b_1,\ldots,b_{i-1}]\rep [a_1,\ldots,a_i]$, (b)
$[b_1,\ldots,b_i]\rep [a_1,\ldots,a_{i+1}]$,\\
(c) $(a_{1,i},b_{1,i},-a_{1,i+1}b_{1,i-1})_\p =1$.
\elm
\pf Let $V_1=[a_1,\ldots,a_{1,i+1}]$,
$V_2=[b_1,\ldots,b_{i-1},a_{1,i}b_{1,i-1},a_{i+1}]$ and
$V_3=[b_1,\ldots,b_i,a_{1,i+1}b_{1,i}]$. We claim that statements (a),
(b) and (c) are equivalent to $V_1\cong V_2$, $V_1\cong V_3$ and
$V_2\cong V_3$. Then obviously any two of these isometries implies the
thir one, so we are done.

By [OM, Theorem 63:21], (a) is equivalent to
$[a_1,\ldots,a_i]\cong [b_1,\ldots,b_{i-1},a_{1,i}b_{1,i-1}]$, which
in turn is equivalent to
$[a_1,\ldots,a_i,a_{i+1}]\cong [b_1,\ldots,b_{i-1},a_{1,i}b_{1,i-1},a_{i+1}]$,
i.e. $V_1\cong V_2$. And (b) is equivalent to
$[a_1,\ldots,a_{i+1}]\cong [b_1,\ldots,b_i,a_{1,i+1}b_{1,i}]$,
i.e. $V_1\cong V_3$.

By Witt cancellation, $V_2\cong V_3$ is equivalent to
$[a_{1,i}b_{1,i-1},a_{i+1}]\cong [b_i,a_{1,i+1}b_{1,i}]$. Now for
every $a,b,c,d\in\fff/\ffs$ with $ab=cd$ we have $[a,b]\cong [c,d]$
iff $c\rep [a,b]$, which is equivalent to $(ac,-ab)_\p =1$. When
$a,b,c,d$ are $a_{1,i}b_{1,i-1},a_{i+1},b_i,a_{1,i+1}b_{1,i}$ we get that
$V_2\cong V_3$ is equivalent to
$(a_{i,i}b_{1,i},-a_{i,i+1}b_{1,i-1})_\p =1$, i.e. to (c). this
concludes the proof.\qed

We now have all ingredients for the proof of Lemma 8.2, which implies
Theorem 8.1.

{\bf Proof of Lemma 8.2.} We may assume that
$G(M/N)\neq\fff$. Otherwise the statement is trivial. By Lemma 8.3(i),
we have $R_i\leq S_i$, so $R_{i+1}-R_i\geq R_{i+1}-S_i$. (This happens
also when $i=n+1=m-1$ in which case $S_{n+1}\gg 0$, which implies
$R_{n+1}\leq S_{n+1}$.)

First assume that $R_{i+1}-R_i=-2e$. We have $R_{i+1}>S_{i-1}$ so
$R_i-S_{i-1}>R_i-R_{i+1}=2e$. Then, by Lemma 8.3(ii), in $\fff/\ffs$
we have $a_{1,i-1}b_{1,i-1}=1$, so $a_{1,i+1}b_{1,i-1}=a_{i,i+1}$. By
Lemma 2.6(i), since $R_{i+1}-R_i\geq  R_{i+1}-S_i$, this implies
$\bar G(a_{i,i+1},R_{i+1}-R_i)\subseteq
\bar G(a_{1,i+1}b_{1,i-1},R_{i+1}-S_i)\subseteq G(M/N)$.

Suppose now that $R_{i+1}-R_i\neq -2e$. By Lemma 2.16, we have 
\begin{multline*}
\bar G(a_{i,i+1},R_{i+1}-R_i)\bar G(a_{1,i+1}b_{1,i-1},R_{i+1}-S_i)\\
=\langle a_{1,i-1}b_{1,i-1}\rangle\ups{f(R_{i+1}-R_i)+d(a_{1,i-1}b_{1,i-1})}
\bar G(a_{1,i+1}b_{1,i-1},R_{i+1}-S_i).
\end{multline*}
Since $\bar G(a_{1,i+1}b_{1,i-1},R_{i+1}-S_i)\subseteq G(M/N)$, the
inclusion $\bar G(a_{i,i+1},R_{i+1}-R_i)\subseteq G(M/N)$ is
equivalent to $\langle a_{1,i-1}b_{1,i-1}\rangle
\ups{f(R_{i+1}-R_i)+d(a_{1,i-1}b_{1,i-1})}\subseteq G(M/N)$, i.e. to 
$\ups{f(R_{i+1}-R_i)+d(a_{1,i-1}b_{1,i-1})}\subseteq G(M/N)$ and
$a_{1,i-1}b_{1,i-1}\in G(M/N)$. 

We have two subcases:

(a) $d[-a_{1,i}b_{1,i-2}]>e+(S_{i-1}-R_i)/2$. By Lemma 8.3(ii), we
have $d(a_{1,i-1}b_{1,i-1})\geq (R_i-S_{i-1})/2+e$. Together with
$f(R_{i+1}-R_i)\geq (R_{i+1}-R_i)/2-e$, this implies
$f(R_{i+1}-R_i)+d(a_{1,i-1}b_{1,i-1})\geq (R_{i+1}-S_{i-1})/2
\geq\beta$, so $(1+\p^{f(R_{i+1}-R_i)/2+d(a_{1,i-1}b_{1,i-1})})\ffs
\subseteq (1+\p^\beta )\ffs\subseteq G(M/N)$.

Next we prove that
$a_{1,i-1}b_{1,i-1}\in\bar G(a_{1,i}b_{1,i-2},R_i-S_{i-1})\subseteq G(M/N)$,
so we are done.

If $(a_{1,i}b_{1,i-2},R_i-S_{i-1})=(-1,-2e)$ in $\fff/\ffs\times\RR$,
then $\bar G(a_{1,i}b_{1,i-2},R_i-S_{i-1})=\ooo\ffs$ and
$\ord a_{1,i}b_{1,i-2}$ is even (we have
$a_{1,i}b_{1,i-2}=-1$ in $\fff/\ffs$). It follows that
$\ord a_{1,i-1}b_{1,i-1}=\ord a_{1,i}b_{1,i-2}-R_i+S_{i-1}
=\ord a_{1,i}b_{1,i-2}+2e$ is even, so $a_{1,i-1}b_{1,i-1}\in\ooo\ffs
=\bar G(a_{1,i}b_{1,i-2},R_i-S_{i-1})$. 

If $(a,R)\neq (-1,-2e)$, then we use Corollary 2.13(iii). Since
$d(a_{1,i-1}b_{1,i-1})\geq (R_i-S_{i-2})/2+e$ and
$d(-a_{1,i}b_{1,i-2})\geq d[-a_{1,i}b_{1,i-2}]>e-(R_i-S_{i-1})/2$,
we have $a_{1,i-1}b_{1,i-1}\in\ups{(R_i-S_{i-1})/2+e}\subseteq
\bar G(a_{1,i}b_{1,i-2},R_i-S_{i-1})$.

(b) $d[-a_{1,i}b_{1,i-2}]\leq e+(S_{i-1}-R_i)/2$. By Corollary 8.7, if
$c=-a_{1,i}b_{1,i-2}$, then $d(c)=d[-a_{1,i}b_{1,i-2}]$,
$G(M/N)=\N (c)$, $d[a_{1,i-1}b_{1,i-1}]\geq R_i-S_{i-1}+d(c)$ and
$R_{j+2}-S_j>4e-2d(c)$ for $1\leq j\leq m-2$. Since
$f(R_{i+1}-R_i)\geq R_{i+1}-R_i-2e$,
$d(a_{1,i-1}b_{1,i-1})\geq R_i-S_{i-1}+d(c)$ and
$R_{i+1}-S_{i-1}>4e-2d(c)$, we have
$f(R_{i+1}-R_i)+d(a_{1,i-1}b_{1,i-1})\geq
R_{i+1}-S_{i-1}-2e+d(c)>2e-d(c)$. By Lemma 1.2(i) 
$\ups{f(R_{i+1}-R_i)+d(a_{1,i-1}b_{1,i-1})}\subseteq\N (c)=G(M/N)$.

We have $d[a_{1,i}b_{1,i-2}]\leq e+(S_{i-1}-R_i)/2<\infty$, so 
$d[a_{1,i-1}b_{1,i-1}]=d(c)\leq 2e$. Then, for $1\leq j\leq m-2$,
since $R_{j+2}-S_j>4e-2d(c)$, we have
$2e+S_j-R_{j+2}<2d(c)-2e\leq d(c)$. Since $R_{i+1}>S_{i-1}$ and
$d[-a_{1,i-2}b_{1,i}]+d[a_{1,i+1}b_{1,i-1}]\geq
d(c)+0>2e+S_{i-1}+R_{i+1}$, by Lemma 8.3(iii), we have
$[b_1,\ldots,b_{i-1}]\rep [a_1,\ldots,a_i]$. If $i=2$, this writes as
$b_1\rep [a_1,a_2]$, so $(a_1b_1,c)_\p =(a_1,b_1,-a_{1,2})_\p =1$.
Hence $a_1b_1\in\N (c)=G(M/N)$ and we are done. If $i\geq 3$, then
similarly $R_i>S_{i-2}$ and
$d[-a_{1,i-1}b_{1,i-3}]+d[a_{1,i}b_{1,i-2}]\geq 0+d(c)>2e+S_{i-2}+R_i$,
so $[b_1,\ldots,b_{i-2}]\rep [a_1,\ldots,a_{i-1}]$. Together with 
$[b_1,\ldots,b_{i-1}]\rep [a_1,\ldots,a_i]$, by Lemma 8.8, this
implies
$(a_{1,i-1}b_{1,i-1},c)_\p =(a_{1,i-1}b_{1,i-1},-a_{1,i}b_{1,i-2})_\p =1$. 
Hence $a_{1,i-1}b_{1,i-1}\in\N (c)=G(M/N)$. \qed

\blm If $R_{i+2}-S_i$ is odd for some $1\leq i\leq m-2$, then
$\ups{[(R_{i+2}-S_i)/2]}\subseteq
\bar G(a_{1,i+1}b_{1,i-1},R_{i+1}-S_i)
\bar G(a_{1,i+2}b_{1,i},R_{i+2}-S_{i+1})$.
\elm
\pf By Lemma 8.3(i), $R_{i+1}\leq S_{i+1}$. (This includes the case
$i=n=m-2$, since $S_{n+1}\gg 0$ implies $R_{n+1}\leq S_{n+1}$.)

For convenience, we denote $a=a_{1,i+1}b_{1,i-1}$ and
$b=a_{1,i+2}b_{1,i}$. We have $b=a_{i+2}b_ia$, so
$\ord b=R_{i+2}+S_i+\ord a$. Since $R_{i+2}-S_i$ is odd, $\ord a$ and
$\ord b$ have opposite parities. Thus $\ord ab$ is odd and so
$d(ab)=0$.

We put $G=\bar G(a,R_{i+1}-S_i,a)\bar G(b,R_{i+2}-S_i)$. Since
$[(R_{i+2}-S_i)/2]=(R_{i+2}-S_i-1)/2$, we must prove that
$\ups{(R_{i+2}-S_i-1)/2}\subseteq G$. If $G=\fff$, this holds
trivially, so we will assume that $G\neq\fff$. Note that
$a\in\bar G(a,R_{i+1}-S_i,a)\subseteq G$ and
$b\in\bar G(b,R_{i+2}-S_i)\subseteq G$, so $ab\in G$. But $\ord ab$ is
odd, so $G\not\subseteq\ooo\ffs$. In particular,
$G\neq\ooo\ffs$. Since also $G\neq\fff$, we have
$\ooo\ffs\not\subseteq G$.

We consider the four cases when $R_{i+1}-S_i$ and $R_{i+2}-S_{i+1}$ are
$\leq 2e$ or $\geq 2e+1$.

(1) $R_{i+1}-S_i,R_{i+2}-S_{i+1}\geq 2e+1$. Since
$R_{i+1}\leq S_{i+1}$, we have
$R_{i+2}-S_i\geq (R_{i+2}-S_{i+1})+(R_{i+1}-S_i)>4e+2$. It follows
that $(R_{i+2}-S_i-1)/2>2e$, so
$\ups{(R_{i+2}-S_i-1)/2}=\ffs\subseteq G$ trivially.

(2) $R_{i+1}-S_i,R_{i+2}-S_{i+1}\leq 2e$. By Lemmas 2.12(i) and 2.16,
we have $G\supseteq\bar G(a,2e)\bar G(b,2e)\supseteq
\ups{2e/2+d(ab)-e}=\ups 0=\fff$, i.e. $G=\fff$, contradicting our
assumption. 

(3) $R_{i+1}-S_i\leq 2e$ and $R_{i+2}-S_{i+1}\geq 2e+1$. We have two
subcases:

(a) $d(-a)>e-(R_{i+1}-S_i)/2$. We cannot have $(a,R_{i+1}-S_i)=-2e$,
since this implies $\ooo\ffs =\bar G(a,R_{i+1}-S_i)\subseteq G$. Then,
by Corollary 2.13(iii),
$\ups{(R_{i+1}-S_i)/2+e}\subseteq\bar G(a,R_{i+1}-S_i)\subseteq G$. On
the other hand, $R_{i+1}\leq S_{i+1}$, so
$R_{i+2}-R_{i+1}\geq R_{i+2}-S_{i+1}\geq 2e+1$. It follows that
$(R_{i+2}-S_i-1)/2\geq (R_{i+1}-S_i)/2+e$, and so
$\ups{(R_{i+2}-S_i-1)/2}\subseteq\ups{(R_{i+1}-S_i)/2+e}\subseteq G$.

(b) $d(-a)\leq e-(R_{i+1}-S_i)/2$. By Corollary 2.13(iii),
$\N (-a)=\bar G(a,R_{i+1}-S_i)\subseteq G$. 
But $G\neq\fff$, so $G=\N (-a)$. We cannot have $a=-\Delta$ in
$\fff/\ffs$, since this would imply $G=\N (\Delta )=\ooo\ffs$. Since
$R_{i+2}-S_{i+1}>2e$, by Corollary 2.15,
$\bar G(b,R_{i+2}-S_{i+1})\subseteq G=\N (-a)$ implies
$R_{i+2}-S_{i+1}+d(-a)=R_{i+2}-S_{i+1}+d(ab)+d(-a)>4e$. And, by Lemma
8.3(ii), in $\fff/\ffs$ we have $a_{1,i+1}b_{1,i+1}=1$, so
$a=a_{1,i+1}b_{1,i-1}=b_{i,i+1}$. Since $b_{i,i+1}=a\neq -\Delta$ in
$\fff/\ffs$, the inequality
$S_{i+1}-S_i+d(-a)=S_{i+1}-S_i+d(-b_{i,i+1})\geq 0$ from Lemma 1.4 is
strict, so $S_{i+1}-S_i-1+d(-a)>0$. We add this inequality to
$R_{i+2}-S_{i+1}+d(-a)>4e$ and divide by $2$ and we get
$(R_{i+2}-S_i-1)/2+d(-a)>2e$, which, by Lemma 1.2(i), implies
$\ups{(R_{i+2}-S_i-1)/2}\subseteq\N (-a)=G$.

(4) $R_{i+1}-S_i\geq 2e+1$ and $R_{i+2}-S_{i+1}\leq 2e$. Again, we
have two subcases: 

(a) $d(-b)>e-(R_{i+2}-S_{i+1})/2$. We cannnot have
$(b,R_{i+2}-S_i)=(-1,-2e)$ in $\fff/\ffs$, since this will
$\ooo\ffs =\bar G(b,R_{i+2}-S_{i+1})\subseteq G$. Then, by Corollary
2.13(iii), $\ups{(R_{i+2}-S_{i+1})/2+e}\subseteq\bar
G(b,R_{i+2}-S_{i+1})\subseteq G$. Since $R_{i+1}\leq S_{i+1}$, we have
$S_{i+1}-S_i\geq R_{i+2}-S_i\geq 2e+1$, so
$(R_{i+2}-S_i-1)/2\geq (R_{i+2}-S_{i+1})/2+e$, which implies
$\ups{(R_{i+2}-S_i-1)/2}\subseteq\ups{(R_{i+2}-S_{i+1})/2+e}
\subseteq G$.

(b) $d(-b)\leq e-(R_{i+2}-S_{i+1})/2$. By Corollary 2.13(iii),
$\N (-b)=\bar G(b,R_{i+2}-S_{i+1})\subseteq G$. But $G\neq\fff$, so
$G=\N (-b)$. We cannot have $a=-\Delta$ in $\fff/\ffs$, since this
would imply $G=\N (\Delta )=\ooo\ffs$. Since $R_{i+1}-S_i>2e$, by
Corollary 2.15, $\bar G(a,R_{i+1}-S_i)\subseteq G=\N (-b)$ implies
$R_{i+1}-S_i+d(-b)=R_{i+1}-S_i+d(ab)+d(-b)>4e$. And, by Lemma 8.3(ii),
in $\fff/\ffs$ we have $a_{1,i}b_{1,i}=1$, so
$b=a_{1,i+2}b_{1,i}=a_{i+1,i+2}$. Since $a_{i+1,i+2}=b\neq -\Delta$ in
$\fff/\ffs$, the inequality
$R_{i+2}-R_{i+1}+d(-b)=R_{i+2}-R_{i+1}+d(-b)\geq 0$ from Lemma 1.4 is
strict, so $R_{i+2}-R_{i+1}-1+d(-b)>0$. We add this inequality to
$R_{i+1}-S_i+d(-b)>4e$ and we divide by $2$ and we get
$(R_{i+2}-S_i-1)/2+d(-b)>2e$, which, by Lemma 1.2(i) implies
$\ups{(R_{i+2}-S_i-1)/2}\subseteq\N (-b)=G$. \qed

\bpr (i) If $R_1+\cdots +R_i\not\equiv S_1+\cdots +S_i\pmod 2$ for
some $1\leq i\leq\min\{ n,m-1\}$, then in (i) the factor
$\bar G(a_{1,i+1}b_{1,i-1},R_{i+1}-S_i)$ may be replaced by
$\N (-a_{1,i+1}b_{1,i-1})$.

(ii) If $\ups\beta =\ups{[(R_{i+2}-S_i)/2]}$ for some
$1\leq i\leq m-2$ such that $R_{i+2}-S_i$ is odd, then the factor
$\ups\beta$ can be dropped in the formula from Theorem 8.1(i).
\epr
\pf (i) This is the same as Theorem 7.3(ii). (By Remark 7.4(2), both
statements are equivalent to
$\N (-a_{1,i+1}b_{1,i-1})\subseteq\theta (X(M/N))$ for all $i$ with
$R_1+\cdots +R_i\not\equiv S_1+\cdots +S_i\pmod 2$.)

(ii) By Lemma 8.9, 
$\ups\beta =\ups{[(R_{i+2}-S_i)/2]}\subseteq
\bar G(a_{1,i+1}b_{1,i-1},R_{i+1}-S_{i-1})\bar G(a_{1,i+1}b_{1,i},R_{i+2}-S_i)
\subseteq\prod_{i=1}^{m-1}\bar
G(a_{1,i+1}b_{1,i-1},R_{i+1}-S_i)$. Thus the factor $\ups\beta$
is superfluous in Theorem 8.1(i). \qed

{\bf References}
\bigskip

[B1] C. N. Beli, {\it Integral spinor norm groups over dyadic local
fields and representations of quadratic lattices}, Ohio State
University thesis (2001).

[B2] C. N. Beli, {\it Integral spinor norm groups over dyadic local
fields}, J. Number Theory 102 (2003), 125–182.

[B3] C. N. Beli, {\it Representations of integral quadratic forms over
dyadic local fields}, Electron. Res. Announc. Amer. Math. Soc. 12
(2006), 100–112.

[B4] C. N. Beli, {\it A new approach to classification of integral
quadratic forms over dyadic local fields}, Trans. Amer. Math. Soc. 362
(2010), 1599–1617.

[B5] C. N. Beli, {\it Representations of quadratic lattices over
dyadic local fields}, (2019), arXiv:1905.04552v2.

[H1] J. S. Hsia, {\it Spinor norms of local integral rotations I},
Pacific J. Math., Vol. 57 (1975), 199 - 206.

[H2] J. S. Hsia, {\it Representations by spinor genera}, Pacific
J. Math., 63 (1976), 147 - 152.

[SP] R. Schulze-Pillot, Darstellung durch Spinorgeschlechter
tern\"arer quadratischer Formen, J. Number Theory, Vol. 12 (1980),
529 - 540.

[X1] F. Xu, {\it A remark on spinor norms of local integral rotations
I}, Pacific J. Math., vol. 136 (1989), 81 - 84.

[X2] F. Xu, {\it Integral spinor norms in dyadic local fields I},
Pacific J. Math., Vol. 157 (1993), 179 - 200.

[X3] F. Xu, {\it Generation of integral orthogonal groups over dyadic
local fields}, Pacific J. Math., vol. 167 (1995), 385 - 398.

[X4] F. Xu, {\it Arithmetic Springer theorem on quadratic forms under
field exten­sions of odd degree}, Contemporary Math., Vol. 249 (1999),
175 - 197.

\end{document}